\documentclass[11pt,twoside]{article}
\usepackage{amsmath,bm}
\usepackage{amssymb}
\usepackage{fancyhdr}
\usepackage{latexsym}
\usepackage{bbding}
\usepackage{mathrsfs}
\usepackage{exscale}
\usepackage{relsize}
\usepackage{wasysym}
\usepackage{cite}
\usepackage{multicol,graphics}
\usepackage{booktabs}
\usepackage{soul}  
\usepackage{color}
\usepackage{graphicx}    
\usepackage{caption}     

\allowdisplaybreaks[4]
\def\to{\rightarrow}
\def\e{\left}
\def\r{\right}
\def\l{\leq}
\def\g{\geq}
\def\fr{\frac}

\def\i{\infty}

\def\R{\mathbb{R}}

\def\p{\partial}

\def\dl{\delta}
\def\vp{\varepsilon}

\def\lb{\label}

\def\bar{\overline}

\newcommand{\SE}{\setcounter{equation}{0} \section}
\renewcommand\tilde{\widetilde}
\newcommand{\be}{\begin{equation}}
\newcommand{\ee}{\end{equation}}
\newcommand{\bc}{\begin{cases}}
\newcommand{\ec}{\end{cases}}
\newcommand{\bes}{\begin{equation*}}
\newcommand{\ees}{\end{equation*}}
\newcommand{\bls}{\begin{align*}}
\newcommand{\els}{\end{align*}}
\newcommand{\baa}{\begin{array}}
\newcommand{\eaa}{\end{array}}
\newcommand{\ba}{\begin{eqnarray}}
\newcommand{\ea}{\end{eqnarray}}
\newcommand{\bas}{\begin{eqnarray*}}
\newcommand{\eas}{\end{eqnarray*}}
\newcommand{\bd}{\begin{description}}
\newcommand{\ed}{\end{description}}

\newtheorem{theo}{\bf Theorem}[section]
\newtheorem{lem}[theo]{\bf Lemma}

\newtheorem{defi}[theo]{\bf Definition}
\newtheorem{rem}[theo]{\bf Remark}

\newenvironment{pr}[1][Proof]{\noindent\textbf{#1.} }{\hfill $\Box$}
\allowdisplaybreaks
\begin{document}
\date{}
\title{\bf{Combustion transition fronts in unbounded domains with multiple cylindrical branches}}
\author{Yang-Yang Yan$^{a}$,\  Wei-Jie Sheng$^{a,}$\thanks{Corresponding author
 (E-mail address: shengwj09@hit.edu.cn).}\ , \ Zhi-Cheng Wang$^{b}$
\\
\footnotesize{$^a$ School of Mathematics, Harbin Institute of Technology}, \\
\footnotesize{Harbin, Heilongjiang, 150001, People's Republic of China}\\
\footnotesize{$^b$ School of Mathematics and Statistics, Lanzhou University,} \\
\footnotesize{Lanzhou, Gansu, 730000, People's Republic of China}
}
\maketitle

\begin{abstract}
This paper is concerned with propagation dynamics for combustion reaction-diffusion equations in domains with multiple cylindrical branches. We first establish the existence and uniqueness of a time-increasing entire solution behaving like planar traveling fronts in some branches and converging to $0$ in the remaining part of the domain as $t\to-\infty$. Under the assumption of complete propagation, we then show that this entire solution propagates into the other branches in the form of planar traveling fronts (up to finite shifts) and converges to $1$ elsewhere as $t\to+\infty$. In particular, it is proved that this entire solution is a transition front connecting $0$ and $1$, whose global mean speed coincides with the planar wave speed. 
By assuming complete propagation for front-like solutions originating from single branch, we further prove that every transition front connecting $0$ and $1$ propagates completely. Moreover, we show that the global mean speed is independent of the choice of transition front. Namely, all transition fronts connecting $0$ and $1$ share the same global mean speed.
Finally, we provide two sufficient geometric conditions under which the complete propagation assumptions are satisfied.

\end{abstract}

\section{Introduction and main results}
\subsection{Setting of the problem}\label{ssss1}
Reaction-diffusion equations are fundamental mathematical models for describing propagation phenomena arising in flame propagation, population dynamics, chemical reactions and many other related fields, see \cite{M1,S2,vvv}. A central issue in the study of such equations is to understand transitions between different steady states. In the whole space, such transition processes are typically described by traveling fronts. Depending on the geometry of their level sets, traveling fronts can be divided into planar and nonplanar fronts. A substantial literature has been devoted to the existence, stability and other qualitative properties of these fronts, see \cite{aw,bns,fm,hn1,vvv,bh,hm,hmr,nt,t1}.
A characteristic feature of traveling fronts is that they converge uniformly in time to steady states away from their level sets. Motivated by this property, Berestycki and Hamel introduced in their pioneering works \cite{BH1,BH2} the notion of transition fronts (see also \cite{S} for the one-dimensional case), thereby extending the classical concept of traveling fronts to a more general framework. Over the past two decades, transition fronts have become a major topic in the study of propagation dynamics in complex geometries, and we will review related results in Section \ref{ssss2}.

In view of the importance of transition fronts in the study of propagation phenomena, we investigate the existence and qualitative properties of transition fronts for reaction-diffusion equations of the form
\begin{equation}\label{te1.1}
\begin{cases}
u_t=\Delta u+f(u), & t\in\mathbb{R},\ x\in\Omega,\\
u_\nu=0, & t\in\mathbb{R},\ x\in\partial\Omega
\end{cases}
\end{equation}
posed in certain  smooth unbounded domains $\Omega\subset\mathbb{R}^N$ $(N\ge2)$.
Here, $x=(x_1,\ldots,x_N)$, $u_t=\partial u/\partial t$, $\Delta=\sum_{j=1}^N \partial^2/\partial x_j^2$ denotes the Laplacian operator, $\nu=\nu(x)$ stands for the outward unit normal vector on $\partial\Omega$ and $u_\nu=\partial u/\partial \nu$.
Throughout this paper, we assume that $f$ is a combustion nonlinearity, that is,
$f\in C^1([0,1])$ and there exists ignition temperature $\theta\in(0,1)$ such that
\be\lb{coma}
 f(s) \equiv 0\ \text { for }\ s\in[0, \theta] \cup\{1\},\quad
f(s)>0\ \text { for }\ s\in(\theta, 1),\quad  \ f'(1)<0.
\ee
For mathematical purpose, we further extend $f$ to a $C^1(\mathbb R)$ function by defining
\begin{equation}\label{ext}
f(s)=0 \quad \hbox{for } s<0,
\qquad
f(s)=f'(1)(s-1) \quad \hbox{for } s>1.
\end{equation}
This type of nonlinearity is commonly used in combustion models based on Arrhenius kinetics with a low-temperature cut-off and on the law of mass action.

In the case $\Omega=\mathbb R^N$ (without boundary condition), there exist several types of transition patterns between the equilibria $0$ and $1$ realized by traveling fronts. As shown in \cite{aw,bns}, problem \eqref{te1.1} admits a unique planar traveling front (up to shifts) connecting $0$ and $1$ of the form
$
u(t, x)=\phi\left(x \cdot e-c_f t\right)$,
where $e\in \mathbb S^{N-1}$ is the direction of propagation ($\mathbb{S}^{N-1}$ denotes the unit sphere of $\mathbb{R}^{N}$), $c_f\in\R$ is the speed of the front, and the profile $\phi: \mathbb{R} \rightarrow[0,1]$ is given by
\begin{equation}\label{p-front}
\begin{cases}
\phi ''(\xi)+c_f \phi '(\xi)+f\left(\phi (\xi)\right)=0, &\xi\in \mathbb{R},\\
\phi'(\xi)<0, &\xi\in\R,\\
\phi(-\infty)=1,\ \phi(+\infty)=0.
\end{cases}
\end{equation}
 Fife and McLeod \cite{fm} proved that the sign of $c_f$ coincides with that of $\int_0^1 f(s)\,ds$, and hence $c_f>0$.
For any $\alpha\in(0,1)$, the level set of  $\phi$ is characterized by
$$
E_\alpha(t)
:=\{x\in\mathbb R^N:\ u(t,x)=\alpha\}
=\{x\in\mathbb R^N:\ x\cdot e=\phi^{-1}(\alpha)+c_f t\},\quad t\in\mathbb R,
$$
where $\phi^{-1}$ denotes the inverse function of $\phi$. In particular, these level sets form a family of parallel hyperplanes orthogonal to $e$. This reflects the planar structure of the front in the moving frame with speed $c_f$.
According to \cite{aw}, the front $\phi$ admits the following exponential estimates and asymptotic behavior: there exist positive constants $K_i\ (i=1,2,3,4)$ and $\Lambda$ such that
\begin{equation} \label{estimate2}
 \begin{split}
 K_1 e^{\Lambda\xi}\leq(1-\phi(\xi)),|\phi'(\xi)|,|\phi''(\xi)|\leq  K_2 e^{\Lambda\xi},&\quad
 \forall  \  \xi<0,\\
 K_3 e^{-c_f\xi}\leq\phi(\xi),|\phi'(\xi)|,|\phi''(\xi)|\leq K_4 e^{-c_f\xi},&\quad
 \forall\ \xi>0
 \end{split}
\end{equation}
and
\begin{align}\label{estimates4}
\lim\limits_{\xi\to+\infty}\frac{\phi'(\xi)}{\phi(\xi)}=-c_f,~
\lim\limits_{\xi\to+\infty}\left|\frac{\phi''(\xi)}{\phi(\xi)}\right|=c_f^2.
\end{align}
Furthermore, Wang and Bu \cite{WBu} established the existence of V-shaped and pyramidal traveling fronts connecting $0$ and $1$ for problem \eqref{te1.1} in $\mathbb{R}^N$, whose level sets are no longer hyperplanes. In addition, Bu, Guo and Wang \cite{BGW} discovered a new type of transition  pattern between the equilibria $0$ and $1$ in $\mathbb{R}^N$ beyond classical traveling fronts. In fact, they showed that problem \eqref{te1.1} admits an entire solution which propagates as three planar traveling fronts in three directions. In particular, this entire solution is a transition front connecting $0$ and $1$. 

The situation in the presence of an obstacle $K$ is also of great interest, where $K\subset\mathbb{R}^N$ is a compact set with smooth boundary. Usually, the domain $\Omega=\mathbb{R}^N\backslash K$ is called an exterior domain. Jia, Wang and Zhang \cite{jwz} showed that problem \eqref{te1.1} admits an entire solution originating from any given transition front defined in $\mathbb{R}\times\mathbb{R}^N$ as $t\to-\infty$, and it converges to a stationary solution as $t\to+\infty$. Later, Yan and Sheng \cite{ys} established the existence of transition fronts of problem \eqref{te1.1} in exterior domains. More precisely, they indicated that for a star-shaped or directionally convex obstacle, there exists a unique transition front $u$ connecting $0$ and $1$ such that $u(t,x)$ converges to a V-shaped traveling front uniformly in $x\in\bar\Omega$ as $t\to\pm\infty$. Moreover, they obtained that all transition fronts connecting $0$ and $1$ for problem \eqref{te1.1}  share the same global mean speed, which is exactly the speed $c_f$ of the planar traveling front $\phi$. These results imply that the interaction between traveling fronts and the geometry of the domain may give rise to new phenomena in the transition between steady states.
It is therefore natural to understand how such phenomena are affected by more complex geometries beyond exterior domains. 

\begin{figure}[htbp]
\centering
\includegraphics[width=7cm]{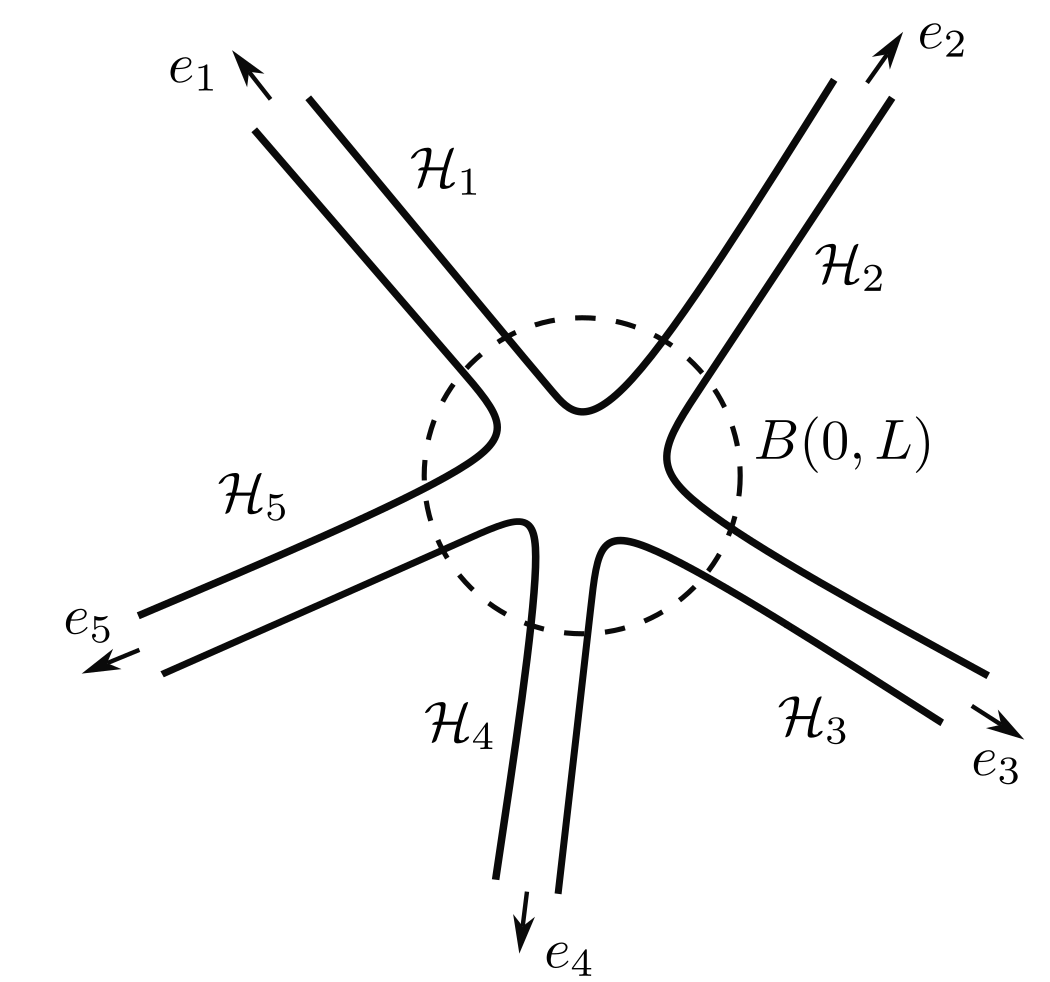}
\caption{A domain with five cylindrical branches in $\R^2$.}
\label{fig1}
\end{figure}

In this paper,  we are interested in unbounded domains consisting of several semi-infinite cylindrical branches connected through a bounded region. A fundamental problem is whether each branch acts as a propagation channel through which transition fronts can enter or leave the domain, and how the geometry of the domain influences the dynamical behavior.
Motivated by this problem, we consider a class of unbounded domains with multiple cylindrical branches. More precisely, we assume that $\Omega\subset\R^N$ is an unbounded domain with smooth boundary such that there exist a real number $L>0$, an integer $m\ge2$ and $m$ cylindrical branches $\{\mathcal H_i\}_{i=1}^m$  satisfying
\be\lb{defmulti}
\left\{\begin{array}{l}
\mathcal{H}_{i} \cap \mathcal{H}_{j}
=\emptyset\quad  \text { for every } i \neq j \in\{1, \cdots, m\}, \\
\mathcal{H}_{i} \backslash B(0, L)\quad 
\text { is connected for every } i \in\{1, \cdots, m\},\\
\Omega \backslash B(0, L)=\bigcup_{i=1}^{m}
\mathcal{H}_{i} \backslash B(0, L),
\end{array}\right.
\ee
where  $B(0,L)$ denotes the open ball of radius $L$ centered at the origin of $\R^N$.
For each $i=1,\cdots,m$, the $i$-th cylindrical branch is defined by
\be\lb{defhi1}
\mathcal H_i
=\left\{
x\in\R^N:\,
x-(x\cdot e_i)e_i\in\omega_i,\,
x\cdot e_i>0
\right\}+x_i,
\ee
the cross section $\omega_i\subset\R^{N-1}$  is a nonempty bounded domain with smooth boundary, axial direction $e_i\in\mathbb S^{N-1}$ is a unit vector, and shift $x_i\in\R^N$ is a given point determining the location of  $\mathcal H_i$.
Figure \ref{fig1} illustrates a simple domain with five cylindrical branches in $\mathbb R^2$.

\subsection{Transition fronts: some known results}\label{ssss2}

Since the introduction by Berestycki and Hamel \cite{BH1,BH2}, transition fronts have been extensively studied in various settings. We begin by recalling the definition of transition fronts for problem  \eqref{te1.1}, and then briefly review some previous results relevant to the present work.
For any subsets $A,B\subset\Omega$, we set
\begin{equation*}
d_{\Omega}(A,B)=\inf\{d_{\Omega}(x,y): (x,y)\in A\times B\},\quad
d_{\Omega}(x,A)=d_{\Omega}(\{x\},A)\ \text{ for }\ x\in\Omega,
\end{equation*}
 where $d_{\Omega}$ denotes the geodesic distance in $\Omega$. 
 Assume that for each $t\in\mathbb R$, there exist an interface $\Gamma_t$ and two open subsets $\Omega_t^\pm\subset\Omega$ satisfying the following properties.
\begin{itemize}
\item[(A1)] The interface $\Gamma_t$ separates $\Omega$ into two disjoint unbounded components $\Omega_t^+$ and $\Omega_t^-$, that is,
\begin{equation}\label{eq1.3}
\begin{cases}
\Omega_t^-\cap \Omega_t^+=\emptyset, \
\partial \Omega_t^-\cap \Omega=\partial \Omega_t^+\cap \Omega=:\Gamma_t, \
\Omega_t^-\cup \Gamma_t \cup \Omega_t^+=\Omega,\\
\sup\big\{d_{\Omega}(x,\Gamma_t): x\in \Omega_t^+\big\}=\sup\big\{d_{\Omega}(x,\Gamma_t): x\in \Omega_t^-\big\}=+\infty.
\end{cases}
\end{equation}

\item[(A2)] For any point $x\in\Gamma_t$, it is not too far from the centers of two large balls contained in $\Omega_t^+$ and $\Omega_t^-$, in the sense that
\begin{equation}\label{eq1.4}
\begin{cases}
\inf \{\sup\{d_{\Omega}(y,\Gamma_t): y\in \Omega_t^+, d_{\Omega}(y,x)\leq r\}: t\in \mathbb{R},\ x\in \Gamma_t \} \rightarrow  +\infty,\\
\inf\{\sup\big\{d_{\Omega}(y,\Gamma_t): y\in \Omega_t^-, d_{\Omega}(y,x)\leq r\big\}: t\in \mathbb{R},\ x\in \Gamma_t \} \rightarrow+\infty
\end{cases}
\ \ \text{as } r\to+\infty.
\end{equation}

\item[(A3)] The interface $\Gamma_t$ is included in finitely many graphs. Namely,
 there exists an integer $\bar n\ge1$ such that, for every $t\in\mathbb R$, there exist open subsets $\omega_{i,t}\subset\mathbb R^{N-1}$, continuous functions $\psi_{i,t}:\omega_{i,t}\to\mathbb R$ and rotations $R_{i,t}$ of $\mathbb R^N$ satisfying
\begin{equation}\label{eq1.6}
\Gamma_t \subset \bigcup_{1\leq i\leq\bar n} R_{i,t}(\{x=(x',x_N)\in \mathbb{R}^N:\ x'\in \omega_{i,t},\ x_N=\psi_{i,t}(x')\}).
\end{equation}
\end{itemize}
 
 Condition (A1) implies in particular that the interface $\Gamma_t$ is nonempty for every $t\in\mathbb R$. From condition (A2), it follows that for every $M>0$, there exists $r_M>0$ such that for all $t\in\mathbb R$ and $x\in\Gamma_t$, there are points $y^\pm\in\Omega_t^\pm$ satisfying
\begin{equation}\label{eq1.5}
y^\pm\in\Omega_t^\pm,\quad
d_\Omega(x,y^\pm)\le r_M,\quad
d_\Omega(y^\pm,\Gamma_t)\ge M.
\end{equation}
Without loss of generality, we can assume that the map $M\mapsto r_M$ is nondecreasing. Moreover, condition (A3) is imposed to rule out interfaces with infinitely many twists.

\begin{defi}[\!\!{\rm\cite{BH1,BH2}}]\label{td1.1}
Let $p^\pm:\R\times\bar\Omega\to\R$ be two classical solutions of problem \eqref{te1.1}.
A transition front connecting $p^-$ and $p^+$ for problem \eqref{te1.1} is a classical solution
$u\not\equiv p^\pm$ for which there are the families of sets $(\Omega_t^\pm)_{t\in\mathbb R}$ and $(\Gamma_t)_{t\in\mathbb R}$ satisfying conditions {\rm(A1)-(A3)}, and such that for any $\varepsilon>0$, there exists $M_\varepsilon>0$ satisfying
\begin{equation*}
\begin{cases}
\forall\,t\in \mathbb{R},\ \ \forall\,x\in \overline{\Omega_t^+}, \ \ \left(d_{\Omega}(x,\Gamma_t)\geq M_{\varepsilon}\right)\Rightarrow |u(t,x)-p^+(t,x)|\leq \varepsilon ,\\
\forall\,t\in \mathbb{R},\ \ \forall\,x\in \overline{\Omega_t^-}, \ \ \left(d_{\Omega}(x,\Gamma_t)\geq M_{\varepsilon}\right)\Rightarrow |u(t,x)-p^-(t,x)|\leq \varepsilon .
\end{cases}
\end{equation*}
Furthermore, $u$ is said to have a global mean speed $\gamma(\geq 0)$ if
$$
\frac{d_{\Omega}(\Gamma_t,\Gamma_s)}{|t-s|}\rightarrow \gamma \ \ \text{as}\ \ |t-s|\rightarrow +\infty.
$$
\end{defi}

Classical traveling fronts are typical examples of transition fronts for $\Omega=\R^N$. For instance, the planar traveling front $u(t,x)=\phi(x\cdot e-c_f t)$ satisfying \eqref{p-front} is a transition front connecting the equilibria $0$ and $1$ of \eqref{te1.1}, with $\Gamma_t=\{x\in\R^N:x\cdot e=c_ft\}$, $\Omega_t^\pm=\{x\in\R^N:\pm(x\cdot e-c_ft)<0\}$ and the global mean speed $c_f$.
However, transition fronts in $\R^N$ are not restricted to traveling fronts. Consider problem \eqref{te1.1} with a bistable nonlinearity, namely,
 \begin{equation*}
 \begin{cases}
 \exists\ \theta\in(0,1),\quad f(0)=f(\theta)=f(1)=0,\quad 
 f'(0)<0, \quad f'(1)<0,\\
f<0 \text{ on }(0, \theta),\quad
   f>0\text{ on }(\theta, 1),\quad \int_{0}^{1} f(s) d s>0.
   \end{cases}
\end{equation*}
Hamel \cite{h1} constructed a class of nonstandard transition fronts of \eqref{te1.1} by rotating a V-shaped traveling front and attaching the reflected V-shaped front to it. These transition fronts behave like a combination of three planar traveling fronts as $t\to -\infty$ and converge to a V-shaped traveling front as $t\to +\infty$. In the same work, he further proved that every transition front connecting $0$ and $1$ possesses the same global mean speed $c_f$.
Subsequently, similar existence results were extended to several other frameworks, see \cite{BGW,sg,sw,sww}. More recently, Guo and Wang \cite{GWang} constructed another class of transition fronts of \eqref{te1.1} by mixing finitely many planar traveling  fronts, whose interfaces look like unclosed polyhedra.  
For more results on the existence and stability of transition fronts in $\R^N$,
 we can refer to \cite{dhz,hr1,mnrr,NR,NR1,nrrz,ss}
 and references therein.

Transition fronts have also been investigated in unbounded domains with complex geometry. If $\Omega$ is an exterior domain, it has been shown that bistable problem \eqref{te1.1} admits several distinct types of transition fronts. Berestycki, Hamel and Matano \cite{bhm} established the existence, monotonicity and uniqueness of an entire solution  originating from a planar traveling front as $t\to-\infty$. For a general compact obstacle, they proved that this entire solution is a transition front connecting $0$ and $u_\infty$, where $u_\infty\in(0,1]$ is the classical solution of
\begin{equation}\label{uwuqiong}
\Delta u_\i+f(u_\i)=0 \ \text{ in }\ \bar\Omega,\quad
(u_\i)_\nu=0 \ \text{ on }\partial\Omega\quad
\text{ and }\quad
\lim_{|x|\to+\i}u_\i(x)=1.
\end{equation}
On the other hand, under suitable geometric assumptions on the obstacle, they showed that this entire solution is a transition front connecting $0$ and $1$, and converges to the same planar traveling front as $t\to+\infty$. Furthermore, Guo, Hamel and Sheng \cite{ghs} proved that the global mean speed of any transition front of  bistable problem \eqref{te1.1} connecting $0$ and $u_\infty$ in exterior domains coincides with the unique planar wave speed.
Moreover, Guo and Monobe \cite{gm} established the existence of entire solutions in exterior domains originating from any transition front of bistable problem \eqref{te1.1} defined in $\mathbb R\times\mathbb R^N$. In particular, they obtained a transition front converging to a V-shaped traveling front as $t\to\pm\infty$ for certain obstacles.
Following these foundational works, a number of further studies on propagation phenomena for reaction-diffusion equations in exterior domains have appeared in recent years. As mentioned in Section \ref{ssss1}, the interaction of planar and V-shaped traveling fronts with obstacles for combustion reaction-difusion equations was investigated in \cite{jwz} and \cite{ys}, respectively. For more results on transition fronts in exterior domains, we refer to \cite{Li,jsw,ys2,JBB} and the references therein.

Cylindrical domains provide another important framework for studying the interaction between geometry and propagation. Berestycki, Bouhours and Chapuisat \cite{bbc} investigated how cylindrical domains with varying cross sections affect the propagation of planar traveling fronts for bistable problem \eqref{te1.1}. More precisely, they considered domains of the form
$$
\Omega=\{(x_1,x')\in\R^N:\,x_1\in\R,\,x'\in\omega (x_1)\subset\R^{N-1}\},
$$
where $\omega$ is independent of $x_1$ for $x_1<0$. They established the existence and uniqueness of a time-increasing entire solution $u$ that propagates from the left  part of the cylinder as a planar traveling front and locally converges  to $u_\infty$ as $t\to+\infty$, where $u_\i$ is the stationary solution  given by \eqref{uwuqiong}. In the case of domains with decreasing cross section, they showed that this entire solution is a transition front connecting $0$ and $1$, and that complete propagation occurs in the sense that
$$
u(t,x)\to1\quad \text{ locally uniformly in }\ x\in \bar\Omega\ \text{ as }\ t\to+\infty.
$$
Moreover, they proved that if the cylinder contains a sufficiently narrow passage, then propagation may be blocked in the sense that
$$
\lim_{t\to+\i}u(t,\cdot)= u_\i\text{ in } C_{loc}^2(\Omega),
\quad
u_\i(x)\to0 \text{ as }x_1\to+\i,\quad x=(x_1,x')\in\Omega.
$$
For domains with multiple cylindrical branches defined in \eqref{defmulti}, Guo, Hamel and Sheng \cite{ghs} established the existence of transition fronts for the bistable problem \eqref{te1.1}. In fact, these transition fronts take the form of planar traveling fronts in some branches as $t\to -\infty$, whose propagation directions are opposite to the axial directions of the corresponding  branches. And they converge in the remaining branches to planar traveling fronts  as $t\to +\infty$, whose propagation directions are consistent with the axial directions of the corresponding  branches. Moreover, they showed that every transition front connecting $0$ and $1$ has the same global mean speed.
Furthermore, Guo \cite{g1} investigated the existence and uniqueness of transition fronts for general bistable reaction-diffusion-advection equations in domains with multiple cylindrical branches of varying cross sections.
Since then, the above results in cylindrical domains have been extended to a variety of situations, including time-periodic bistable reaction-diffusion equations, monotone bistable reaction-diffusion systems, and models arising from biology, see \cite{gfl,SHENGW,swwcv,LXL} for further details.
To the best of the authors' knowledge, transition fronts of combustion equations have not been considered in cylindrical domains. This problem is addressed in the present paper.

Besides exterior and cylindrical domains, propagation phenomena for reaction-diffusion equations have also been studied in other classes of unbounded domains. For bistable problem \eqref{te1.1}  without Neumann boundary conditions, Jimbo and Morita \cite{JM} considered a domain made up of a finite union of half-lines emanating from a common end point. They proved the existence of entire solutions which converge to planar traveling fronts in some of the half-lines and to $0$ in the remaining ones as $t \to -\infty$. Moreover, they provided a condition under which these entire solutions are blocked as $t \to +\infty$.
They further showed in \cite{JM2} that, under suitable conditions, the entire solution converges to planar traveling fronts in each output branch far from the junction as $t\to+\i$. In a subsequent work \cite{JM3}, they verified that in an unbounded metric graph consisting of half-lines and triple junctions, the blocking of propagation depends not only on the geometry of the junctions but also on the distances between neighboring junctions.
In funnel-shaped domains consisting of a straight part and a conical part with a positive opening angle, Hamel and Zhang \cite{FZ} proved the existence and uniqueness of an entire solution to bistable problem  \eqref{te1.1} that behaves like a planar traveling front propagating from the straight part of the domain as $t \to -\infty$. They further established a dichotomy between complete propagation and blocking, and provided sufficient geometric conditions under which each of these two cases occurs. Moreover, they proved that in the complete propagation case, the solution is a transition front connecting $0$ and $1$ with global mean speed equal to the unique speed of planar traveling fronts. They also showed that in the blocking case, the solution remains a transition front connecting  $0$ and $1$, but possesses no global mean speed.
Later, Li, Zhang and Zhang \cite{Li2} extended these results to time-periodic bistable reaction-diffusion equations.
Propagation and blocking phenomena for combustion reaction-diffusion equations in these classes of unbounded domains remain not fully understood and will be studied in a forthcoming work.

\subsection{Main results}\label{ssss3}
The purpose of this paper is to establish the existence of transition fronts for the combustion reaction-diffusion equation \eqref{te1.1}  in domains with multiple cylindrical branches, and to investigate the influence of the geometry of the domain on propagation. In what follows, we assume that $\Omega$ is a domain with $m\,(\ge2)$ cylindrical branches defined by \eqref{defmulti}. We now state the main results. 

Our first result concerns the existence of front-like entire solutions. More precisely, we show that problem \eqref{te1.1} admits entire solutions that converge, as $t\to -\infty$, to planar traveling fronts (moving opposite to the axial directions of the corresponding branches) in some branches, and to $0$ in the remaining parts of the domain.
\begin{theo}\lb{existence-}
Let
 $I$ and
$J$ be two non-empty subsets of index set $\{1, \cdots, m\}$ such that
 $I \cap J=\emptyset$ and $I \cup J=\{1, \cdots, m\}$.
For any $\beta\in(0,1)$,  problem \eqref{te1.1} admits
 a unique time increasing entire solution $u(t, x)$ satisfying
\begin{equation}\lb{dm0}
\begin{cases}
\fr{u(t, x)-\phi\left(-x \cdot e_{i}-c_f t\right)}{\phi^\beta\left(-x \cdot e_{i}-c_f t\right)}{\rightarrow}  0 \ \
&\text {uniformly in }\overline{\mathcal{H}_{i}} \cap \bar{\Omega} \text { for all } i \in I, \\
\fr{u(t, x)}{\phi^\beta\left(x \cdot e_{j}-c_f t\right)}{\rightarrow}  0 \ \
&\text {uniformly in }\overline{\mathcal{H}_{j}} \cap \bar{\Omega}\text { for all } j \in J, \\
\fr{u(t, x)}{\phi^\beta\left(x\cdot e_j-c_f t\right)} {\rightarrow}0 \ \
&\text {uniformly in } \overline{\Omega \backslash\bigcup_{s=1}^m \mathcal{H}_{s}}\text{ for all $j\in J$}
\end{cases}
\ \ \text{as}\ t \rightarrow-\infty.
\end{equation}
\end{theo}

The next result is devoted to the large time behavior of the entire solution $u(t,x)$ obtained above and the existence of transition fronts. 
More precisely, we show that if $u(t,x)$ propagates completely in the sense of
\begin{equation}\label{coms}
u(t,x)\to 1 \quad \text{as } t\to +\infty \text{ locally uniformly in } x\in\bar\Omega,
\end{equation}
then $u(t,x)$ is a transition front and converges to the planar traveling front $\phi(x\cdot e_j - c_f t)$ (up to a finite shift) in branch $\mathcal H_j$  ($j\in J$) as $t\to +\infty$.

\begin{theo}\lb{existence}
Let $\beta\in(0,1)$ be any constant and $u(t,x)$ be the entire solution given
in Theorem \ref{existence-}.
If $u$ propagates completely in the sense of \eqref{coms},
then it is a transition front connecting $ 0$ and $ 1$ in the sense of Definition \ref{td1.1} with   $\left(\Gamma_{t}\right)_
{t \in \mathbb{R}}$ and $\left(\Omega_{t}^{\pm}\right)_{t \in
\mathbb{R}}$ defined by
\be\lb{dm1}
\bc
\Gamma_{t}=\bigcup_{i \in I}\left\{x \in \mathcal{H}_{i}
\cap \Omega: x \cdot e_{i}=-c_{f}t+A\right\},\ \  & t \leq 0,\\
\Gamma_{t}=\bigcup_{j \in J}\left\{x \in \mathcal{H}_{j}
\cap \Omega: x \cdot e_{j}=c_{f} t+A\right\},\  \  & t>0
\ec
\ee
and
\be\lb{dm2}
\begin{cases}\Omega_{t}^{+}=\bigcup_{i \in I}\left\{x
\in \mathcal{H}_{i} \cap \Omega: x \cdot
e_{i}>-c_{f}t+A\right\}, \ \ & \Omega_{t}^{-}
=\Omega \backslash \overline{\Omega_{t}^{+}},
\quad  t \leq 0, \\
\Omega_{t}^{-}=\bigcup_{j \in J}\left\{x \in \mathcal{H}_{j}
\cap \Omega: x \cdot e_{j}>c_{f} t+A\right\},\ \ &
\Omega_{t}^{+}=\Omega \backslash \overline{\Omega_{t}^{-}},
\quad t>0
\end{cases}
\ee
for some $A>0$. Moreover, there exist some real
numbers $(\tau_{j}^*)_{j \in J}$ such that
\be\lb{dm3}
\bc
\fr{u(t, x)- \phi\left(x \cdot e_{j}-c_{f} t+\tau_{j}^*\right)}{\phi^\beta\left(x \cdot e_{j}-c_{f} t+\tau_{j}^*\right)}
{\rightarrow}  0 \ &
\text { uniformly in } \overline{\mathcal{H}_{j}} \cap \bar{\Omega}\text{ for all }j\in J, \\
\frac{1-u(t, x)}{\phi^\beta(-x\cdot e_i-c_f t)} {\rightarrow}
 0 \ &\text { uniformly in } \overline{\mathcal{H}_{i}} \cap \bar{\Omega}\text{ for all }i\in I, \\
\frac{1-u(t, x)}{\phi^\beta(-x\cdot e_i-c_f t)} {\rightarrow}
 0 \ &\text { uniformly in } \overline{\Omega \backslash
\bigcup_{s=1}^m \mathcal{H}_{s}} \text{ for all } i\in I
\ec
\ \ \text{as}\ t \rightarrow+\infty.
\ee

\end{theo}

Since $\phi \in (0,1)$ and in view of \eqref{dm0} and \eqref{dm3}, we summarize the asymptotic behavior of the entire solution
$u(t,x)$ under the complete propagation condition in Table \ref{asymp}.
\begin{table}[htbp]
\centering
\renewcommand{\arraystretch}{1.3}
\setlength{\tabcolsep}{10pt}
\begin{tabular}{lcc}
\toprule
 & $u(t,x)$ $(t\to -\infty)$ 
 & $u(t,x)$ $(t\to +\infty)$ \\ 
 \midrule

 $x\in\bar{\mathcal H_i}\,(i\in I)$  
 & $\phi(-x\cdot e_i-c_f t)$
 &1\\

 $x\in\bar{ \Omega \setminus \bigcup_{s=1}^m \mathcal H_s}$ 
 &0
 &1\\ 
 
 $x\in\bar{\mathcal H_j}\,(j\in J)$  
  &0
 & $\phi(x\cdot e_j-c_f t+\tau_j^*)$\\
\bottomrule
\end{tabular}
\caption{Asymptotic behavior of the entire solution under complete propagation.}
\label{asymp}
\end{table}
It is easy to see that the set $I$ corresponds to the branches from which propagation originates as $t \to -\infty$, whereas $J$ corresponds to the outgoing branches. In addition, we will provide geometric conditions ensuring that the complete propagation assumption holds in Theorems \ref{suffc1} and \ref{suffc2}.

The above result constructs a transition front connecting $0$ and $1$ with global mean speed $c_f$. However, in domains with multiple cylindrical branches, it is not known whether the global mean speed depends on the choice of transition front. The following theorem shows that  all transition fronts connecting $0$ and $1$ propagate with the same global mean speed.
For each $i \in\{1, \cdots, m\}$,
let $u^{i}: \mathbb{R} \times \bar{\Omega}
 \rightarrow(0,1)$ be
the time-increasing front-like solution of \eqref{te1.1} coming
from the branch $\mathcal{H}_{i}$, that is,  for any $\beta\in(0,1)$, there holds
\be\lb{defhi}
\bc
\fr{u^{i}(t, x)-\phi\left(-x \cdot e_{i}-c_{f} t\right)}{\phi^\beta\left(-x \cdot e_{i}-c_{f} t\right)}
 {\rightarrow} 0
\ \ &\text {uniformly in } \overline{\mathcal{H}_{i}}
\cap \bar{\Omega}, \\
\fr{u^{i}(t, x)}{\phi^\beta\left(x \cdot e_{j}-c_{f} t\right)}{\rightarrow}
0\ \ &\text{uniformly in }\overline{\mathcal H_j}\cap\bar{\Omega} \text{ for all $j\neq i$},\\
\fr{u^{i}(t, x)}{\phi^\beta\left(x \cdot e_{j}-c_{f} t\right)}  {\rightarrow}
0 \ \ &\text {uniformly in } \overline{\Omega \backslash\bigcup_{s=1}^m
\mathcal{H}_{s}} \text{ for all }j\neq i
\ec  \ \ \text{as $t \rightarrow-\infty$}.
\ee

\begin{theo}\lb{gms}
Assume that,
for every $i \in$ $\{1,\cdots,m\}$, the time increasing
solution $u^{i}$ of \eqref{defhi} propagates completely
in the sense of \eqref{coms}. Then any transition front $u$
of \eqref{te1.1} connecting $0$ and $1$ in the sense of Definition \ref{te1.1} also propagates completely
 and admits $c_{f}$ as the unique global mean
speed.
\end{theo}

The assumption of complete propagation appearing in Theorems \ref{existence} and \ref{gms} is nontrivial. Indeed, propagation may be blocked even in domains with only two cylindrical branches \cite{bbc}. Our final results provide two geometric conditions ensuring complete propagation. The first concerns star-shaped domains with suitable connectivity assumptions. The second shows that complete propagation holds for sufficiently large scalings of any given domain with multiple cylindrical branches.
Before stating the following theorems, we introduce some notation. For every $i\neq j\in\{1,\cdots,m\}$, let
$P_{i,j}:\mathbb{R}\to\Omega$ be a continuous path connecting the branches $\mathcal H_i$ and $\mathcal H_j$, in the sense that, for every $A\in\mathbb{R}$, there exists $\sigma>0$ such that
\[
\begin{cases}
P_{i,j}(s)\in\{x\in\mathcal H_i\cap\Omega:\,x\cdot e_i\ge A\}, & s\le-\sigma,\\
P_{i,j}(s)\in\{x\in\mathcal H_j\cap\Omega:\,x\cdot e_j\ge A\}, & s\ge \sigma.
\end{cases}
\]
Without loss of generality, we assume that $P_{i,j}(\mathbb R)=P_{j,i}(\mathbb R)$ for each $i\neq j$.
When $m\ge3$, any two paths $P_{i,j}$ and $P_{i,k}$ with $j\neq k\in\{1,\ldots,m\}\setminus\{i\}$ are assumed to coincide inside the branch $\mathcal H_i$ sufficiently far from the origin. Namely, there exists $s_i\in\mathbb R$ such that for each $i$,
\[
P_{i,j}(s)=P_{i,k}(s) \quad
\text{ for all }
s\le s_i\ \text{ and }\
j\neq k\in\{1,\ldots,m\}\setminus\{i\}.
\]

\begin{theo}\lb{suffc1}
Assume that the domain $\Omega$ with $m\,(\geq2)$ cylindrical
branches satisfying the following conditions:
\begin{itemize}
\item[\rm(1)] There is a constant $R>0$ for which the paths $P_{i, j}\,(1\leq i\neq j\leq m)$  can be chosen so that
$$
B\left(P_{i, j}(s), R\right) \subset \Omega \  \text { for all }
s \in \mathbb{R} \quad \text { and }\quad 0 \in \bigcup_{1 \leq i
\neq j \leq m} P_{i, j}(\mathbb{R}).
$$
\item[\rm(2)] $\Omega$ is star-shaped with respect to $0$, that is,  $0\in\Omega$ and
$sy\in\bar\Omega$ for all $y\in\partial\Omega$ and $s\in[0,1)$.
\end{itemize}
Then for every $i \in\{1,\cdots,m\}$, the front-like solution $u^{i}$ satisfying
\eqref{defhi} propagates completely in the sense of \eqref{coms}. Furthermore, any entire solution $u$ given in Theorem \ref{existence-}  propagates
completely.
\end{theo}

\begin{theo}\lb{suffc2} Let $\Omega$ be any domain with $m\,(\geq2)$ cylindrical branches satisfying \eqref{defmulti}.
Then there exists a constant $R_0>0$ such that, for every $R\ge R_0$, every
$x_0\in\mathbb R^N$ and every $i\in\{1,\ldots,m\}$, the front-like
solution $u^i$ satisfying \eqref{defhi} in
$R\Omega_0+x_0$ 
propagates completely. Furthermore, any entire solution constructed in
Theorem \ref{existence-} for the domain $R\Omega_0+x_0$  also propagates completely.
\end{theo}

The proofs in this paper rely mainly on the method of sub- and supersolutions. However, unlike the bistable case, the combustion nonlinearity is degenerate at the unstable equilibrium $0$ (namely,  $f'(0)=0$), which prevents us from verifying the required differential inequalities by linearization at $0$. To overcome this difficulty, we introduce perturbation terms involving fractional powers of the planar traveling front in the construction of sub- and supersolutions. Combined with the precise asymptotic behavior of the planar traveling front presented in Section \ref{ssss1}, this enables us to verify the differential inequalities associated with the constructed sub- and supersolutions. This also explains the appearance of the fractional powers $\phi^\beta$ in Theorems \ref{existence} and \ref{gms}.

Another difficulty is that the desired entire solution exhibits different asymptotic behaviors in the incoming branches, the outgoing branches and the intermediate part of the domain. To deal with these different behaviors, we construct refined sub- and supersolutions for different parts of the domain and different time intervals. Although the construction and verification of these sub- and supersolutions are complex, they provide an effective approach to the study of combustion reaction-diffusion equations in domains with multiple cylindrical branches.

We organize this paper as  {follows}. Section \ref{s2} is devoted to the existence and uniqueness of entire solutions emanating from planar traveling fronts in some branches, that is, we prove Theorem \ref{existence-}.  In Section \ref{s3}, we study the large time behavior of the entire solution given in Theorem \ref{existence-}, namely, we prove Theorem \ref{existence}. In Section \ref{s4}, we focus on the proof of the existence and  uniqueness of global mean speed, that is, Theorem \ref{gms}. In Section \ref{s6}, we show some geometrical conditions to guarantee the complete propagation, that is, we prove  Theorems \ref{suffc1}  and \ref{suffc2}.

\SE{Entire solution emanating from planar traveling fronts}\label{s2}
In this section, we shall consider the existence  and uniqueness of entire solutions emanating from  planar traveling fronts in some branches
with the help of the method of sub- and supersolutions.

\subsection{Construction of sub- and supersolutions}
We first state the definitions of sub- and supersolutions of \eqref{te1.1}.
\begin{defi}\label{dssd}
A function $u\in C^{0,1}((0,T)\times\bar\Omega)\cap C^{1,2}((0,T)\times\Omega) $ $(T>0)$
 is called a subsolution (supersolution)
of \eqref{te1.1} in $(0,T)\times\Omega$, if
\be\label{dss}\bc
u_t-\Delta u-f(u)\leq(\geq)0, &(t,x)\in(0,T)\times\Omega,\\
u_\nu\leq(\geq)0 , &(t,x)\in(0,T)\times\p\Omega.
\ec
\ee
If $ u $ and $v$ both are subsolutions (supersolutions) of \eqref{te1.1} in $(0,T)\times\Omega$, then
$\max(u,v)$ $(\min(u,v))$ is still  a  subsolution (supersolutions) of \eqref{te1.1} in $(0,T)\times\Omega$.
\end{defi}

The following comparison principle is obtained from \cite{FT,pw}.
 \begin{lem}\label{cp}
Let $\underline u$ and $\bar u$ be sub- and supersolutions of \eqref{te1.1} in $(0, T)\times\Omega$, $T>0$,
 respectively. If there holds $\underline u(0,x)\leq\bar u(0,x)$  for $x\in\Omega$,
then $\underline u(t,x)\leq \bar u(t,x)$ for $t\in(0,T)$ and $x\in\Omega$.
\end{lem}

For the sake of convenience,
we first introduce some preliminaries.
 Write
$$
\Lambda(\lambda)=\lambda^2+c_f \lambda.
$$
Take any $\beta\in(0,1)$. It is obvious that $\Lambda(-\beta c_f)<0$   due to
$\Lambda(0)=\Lambda(-c_f)=0$.
Choose a nonnegative $C^2$ function $\hat \zeta_i(x):\bar{\mathcal H_i}\to\mathbb R$ ($i\in\{1,\cdots,m\}$) with compact support in $\bar{\mathcal H_i}$ such that $\nabla \hat\zeta_i(x)\cdot \nu= 1$ for $x\in\p\mathcal H_i$. Then continuous functions  $\hat \zeta_i(x)$, $\nabla \hat \zeta_i(x)$, $\Delta\hat \zeta_i(x)$ are bounded in the sense of $L^\i$ norm and compactly supported. In fact, we can construct a truncated function as a such $\hat\zeta$ by applying the classical distance function in  \cite{gt} around the boundary $\partial\mathcal H_i$. Then there exists a large constant $\hat C$ such that  $\zeta_i(x)=\hat\zeta_i(x)+\hat C\gg0$ satisfies
\be\lb{defzeta}
\bc
\Delta \zeta_i(x)+2\e(1+ \e\| \phi'/\phi \r\|_{L^\i(\R)}\r)|\nabla \zeta_i(x)|\l \min\e(-\fr{f'(1)}{8},-\fr{\Lambda(-\beta c_f)}{8}\r) \zeta_i(x),\ \ x\in\bar{\mathcal H_i},\\
\nabla \zeta_i(x)\cdot \nu>0,\ \ x\in\p \mathcal H_i\cap\p\Omega,
\ec
\ee
where $|\nabla \zeta_i(x)|$ is the modulus of the vector $\nabla \zeta_i(x)$.
Denote
\begin{align}\label{sizeta}
\bar\zeta_0=\max_{1\l i\l m}\|\zeta_i\|_{L^\i(\bar{\mathcal H_i})}(>0)
\text{ and }
\tilde C=\min_{1\leq i\leq m}\inf_{x\in\bar{\mathcal H_i}}\zeta_i(x)(>0).
\end{align}
By \eqref{coma}  {and \eqref{ext}}, there exists a small constant $\varepsilon_0\in(0,\theta/2)$ such that
\be\lb{vp0}
\frac{3}{2}f'(1)\leq f'(u)\leq\frac{3}{4}f'(1), \ \forall u\in(1-2\varepsilon_0,1+2\varepsilon_0).
\ee
Define
\be\lb{deflf}
L_f= \sup _{u \in {(-2\vp_0,  {1+2\vp_0})}}|f'({u})|.
\ee
Recall that $L>0$ is given in \eqref{defmulti}. Take any $\vp>0$ such that
\begin{align}\label{defvp}
\vp<\min\e( c_f,\frac{\theta}{2},\frac{\vp_0}{2},
 \fr{c_f|f'(1)|}8, \fr{ c_f|\Lambda(-\beta c_f)|}8,\frac{|\Lambda(-\beta c_f)|}{24},
{\frac{\vp_0\bar\zeta_0}{\tilde C}}\r)(<1).
\end{align}
Define
\begin{align}\label{mu}
\mu= \frac{\vp}{c_f}(<1)\ \text{ and }\ \tilde\vp=\frac{\vp}{\bar\zeta_0}.
\end{align}
 {Hereafter, let
 $I$ and
$J$ be two non-empty subsets of $\{1, \cdots, m\}$ such that
 $I \cap J=\emptyset$ and $I \cup J=\{1, \cdots, m\}$.}
\begin{lem}\lb{exsub}
There exist $T<0$ and $\rho>0$ such that
 the function
$\underline{u}(t, x)$ defined by
$$
\underline{u}(t, x)= \begin{cases}\max
\left( \phi(\underline{\xi}_i(t,x))-
\tilde{\vp} e^{-\mu \left(x \cdot
e_{i}-L\right)} \zeta_i(x)  \phi^\beta(\underline{\xi}_i(t,x)),
  0\right)&
\text{in }\overline{\mathcal{H}_{i}}\cap\overline\Omega\text{ for all }i\in I, \\
  0 &
\text {in } \overline{\Omega\backslash\bigcup_{i\in I}
\mathcal{H}_{i}}\end{cases}
$$
is a subsolution of \eqref{te1.1} for all $t \leq T$ and
$x \in {\Omega}$,
where  $\underline{\xi}_i(t,x)=-x\cdot e_i
-c_ft+\rho e^{\vp t} {+1}$.
\end{lem}

\begin{pr}
Fix any $i\in I$. Let $ \vp' =\tilde C\tilde \vp$.
 Since $\tilde C\leq\bar\zeta_0$ by \eqref{sizeta}, it then follows from \eqref{mu} that $\vp'=\vp\tilde C/\bar\zeta_0e^{-\mu}\leq\vp$.
Since $\phi(-\i)=1$ and $\phi(+\i)=0$, there exists a constant $C>0$ such that
\begin{align}\lb{defC'}
 \phi^{1-\beta}\leq\vp' {\le \vp} \ \ \text{in $[C,+\i)$} \ \ \text{and}
\ \  \phi\g  1-\vp' { \g1-\vp} \ \ &\text{in $(-\i,-C]$}.
\end{align}
By \eqref{estimates4} and  $\phi'(-\i)=\phi''(-\i)=0$, even if it means increasing $C$,
one can assume that
\begin{align}
&-\fr{3}2c_f<\fr{\phi'}{\phi}<-\fr{1}2 c_f\ \text{ in }\ [C,+\i), \lb{betac1}\\
\fr{\phi''}{\phi}-\e(\fr{\phi'}{\phi}\r)^2\l -\fr14 &\Lambda(-\beta c_f),\
c_f\beta \fr{\phi'}{\phi}
+\beta^2 \e(\fr{\phi'}{\phi}\r)^2\l \fr34 \Lambda(-\beta c_f) \ \text{ in }\ [C,+\i) \lb{betac2}
\end{align}
and
\be\lb{defphi'}
\fr{\phi''}{\phi}\l -\fr18 f'(1) \
  \ \ \text{in $(-\i, -C]$}.
\ee
Define $\kappa>0$ such that
\be\lb{defk}
\phi '<-\kappa \ \ \text{in $[-C,C]$} .
\ee
Pick $\rho>0$ large enough such that
\be\lb{defrho}
\rho\kappa\geq e^{\mu(C+L)}\e(\mu ^2-f'(1)/8
+\beta \e\| \phi''/\phi \r\|_{L^\i(\mathbb R)} +L_f\r).
\ee
Take $T<0$ negative enough  such that
\be\lb{deft}
 c_f  T\l-L-C\ \ \text{and}\ \ \rho e^{\vp T}\l\min\left(1,-\frac{\Lambda(-\beta c_f)}{12c_f}\right).
\ee

Let us first verify that the function
$\underline{u}(t, x)$ is well defined.
For $t\l T$ and $x\in\bar {\mathcal H_i}\cap\bar\Omega$ such
that $x\cdot e_i\l L$, one gets from \eqref{deft} that
$\underline{\xi}_i(t,x)\geq-L -c_f T \g C$, hence
$ \phi^{1-\beta}(\underline{\xi}_i(t,x))\l \vp' $ by \eqref{defC'}.
In addition, for $x\in\bar {\mathcal H_i}\cap\bar\Omega$
such that $x\cdot e_i\leq L$, there holds
$\tilde \vp  e^{-\mu \left(x \cdot
e_{i}-L\right)} \zeta_i(x)
\g  \tilde \vp \tilde C =\vp'$ by \eqref{sizeta}. Hence $\underline{u}(t, x)=  0$ for
$(t,x)\in (-\i,T]\times\bar {\mathcal H_i}\cap\bar\Omega$ such
that $x\cdot e_i\leq L$. As a result, the function $\underline{u}$ is  continuous in $(-\i,T] \times \bar{\Omega}$. It is easy to see that $\underline{u}$ is of class $C^{1,2}$ in the set where it is positive.

We then check that $\underline u$ is a subsolution of \eqref{te1.1} in $(-\i,T]\times\bar\Omega$. Let us first check the boundary conditions.
By $\nu\cdot e_i=0$  on  $\partial \mathcal H_i\cap\p\Omega$ and
 \eqref{defzeta}, one obtains that
\[\underline u_\nu(t,x)=
-\tilde \vp  e^{-\mu \left(x \cdot
e_{i}-L\right)}   \phi^\beta(\underline\xi _i(t,x)) \nabla\zeta_i (x)\cdot \nu\l   0
\]
for $t\l T$ and $x\in\p\mathcal H_i\cap\p\Omega$.
It suffices to show that
\[
\mathscr L [\underline u](t,x):=\underline u_t(t,x)
-  \Delta \underline u  (t,x)-f (\underline u(t,x))\l 0
\]
for $t\l T$ and $x\in\bar\Omega$ such that $\underline u(t,x)> 0$.
After some computations, it follows  that
\begin{align*}
\mathscr L[\underline u](t,x)=
&\rho\vp e^{\vp t} \phi'(\underline{\xi}_i(t,x))
-\tilde \vp \beta e^{-\mu \left(x \cdot
e_{i}-L\right) }\zeta_i(x) \phi^{\beta-1}(\underline{\xi}_i(t,x))\phi'(\underline{\xi}_i(t,x))
\e(-c_f +\rho\vp e^{\vp t}\r)\\
&+ \tilde \vp  \mu ^2 e^{-\mu \left(x \cdot
e_{i}-L\right) } \zeta_i(x) \phi^{\beta}(\underline{\xi}_i(t,x)) -2 \tilde \vp  \mu   e^{-\mu \left(x \cdot
e_{i}-L\right)  } \nabla\zeta_i(x)\cdot e_i  \phi^{\beta}(\underline{\xi}_i(t,x))\\
&+ \tilde \vp    e^{-\mu \left(x \cdot
e_{i}-L\right) } \Delta\zeta_i(x) \phi^{\beta}(\underline{\xi}_i(t,x))
+ 2 \tilde \vp  \mu \beta  e^{-\mu \left(x \cdot
e_{i}-L\right) } \zeta_i(x)\phi^{\beta-1}(\underline{\xi}_i(t,x))\phi'(\underline{\xi}_i(t,x))\\
&-2 \tilde \vp    \beta  e^{-\mu \left(x \cdot
e_{i}-L\right) } \nabla\zeta_i(x)\cdot e_i \phi^{\beta-1}(\underline{\xi}_i(t,x))\phi'(\underline{\xi}_i(t,x))\\
&+\tilde \vp\beta(\beta-1) e^{-\mu \left(x \cdot
e_{i}-L\right) } \zeta_i(x) \phi^{\beta-2}(\underline{\xi}_i(t,x))(\phi'(\underline{\xi}_i(t,x)))^2\\
&+\tilde \vp \beta  e^{-\mu \left(x \cdot
e_{i}-L\right) } \zeta_i(x)\phi^{\beta-1}(\underline{\xi}_i(t,x))\phi''(\underline{\xi}_i(t,x))\\
&+f ( \phi(\underline{\xi}_i(t,x)))-f (\underline u(t,x))
.
\end{align*}
Since $\mu\in(0,1)$, we obtain from \eqref{defzeta} that
\begin{align}\label{zzzz}
&-2 \tilde \vp  \mu   e^{-\mu \left(x \cdot
e_{i}-L\right)  } \nabla\zeta_i(x)\cdot e_i  \phi^{\beta}(\underline{\xi}_i(t,x))
+ \tilde \vp    e^{-\mu \left(x \cdot
e_{i}-L\right) } \Delta\zeta_i(x) \phi^{\beta}(\underline{\xi}_i(t,x))\nonumber\\
&-2 \tilde \vp    \beta  e^{-\mu \left(x \cdot
e_{i}-L\right) } \nabla\zeta_i(x)\cdot e_i \phi^{\beta-1}(\underline{\xi}_i(t,x))\phi'(\underline{\xi}_i(t,x))\nonumber\\
\l& \tilde \vp    e^{-\mu \left(x \cdot
e_{i}-L\right) }  \phi^{\beta}(\underline{\xi}_i(t,x))\e(\Delta\zeta_i(x)+2\e(1+\e\|\phi'/\phi\r\|_{L^\i(\R)}\r)|\nabla\zeta_i(x)|  \r)\nonumber\\
\l&\min\left(-\fr{f'(1)}8,-\fr{\Lambda(-\beta c_f)}{8}\right) \tilde \vp    e^{-\mu \left(x \cdot
e_{i}-L\right) }  \phi^{\beta}(\underline{\xi}_i(t,x))\zeta_i(x).
\end{align}

If
$\underline{\xi}_i(t,x)\l -C$, it then can be deduced from \eqref{defC'} that $\phi(\underline{\xi}_i(t,x))\g 1-\vp'$,
 hence one gets from \eqref{defvp} and $x\cdot e_i\geq L$ that
\[\underline u(t,x)= \phi(\underline{\xi}_i(t,x))-
\tilde \vp  e^{-\mu \left(x \cdot
e_{i}-L\right) } \zeta_i(x)  \phi^\beta(\underline{\xi}_i(t,x))\g 1-\vp'-\vp\g 1-2\vp\g 1-\vp_0.\]
By the mean value theorem and \eqref{vp0}, there holds
\begin{align*}
&f ( \phi(\underline{\xi}_i(t,x)))-f (\underline u(t,x))
\l \fr{3}4f'(1) \tilde \vp e^{-\mu \left(x \cdot
e_{i}-L\right)} \zeta_i(x) \phi^\beta(\underline{\xi}_i(t,x)).
\end{align*}
Note that \eqref{defvp} and \eqref{deft} yield that $-c_f +\rho\vp e^{\vp t}<0$. Since $\phi'<0$, $\beta<1$, $\mu^2\leq\mu=\vp/c_f\l -f'(1)/8$, it then  follows from \eqref{defvp}, \eqref{defphi'}, \eqref{deft} and \eqref{zzzz} that
\begin{align*}
\mathscr L[\underline u](t,x)\l
& \fr34{f'(1)}  \tilde \vp  e^{-\mu \left(x \cdot
e_{i}-L\right)} \zeta_i(x) \phi^\beta(\underline{\xi}_i(t,x)) \\
&+\tilde \vp  e^{-\mu \left(x \cdot e_{i}-L\right) } \zeta_i(x) \phi^{\beta}(\underline{\xi}_i(t,x)) \e( \mu ^2 -\fr{f'(1)}8
+\beta \fr{\phi''(\underline{\xi}_i(t,x))}{\phi(\underline{\xi}_i(t,x))}\r)\\
\l& \tilde \vp e^{-\mu \left(x \cdot
e_{i}-L\right) }\zeta_i(x) \phi^{\beta}(\underline{\xi}_i(t,x))\e(\fr{3f'(1)}4-\fr{f'(1)}8-\fr{f'(1)}8-\fr{f'(1)}8 \r)\\
\l & 0.
\end{align*}

If
$\underline{\xi}_i(t,x)\g C$, then one gets from \eqref{defvp} and \eqref{defC'} that $\underline u (t,x)\leq\phi(\underline{\xi}_i(t,x))\l \vp'\l \theta$,  which implies that  $f ( \phi(\underline{\xi}_i(t,x)))=f (\underline u(t,x))=0$.
Since $\phi'<0$ and $\mu^2\leq\mu=\vp/c_f\l -\Lambda(-\beta c_f)/8$, one concludes from \eqref{defvp}, \eqref{betac1}, \eqref{betac2}, \eqref{deft} and \eqref{zzzz} that
\begin{align*}
\mathscr L[\underline u](t,x)
&\leq\tilde\varepsilon e^{-\mu(x\cdot e_i-L)}\zeta_i(x)\phi^\beta(\underline{\xi}_i(t,x))\left(
\mu^2-\rho\varepsilon\beta e^{\varepsilon t}\frac{\phi'(\underline{\xi}_i(t,x))}{\phi(\underline{\xi}_i(t,x))}
+c_f\beta\frac{\phi'(\underline{\xi}_i(t,x))}{\phi(\underline{\xi}_i(t,x))}\right.\\
&~~\left.+\beta^2\left(\frac{\phi'(\underline{\xi}_i(t,x))}{\phi(\underline{\xi}_i(t,x))}\right)^2
+\beta\left(\frac{\phi''(\underline{\xi}_i(t,x))}{\phi(\underline{\xi}_i(t,x))}-\left(\frac{\phi'(\underline{\xi}_i(t,x))}{\phi(\underline{\xi}_i(t,x))}\right)^2\right)
-\frac{\Lambda(-\beta c_f)}{8}\right)\\
&\leq\tilde\varepsilon e^{-\mu(x\cdot e_i-L)}\zeta_i(x)\phi^\beta(\underline{\xi}_i(t,x))\left(
-\frac{\Lambda(-\beta c_f)}{8}-\frac{\Lambda(-\beta c_f)}{8}\right.\\
&~~\left.+\frac{3\Lambda(-\beta c_f)}{4}-\frac{\Lambda(-\beta c_f)}{4}-\frac{\Lambda(-\beta c_f)}{8}\right)\\
&\leq0.
\end{align*}

If
$-C\l\underline{\xi}_i(t,x)\l C$, then $\phi'(\underline{\xi}_i(t,x))< -\kappa$ by \eqref{defk}.
From \eqref{deflf}, one obtains that
\[
\begin{aligned}
&~~~f ( \phi(\underline{\xi}_i(t,x)))-f (\underline u(t,x))
\l L_f \tilde \vp  e^{-\mu \left(x \cdot
e_{i}-L\right)} \zeta_i(x) \phi^\beta(\underline{\xi}_i(t,x)).
\end{aligned}
\]
Moreover, it can be deduced  from $\underline{\xi}_i(t,x)\l C$ that
 $x\cdot e_i\g
-c_ft {+1}+\rho e^{\vp t}-C\g -c_ft -C$, hence
\[e^{-\mu \left(x \cdot
e_{i}-L\right) }\l e^{-\mu \left(-c_f t -C-L\right) }= e^{ \mu  c_f t } e^{ \mu (C+L)}=  e^{ \vp t} e^{ \mu (C+L)}.
\]
Since $\beta<1$, $0<\phi<1$ and $\phi'<0$, it then follows from \eqref{sizeta}, \eqref{deflf}, \eqref{defvp}, \eqref{mu}, \eqref{defk}, \eqref{defrho} and \eqref{zzzz} that
\begin{align*}
\mathscr L[\underline u](t,x)
&\l \rho\vp e^{\vp t} \phi'(\underline{\xi}_i(t,x))+ \tilde \vp  \mu ^2 e^{-\mu \left(x \cdot
e_{i}-L\right) } \zeta_i(x) \phi^{\beta}(\underline{\xi}_i(t,x))\\
&~~~-2 \tilde \vp  \mu   e^{-\mu \left(x \cdot
e_{i}-L\right)  } \nabla\zeta_i(x)\cdot e_i  \phi^{\beta}(\underline{\xi}_i(t,x))
+ \tilde \vp    e^{-\mu \left(x \cdot
e_{i}-L\right) } \Delta\zeta_i(x) \phi^{\beta}(\underline{\xi}_i(t,x))\\
&~~~-2 \tilde \vp    \beta  e^{-\mu \left(x \cdot
e_{i}-L\right) } \nabla\zeta_i(x)\cdot e_i \phi^{\beta-1}(\underline{\xi}_i(t,x))\phi'(\underline{\xi}_i(t,x))\\
&~~~+\tilde \vp \beta  e^{-\mu \left(x \cdot
e_{i}-L\right) } \zeta_i(x)\phi^{\beta-1}(\underline{\xi}_i(t,x))\phi''(\underline{\xi}_i(t,x))
+f ( \phi(\underline{\xi}_i(t,x)))-f (\underline u(t,x))\\
&\l -\rho\vp e^{\vp t} \kappa+ \tilde\vp   e^{\vp t} e^{\mu (C+L)}\bar\zeta_0\e(\mu ^2 -f'(1)/8
+\beta  \e\| \phi''/\phi \r\|_{L^\i(\R)} +L_f\r)\\
&\l 0.
\end{align*}

As a result, the function $\underline u(t,x)$ is a subsolution of \eqref{te1.1} for $t\leq T$ and $x\in\Omega$.
This completes the proof.
\end{pr}
\vspace{0.3cm}

 {Let $\eta_\vp:\R\to\R$ be a nonincreasing $C^2$ function such that
\begin{align}\label{eet}
\eta_\vp(s)=
\begin{cases}
e^{-\mu(s-L)},&s\geq L+3,\\
1,&s\leq L+2,
\end{cases}
\end{align}
where $\mu=\vp/c_f$ is defined by \eqref{mu}.
Even if it means decreasing $\vp$, one can assume without loss of generality that
\begin{align}\label{etv}
\|\eta_\vp''/\eta_\vp\|_{L^\i(\R)}
\leq-\frac{\Lambda(-\beta c_f)}{4}.
\end{align}
Let $\pi_\vp:\R\to\R$ be a nondecreasing $C^2$ function such that
\begin{align}\label{pii}
\pi_\vp(s)=
\begin{cases}
1,&s\geq L+2,\\
0,&s\leq L+1.
\end{cases}
\end{align}
Since $\phi(-\i)=0$, even if it means increasing $C$ defined as in \eqref{defC'}, we can assume that}
\begin{align}\label{pv}
 {\phi^{1-\beta}(\xi)\left(\|\pi_\vp''\|_{L^\i(\R)}+2\|\pi_\vp'\|_{L^\i(\R)}\|\phi'/\phi\|_{L^\i(\R)}\right)\leq-\frac{\vp\Lambda(-\beta c_f)}{4}
\ \ \text{ for all}\ \xi\geq C.}
\end{align}

\begin{lem}\lb{exsup}
 {There exist
$\rho>0$ and $T_\vp<0$ such that the function $\bar u(t,x)$ defined by
\be\lb{baru}
\bar u(t,x)=\bc
\bar u_1(t,x) \ \ &\text{in }(-\infty,T_\varepsilon]\times\bar{\Omega\backslash\bigcup_{j\in J}\{ x\in{\mathcal H_j}:x\cdot e_j\g L\}},\\
\min(\bar u_1(t,x),\bar u_2(t,x))\ \ &\text{in }(-\infty,T_\varepsilon]\times\bigcup_{j\in J}\bar{\{x\in{\mathcal H_j}:L\l x\cdot e_j\l L+1\}},\\
\bar u_2(t,x) \ \ &\text{in }(-\infty,T_\varepsilon]\times\bigcup_{j\in J}\bar{\{ x\in{\mathcal H_j}:x\cdot e_j\g L+1\}}
\ec
\ee
is a supersolution of \eqref{te1.1} for $t\l T_\vp$ and
$x\in \Omega$, where
$$
\bar u_1(t,x):=
\begin{cases}
\pi_\vp(x\cdot e_i) \phi(\bar\xi_i(t,x))+\vp\eta_\vp(x\cdot e_i)\phi^\beta(\xi(t))\\
~~~~~~~~~~~~~~~~~~~~~~~~
\text{in $(-\infty,T_\varepsilon]\times \bar{\{x\in{\mathcal H_i}:x\cdot e_i\geq L+1\}}$  for each $i\in I$},
\\
\vp\phi^\beta(\xi(t))
~~~~~~~~~~~~~\text{in $(-\infty,T_\varepsilon]\times \bar{\Omega\backslash\bigcup_{s=1}^m\{x\in\mathcal H_s:x\cdot e_s\geq L+1\}}$}
\end{cases}
$$
and
$$
\bar u_2(t,x):=\vp e^{-\mu(x\cdot e_j-L)}\phi^\beta(\hat\xi_j(t,x))\ \  \text{ in } (-\infty,T_\varepsilon]\times \{x\in\bar{\mathcal H_j}:x\cdot e_j\geq L\}\ \text{ for each }j\in J
$$
with $\bar\xi_i(t,x)=-x\cdot e_i -c_f  t-\rho e^{\vp t}$, $\hat\xi_j(t,x)=x\cdot e_j -c_f  t-\rho e^{\vp t}-2L-3$ and
$\xi(t)=-L-c_ft-\rho e^{\vp t}-2$.}
\end{lem}

\begin{pr}
 {Take $\rho>0$ such that
\begin{align}\label{r-k}
\rho \kappa \geq e^{\mu(C+L+1)}\e(\frac{|f'(1)|}{8}+L_f\r),
\end{align}
where $L_f$ is defined by \eqref{deflf}.
Choose $T_\vp<0$ negative enough such that}
\be\lb{Tvp}
 {\rho e^{\vp T_\vp} \l 1
\ \ \text{and}\ \ c_f T_\vp\l-C-2L-4.}
\ee

 We first verify that $\bar u(t,x)$ is well defined. Note that the surfaces $\{x\in\mathcal H_j:x\cdot e_j=L\}$ and
 $\{x\in\mathcal H_j:x\cdot e_j=L+1\}$ are bounded for each $j\in J$
 since the section $\omega_j\subset\R^{N-1}$
given by \eqref{defhi1} is  bounded. For $x\in\bar{\mathcal H_i}$ such that $x\cdot e_i=L+1$ with $i\in I$,
it can be deduced from \eqref{eet} and \eqref{pii} that $\pi_\vp(x\cdot e_i)=0$ and $\eta_\vp(x\cdot e_i)=1$, hence
$\bar u_1(t,x)=\vp\phi^\beta(\xi(t))$.
It is evident that $\bar u_1(t,x)$ is of $C^{1,2}$ in
$(-\i,T_\vp]\times\e(\Omega\backslash\bigcup_{j\in J} \{x\in\mathcal H_j:x\cdot e_j\geq L+1\}\r)$
and  that $\bar u_2(t,x)$ is of $C^{1,2}$ in
$(-\i,T_\vp]\times\e(\bigcup_{j\in J} \{x\in\mathcal H_j:x\cdot e_j\geq L\}\r)$.
 Take any $j\in J$.
For $t\l T_\vp$ and
$x\in \bar{\mathcal H_j}$ such that $x\cdot e_j=L$,
it follows from $\phi'<0$ that
$$
\bar u_2(t,x)\g\vp\phi^\beta(L-c_f t-\rho e^{\vp t}-2L-3)\g\vp\phi^\beta(-L-c_f t-\rho e^{\vp t}-2)=
\vp \phi^\beta(\xi(t))=\bar u_1(t,x).
$$
For $t\l T_\vp$ and
$x\in \bar{\mathcal H_j}$ such that $x\cdot e_j=L+1$, one gets from $\phi'<0$ that
$$
\bar u_2(t,x)\l\vp\phi^{\beta}(L+1-c_ft-\rho e^{\vp t}-2L-3)= \vp \phi^\beta(\xi(t))=\bar u_1(t,x).
$$
As a conclusion, the function $\bar{u}$ is  continuous in $(-\i,T_\vp] \times \bar{\Omega}$.

 We then prove that the function $\bar u$ defined by \eqref{baru} is a supersolution of \eqref{te1.1}
in $(-\i,T_\vp]\times\bar\Omega$.
Let us first check the boundary condition.
It is easy to see that $(\bar u_1)_\nu(t,x)=0$  for $t\leq T_\varepsilon$ and $x\in\p{\Omega}\cap\p\bar{\e(\Omega\backslash \bigcup_{s=1 }^m\{x\in\mathcal H_s:x\cdot e_s\geq L+1 \}\r)}$.
Since $ \nu\cdot e_s=0$ on $\partial\mathcal H_s\cap\p\Omega$ for all $s\in\{1,\cdots,m\}$,
it then follows that
\begin{align*}
(\bar u_1)_\nu(t,x)=&\pi_\vp'(x\cdot e_i) \phi(\bar\xi_i(t,x))\nu\cdot e_i-\pi_\vp(x\cdot e_i) \phi'(\bar\xi_i(t,x))\nu\cdot e_i
+\vp\eta_\vp'(x\cdot e_i)\phi^\beta(\xi(t))\nu\cdot e_i=0
\end{align*}
for $t\leq T_\varepsilon$ and $x\in\p{\mathcal H_i}\cap\p\Omega$ such that $x\cdot e_i\geq L+1$ with any $i\in I$,
and
\begin{align*}
(\bar u_2)_\nu(t,x)=-\vp \mu e^{-\mu(x\cdot e_j-L)}\phi^\beta(\hat\xi_j(t,x))\nu\cdot e_j
+\vp e^{-\mu(x\cdot e_j-L)}\beta\phi^{\beta-1}(\hat\xi_j(t,x))\phi'(\hat\xi_j(t,x))\nu\cdot e_j
= 0
\end{align*}
for $t\l T_\vp$ and $ x\in \p{\mathcal H_j}\cap\p\Omega$ such that $x\cdot e_j\geq L$ with any $j\in J$.

We now turn to  prove that
$\mathscr L[\bar u_1](t,x)=(\bar u_{1})_t(t,x)- \Delta\bar u_{1}(t,x)-f(\bar u_1(t,x))\g 0$
for $t\l T_\vp$ and
$ x\in \bigcup_{i\in I}\bar{\{x\in\mathcal H_i:x\cdot e_i\geq L+1\}}\cup\bar{\Omega\backslash\bigcup_{s=1}^m\{x\in\mathcal H_s:x\cdot e_s\geq L+1\} }$.
Fix any $i\in I$.
For $t\l T_\vp$ and $ x\in \bar{\mathcal H_i}\cap\overline\Omega$ such that $x\cdot e_i\geq L+3$,
one infers from \eqref{eet} and \eqref{pii} that $\pi_\vp(x\cdot e_i)=1$ and $\eta_\vp(x\cdot e_i)= e^{-\mu(x\cdot e_i-L)}$,
hence $\bar u_1(t,x)=\phi(\bar\xi_i(t,x))+\vp e^{-\mu(x\cdot e_i-L)}\phi^\beta(\xi(t))$.
After some computations, it follows that
 \begin{align*}
\mathscr L [\bar u_1](t,x)
=&-\rho\vp  e^{\vp  t}  \phi'(\bar\xi_i(t,x))
-\vp  c_f \beta e^{-\mu \left(x \cdot
e_{i}-L\right)} \phi^{\beta-1}(\xi(t)) \phi'(\xi(t))\\
&-\vp  \beta \rho\vp e^{\vp t}e^{-\mu \left(x \cdot
e_{i}-L\right)}  \phi^{\beta-1}(\xi(t)) \phi'(\xi(t))
-  \vp  \mu ^2 e^{-\mu \left(x \cdot
e_{i}-L\right)}  \phi^{\beta}(\xi(t)) \\
&+f ( \phi(\bar\xi_i(t,x)))-f (\bar u_1(t,x)).
\end{align*}

If
$\bar \xi_i(t,x)\l -C$,
one then gets from $\phi<1$, \eqref{defvp} and \eqref{defC'} that
$$
1-\vp_0 \l1-\vp\l\phi(\bar\xi_i(t,x))\l \bar u_1(t,x) \l  1+\vp\leq1+\vp_0
.$$
Together with \eqref{vp0} and the mean value theorem, there holds
\begin{align*}
&f ( \phi(\bar\xi_i(t,x)))-f (\bar u_1(t,x))
\g -\fr{3 f'(1)}4\vp e^{-\mu(x\cdot e_i-L)}\phi^\beta(\xi(t)).
\end{align*}
Therefore, it can be deduced from $\phi'<0$, \eqref{defvp} and \eqref{mu}   that
\begin{align*}
\mathscr L [\bar u_1](t,x)
&\g \vp e^{-\mu \left(x \cdot
e_{i}-L\right)} \phi^{\beta}(\xi(t)) \e(-  \mu ^2
-\frac{3f'(1)}{4}\r) \\
&\g\vp  e^{-\mu(x\cdot e_i-L) }  \phi^\beta(\xi(t))
\e(-\fr{3f'(1)}4 +\fr{f'(1)}8\r)\\
&\g  0.
\end{align*}

If $\bar\xi_i(t,x)\g C$,  then
$ 0\leq\phi(\bar\xi_i(t,x))\leq \bar u_1 (t,x)\l 2\vp\l\theta $ by  $x\cdot e_i\geq L+3$, \eqref{defvp} and \eqref{defC'},
hence $f ( \phi(\bar\xi_i(t,x)))=0=f (\bar u_1(t,x))$. By \eqref{Tvp}, there holds $\xi(t)\geq-L-c_fT_\vp-\rho e^{\vp T_\vp}-2\geq C$.
Thus, it can be inferred from $\phi'<0$, \eqref{defvp}, \eqref{mu} and \eqref{betac2} that
\begin{align*}
\mathscr L [\bar u_1](t,x)
&\g  \vp e^{-\mu \left(x \cdot
e_{i}-L\right)} \phi^{\beta}(\xi(t)) \e(-  \mu ^2
-c_f\beta\frac{\phi' (\xi(t))}{\phi(\xi(t))}\r) \\
&\g \vp e^{-\mu(x\cdot e_i-L) }  \phi^\beta(\xi(t))
\e(-\fr{3\Lambda(-\beta c_f)}4 +\fr{\Lambda(-\beta c_f)}8\r)\\
&\g  0.
\end{align*}

If $-C\l \bar\xi_i(t,x)\l C$, then
$x\cdot e_i \g -c_f  t-C-1$ by \eqref{Tvp}, hence
$$
e^{-\mu(x\cdot e_i-L)}\l e^{\mu c_f  t}e^{\mu(C+L+1)}
= e^{\vp t}e^{\mu(C+L+1)}.
$$
Since $0\leq\bar u_1(t,x)\leq1+\vp\leq1+\vp_0$ by \eqref{defvp},
  one derives from $\phi<1$, \eqref{sizeta} and \eqref{deflf}  that
\begin{align*}
f ( \phi(\bar\xi_i(t,x)))-f (\bar u_1(t,x))
\g -\vp L_fe^{-\mu(x\cdot e_i-L)} \phi^\beta(\xi(t))\geq-\vp L_fe^{\vp t}e^{\mu(C+L+1)}  .
\end{align*}
By $\phi'<0$, \eqref{mu}, \eqref{defk} and \eqref{r-k}, one gets that
\begin{align*}
\mathscr L [\bar u_1](t,x)
&\geq\rho\vp  \kappa e^{\vp  t}
- \vp  \mu ^2 e^{-\mu \left(x \cdot
e_{i}-L\right)}  \phi^{\beta}(\xi(t))+f ( \phi(\bar\xi_i(t,x)))-f (\bar u_1(t,x))\\
&\geq\rho\vp  \kappa e^{\vp  t}
-\vp e^{\vp t}e^{\mu(C+L+1)}\e(\frac{|f'(1)|}{8}+L_f\r)\\
&\geq0.
\end{align*}

For $t\l T_\vp$ and $ x\in \bar{\mathcal H_i}\cap\overline\Omega$ such that $L+2\leq x\cdot e_i\leq L+3$,
one infers from \eqref{pii} that $\pi_\vp(x\cdot e_i)=1$, hence
$\bar u_1(t,x)=\phi(\bar\xi_i(t,x))+\vp\eta_\vp(x\cdot e_i)\phi^\beta(\xi(t))$.
It follows from \eqref{Tvp} that $\bar\xi_i(t,x)\geq-L-3-c_f T_\vp-\rho e^{\vp T_\vp}\geq C$,
hence $0\leq\phi(\bar\xi_i(t,x))\leq \bar u_1(t,x)\leq2\vp\leq\theta$ by \eqref{defvp},
which implies that $f(\phi(\bar\xi_i(t,x)))=f(\bar u_1(t,x))$=0.
Note that \eqref{Tvp} yields that $\xi(t)\geq -L-c_fT_\vp-\rho e^{\vp T_\vp}-2\geq C$.
After some direct calculations,
one gets from $\phi'<0$, \eqref{betac2} and \eqref{etv} that
\begin{align*}
\mathscr L [\bar u_1](t,x)
=&-\rho\vp  e^{\vp  t}  \phi'(\bar\xi_i(t,x))
-\vp c_f \beta \eta_\vp(x\cdot e_i) \phi^{\beta-1}(\xi(t)) \phi'(\xi(t))\\
&-\vp \beta \rho\vp e^{\vp t}\eta_\vp(x\cdot e_i) \phi^{\beta-1}(\xi(t)) \phi'(\xi(t))
-  \vp\eta_\vp''(x\cdot e_i) \phi^{\beta}(\xi(t)) \\
\geq&\vp \eta_\vp(x\cdot e_i)\phi^\beta(\xi(t))\left(-c_f\beta \frac{\phi'(\xi(t))}{ \phi(\xi(t))}-\frac{\eta_\vp''(x\cdot e_i)}{\eta_\vp(x\cdot e_i)}\right)\\
\geq&\vp \eta_\vp(x\cdot e_i)\phi^\beta(\xi(t))\left(-\frac{3\Lambda(-\beta c_f)}{4}-\left\|{\eta_\vp''}/{\eta_\vp}\right\|_{L^\infty(\R)}\right)\\
\geq&0.
\end{align*}

For $t\l T_\vp$ and $ x\in \bar{\mathcal H_i}\cap\overline\Omega$ such that $L+1\leq x\cdot e_i\leq L+2$,
it can be deduced from \eqref{eet} that $\eta_\vp(x\cdot e_i)=1$, hence
$\bar u_1(t,x)=\pi_\vp(x\cdot e_i)\phi(\bar\xi_i(t,x))+\vp\phi^\beta(\xi(t))$.
By \eqref{Tvp}, one gets that $\bar\xi_i(t,x)\geq-L-2 -c_f T_\vp-\rho e^{\vp T_\vp}\geq C$,
it then can be deduced from $\pi_\vp\in(0,1)$, \eqref{defvp} and \eqref{defC'} that
 $0\leq\phi(\bar\xi_i(t,x))\leq\bar u_1(t,x)\leq2\vp\leq\theta$,
which implies that $f(\phi(\bar\xi_i(t,x)))=f(\bar u_1(t,x))$=0.
Moreover, since $\phi'<0$  and $\bar\xi_i(t,x)\geq\xi(t)=-L-2-c_ft-\rho e^{\vp t}\geq C$ by \eqref{Tvp}, one then arrives at
$\phi(\bar\xi_i(t,x))\leq\phi(\xi(t))$.
After some computations, it follows from $\phi'<0$, $\beta\in(0,1)$, \eqref{betac2} and \eqref{pv} that
\begin{align*}
\mathscr L [\bar u_1](t,x)
=&-\rho\vp  e^{\vp  t}\pi_\vp(x\cdot e_i)  \phi'(\bar\xi_i(t,x))
-\vp  c_f \beta \phi^{\beta-1}(\xi(t)) \phi'(\xi(t))\\
&+\pi_\vp(x\cdot e_i)f(\phi(\bar\xi_i(t,x)))
-f(\bar u_1(t,x))-\vp  \beta \rho\vp e^{\vp t}  \phi^{\beta-1}(\xi(t)) \phi'(\xi(t))\\
&-  \pi_\vp''(x\cdot e_i)\phi(\bar\xi_i(t,x)) +2\pi_\vp'(x\cdot e_i)\phi'(\bar\xi_i(t,x)) \\
\geq&-\phi^\beta(\xi(t))\left(\phi^{1-\beta}(\bar\xi_i(t,x))
\left(\|\pi_\vp''\|_{L^\i(\R)}
+2\|\phi'/\phi\|_{L^\i(\R)}\|\pi_\vp'\|_{L^\i(\R)}\right)+\vp \beta c_f\frac{\phi'(\xi(t))}{\phi(\xi(t))}\right)\\
\geq&-\phi^\beta(\xi(t))\left(-\frac{\vp\Lambda(-\beta c_f)}{4} +\frac{3\vp\Lambda(-\beta c_f)}{4} \right)\\
\geq&0.
\end{align*}

For $t\l T_\vp$ and $ x\in \bar{\Omega\backslash\bigcup_{s=1}^m\{x\in\mathcal H_s:x\cdot e_s\geq L+1\}}$,
we know that $\bar u_1(t,x)=\vp\phi^\beta(\xi(t))$.
By \eqref{defvp} and $\phi<1$, there holds
$\bar u_1(t,x)\leq\vp\leq\theta$, then $f(\bar u_1(t,x))=0$.
Together with $\phi'<0$, one gets that
\begin{align*}
\mathscr L [\bar u_1](t,x)=\vp\beta(-c_f-\rho\vp e^{\vp t})\phi^{\beta-1}(\xi(t))\phi'(\xi(t))-f(\bar u_1(t,x))\geq0.
\end{align*}

Fix any $j\in J$.  Let us now check that
$\mathscr L[\bar u_2](t,x)=(\bar u_{2})_t(t,x)- \Delta \bar u_{2}(t,x)-f(\bar u_2(t,x))\g 0$
for $t\l T_\vp$ and $ x\in \bar{\mathcal H_j}\cap\overline\Omega$ such that $x\cdot e_j\geq L$. After some computations, it  follows that
\begin{align*}
\mathscr L [\bar u_2](t,x)
=&-\vp(c_f+\rho\vp e^{\vp t}) \beta e^{-\mu \left(x \cdot e_{j}-L\right)}  \phi^{\beta-1}(\hat\xi_j(t,x)) \phi'(\hat\xi_j(t,x))\\
&-\vp  \mu ^2 e^{-\mu \left(x \cdot e_{j}-L\right)}  \phi^{\beta}(\hat\xi_j(t,x))
    +2\vp \mu \beta  e^{-\mu \left(x \cdot e_{j}-L\right)}     \phi^{\beta-1} (\hat\xi_j(t,x))\phi' (\hat\xi_j(t,x))\\
&-\vp \beta(\beta-1)   e^{-\mu \left(x \cdot e_{j}-L\right)}  \phi^{\beta-2} (\hat\xi_j(t,x))(\phi' (\hat\xi_j(t,x)))^2\\
&-\vp \beta  e^{-\mu \left(x \cdot e_{j}-L\right)}   \phi^{\beta-1} (\hat\xi_j(t,x))\phi'' (\hat\xi_j(t,x))-f (\bar u_2(t,x)).
\end{align*}
By \eqref{Tvp} and $x\cdot e_j\geq0$, there holds $\hat\xi_j(t,x)\geq-c_f T_\vp-\rho e^{\vp T_\vp}-2L-3\geq C$.
Since \eqref{sizeta}
and \eqref{defvp} imply that $\bar u_2(t,x)\leq\vp\leq\theta$, then $f(\bar u_2(t,x))=0$.
  Thus, it follows from $\phi'<0$, $\beta<1$,  \eqref{defvp}, \eqref{mu}, \eqref{betac1} and \eqref{betac2}  that
\begin{align*}
\mathscr L [\bar u_2](t,x)
\g&
\vp  e^{-\mu \left(x \cdot
e_{j}-L\right)}  \phi^{\beta}(\hat\xi_j(t,x))\e(- \mu ^2
 + \beta \e(\fr{\phi' (\hat\xi_j(t,x))}{\phi(\hat\xi_j(t,x))}\r)^2
-\beta\fr{\phi''(\hat\xi_j(t,x))}{\phi(\hat\xi_j(t,x))}\r.\\
 &\e.+(2\mu\beta-c_f\beta) \fr{\phi' (\hat\xi_j(t,x))}{\phi(\hat\xi_j(t,x))} -\beta^2 \e(\fr{\phi' (\hat\xi_j(t,x))}{\phi(\hat\xi_j(t,x))}\r)^2 \r) \\
\g &\vp e^{-\mu \left(x \cdot
e_{j}-L\right)}    \phi^\beta (\hat\xi_j(t,x))
\e(\fr{\Lambda(-\beta c_f)}{8}-\fr{3\Lambda(-\beta c_f)}{4}+\fr{\Lambda(-\beta c_f)}{8}+\fr{\Lambda(-\beta c_f)}{4}\r)\\
\g& 0.
\end{align*}

 As a result, the function $\bar u(t,x)$ defined by \eqref{baru} is a supersolution of \eqref{te1.1}
for $t\leq T_\vp$ and $x\in\Omega$. This completes the proof.
\end{pr}

\subsection{Existence, monotonity and uniqueness of the entire solution}
This subsection is devoted to existence, monotonicity and uniqueness of the entire solutions emanating
from planar traveling fronts in some branches, that is, we prove Theorem \ref{existence-}.
\vspace{0.3cm}

\begin{pr}[Proof of Theorem \ref{existence-}]
We divide the proof into two steps.

{\it Step 1: existence of entire solutions.}
Let us first consider the existence of entire solutions
emanating from planar traveling fronts in some branches.
Write $\tilde T=\min(T, T_\vp)$,  where $T$ and $T_{\varepsilon}$ are defined by Lemmas \ref{exsub} and \ref{exsup}, respectively.
For $t\l \tilde T$ and $x\in\overline{\mathcal H_i}\cap\overline\Omega$ such that $x\cdot e_i\geq L+2$ with any $i\in I$,
 it is easy to see from $\phi'<0$, \eqref{eet} and \eqref{pii} that
 $$
 \underline{u}(t,x)\leq\phi(\underline \xi_i(t,x))\leq\phi(\overline \xi_i(t,x))
\leq \bar u(t,x).
$$
For $t\l\tilde T$ and $x\in\overline{\mathcal H_i}\cap\overline\Omega$ such that $L\leq x\cdot e_i\leq L+2$ with any $i\in I$,
one infers from \eqref{Tvp} that $\underline\xi_i(t,x)\geq\xi(t)=-L-2-c_f T_\vp-\rho e^{\vp T_\vp}\geq C$,
it then follows from $\phi'<0$, \eqref{defC'} and \eqref{eet} that
$$
\underline u(t,x)\leq \phi^{1-\beta}(\underline\xi_i(t,x))\phi^\beta(\underline\xi_i(t,x))
\leq\vp\phi^\beta(\xi(t))\leq\bar u(t,x).
$$
For $t\leq\tilde T$ and $x\in\bar\Omega\backslash \bigcup_{i\in I}\{x\in\mathcal H_i:x\cdot e_i\geq L\}$, one infers from the proof of Lemma \ref{exsub} that
$\underline u(t,x)=0\leq\bar u(t,x)$.
Hence, one obtains that $\underline{u}(t,x)
\leq \bar u(t,x)$ for $t\leq\tilde T$ and $x\in\bar\Omega$.

Take $n\in\mathbb N\cap(|\tilde T|+1,+\i)$.  Let $u_{n}(t,x)$ be the solution of  \eqref{te1.1} for $t\geq-n$ with initial value
$u_{n}(-n, x)=\underline{u}(-n, x)$.
By Lemma \ref{cp}, there holds
\be\lb{ex1.1}
\underline u(t, x) \leq u_{n}(t, x) \leq \bar u(t, x) \ \text { for }-n \leq t\l \tilde T  \text { and } x \in\Omega
\ee
and
$$
u_{n}(t, x) \geq u_{n-1}(t, x) \ \text { for  }\ t \in[-n+1,+\infty) \ \text{and}\ x \in\Omega.
$$
It is easy to see that the sequence $u_{n}$  is nondecreasing in $n$.
Together with interior Schauder estimates, there is an entire solution
$u(t, x)$ to \eqref{te1.1} such that $u_{n}(t,x)\to  u(t,x)$ locally uniformly for $(t,x)\in\mathbb R\times\bar\Omega$ as $n\to+\i$. By \eqref{ex1.1}, the solution $ u(t, x)$ satisfies
\begin{align}\label{uesi}
\underline u(t, x) \leq  u(t, x) \leq \bar u(t, x)\  \text { for all } t \leq \tilde T \text { and } x \in \overline\Omega.
\end{align}

Let us check that the solution $u$ satisfies \eqref{dm0}. We first show that for any $i\in I$, there holds
\begin{align}\label{limI1}
I_1(t,x):= \fr{u(t, x)- \phi(-x\cdot e_i-c_f t)}{\phi^\beta(-x\cdot e_i-c_f t)}\to0\text{\ \  uniformly for $x\in\bar{\mathcal H_i}\cap\bar\Omega$ as $t\to-\i$.}
\end{align}
 For $x\in\bar{\mathcal H_i}\cap\bar\Omega$ such that $x\cdot e_i\geq L+1$,
 since $\pi(x\cdot e_i)\in[0,1]$ and $\eta(x\cdot e_i)\in[0,1]$ by \eqref{eet} and \eqref{pii}, then
 $$
 \bar u(t,x)=\bar u_1(t,x)\leq \phi(-x\cdot e_i-c_f t-\rho e^{\vp t})+\vp\phi^\beta(-L-c_f t-\rho e^{\vp t}-2).
 $$
 Since $\phi<\phi^\beta$ and $\phi'<0$, it follows from \eqref{sizeta}, \eqref{mu},  \eqref{baru}, \eqref{uesi} and the mean value theorem that there exists $\tau\in(0,1)$ such that
\begin{align*}
I_1(t,x)
&\leq\fr{\bar u_1(t, x)- \phi(-x\cdot e_i-c_f t)}{\phi^\beta(-x\cdot e_i-c_f t)}\nonumber\\
&\leq\frac{-\rho e^{\vp t}\phi'(-x\cdot e_i-c_f t-\tau\rho e^{\vp t})+\vp \phi^\beta(-L-c_f t-\rho e^{\vp t}-2)}
{\phi^\beta(-x\cdot e_i-c_f t)}\nonumber\\
&\leq\rho e^{\vp t}\|\phi'/\phi\|_{L^\i(\R)}\frac{\phi^\beta(-x\cdot e_i-c_f t-\tau\rho e^{\vp t})}{{\phi^\beta(-x\cdot e_i-c_f t)}}
+\frac{\vp\phi^\beta(-x\cdot e_i-1-c_f t-\rho e^{\vp t})}
{\phi^\beta(-x\cdot e_i-c_f t)}\nonumber\\
&\leq\left(\rho e^{\vp t}\|\phi'/\phi\|_{L^\i(\R)}+\vp\right)\frac{\phi^\beta(-x\cdot e_i-1-c_f t-\rho e^{\vp t})}{{\phi^\beta(-x\cdot e_i-c_f t)}}.
\end{align*}
Write $H(t,x):={\phi^\beta(-x\cdot e_i-1-c_f t-\rho e^{\vp t})}/{{\phi^\beta(-x\cdot e_i-c_f t)}}$. If $-x\cdot e_i-c_f t\to-\i$ as $t\to-\i$, then $\lim_{t\to-\i}H(t,x)=1$ since $\phi(-\i)=1$. If $-x\cdot e_i-c_f t\to+\i$ as $t\to-\i$, it then follows from \eqref{estimate2} that
$\lim_{t\to-\i}H(t,x)
\leq{K_4^\beta}e^{c_f}/{K_3^\beta}$.
If there exists a constant $R>0$ such that $|-x\cdot e_i-c_f t|\leq R$ as $t\to-\i$, then  $\lim_{t\to-\i}H(t,x)\leq{1}/{\phi^\beta(R)}$.
Thus, there exists a constant $C_1>0$ such that
\begin{align}\label{I11}
\lim_{t\to-\i}I_1(t,x)\leq C_1\vp\text{ \ \ uniformly for $x\in\bar{\mathcal H_i}\cap\bar\Omega$ such that $x\cdot e_i\geq L+2$.}
\end{align}
For $x\in\bar{\mathcal H_i}\cap\bar\Omega$ such that $0\leq x\cdot e_i\leq L+1$,
there holds
$$
u(t,x)\leq\min\e(\bar u_1(t,x),\bar u_2(t,x)\r)\leq\bar u_1(t,x)=\vp\phi^\beta(-L-c_f t-\rho e^{\vp t}-2).
$$
It follows from $\phi>0$, \eqref{estimate2} and \eqref{baru} that
\begin{align}\label{I12}
\lim_{t\to-\i}I_1(t,x)
\leq&\lim_{t\to-\i}\fr{\bar u_1(t, x)}{\phi^\beta(-x\cdot e_i-c_f t)}
\leq\lim_{t\to-\i}\frac{\vp\phi^\beta(-L-c_f t-\rho e^{\vp t}-2)}
{\phi^\beta(-c_f t)}\leq\frac{K_4^\beta}{K_3^\beta}\vp e^{\beta c_f(L+2)}.
\end{align}
On the other hand,
it can be derived from a similar argument as the proof of \eqref{I11} that there exists a constant $C_2>0$ such that
\begin{align}\label{I14}
I_1(t,x)
\geq\fr{\underline u(t, x)- \phi(-x\cdot e_i-c_f t)}{\phi^\beta(-x\cdot e_i-c_f t)}
\geq -C_2\vp\text{ \ \ uniformly for $x\in\bar{\mathcal H_i}\cap\bar\Omega$ as $t\to+\i$.}
\end{align}
Combining with \eqref{I11}, \eqref{I12},  \eqref{I14} and the arbitrariness of ${\vp}$, one obtains that \eqref{limI1} is true.

We then  prove that for any $j\in J$,
\begin{align}\label{uuj}
I_2(t,x):=\frac{u(t, x)}{\phi^\beta(x\cdot e_j-c_f t)} \to   0 \text{  uniformly in }\overline{\mathcal{H}_{j}} \cap \bar{\Omega}\text{ as }t\to-\i.
\end{align}
For $x\in\bar {\mathcal H_j}\cap\bar\Omega$ such that $x\cdot e_j\geq L$, it follows from \eqref{estimate2} that
\begin{align*}
0\leq\lim_{t\to-\i}I_2(t,x)
&\leq\lim_{t\to-\i}\frac{\bar u_2(t,x)}{\phi^\beta(x\cdot e_j-c_f t)}\\
&\leq\lim_{t\to-\i}\frac{\vp \phi^\beta(x\cdot e_j-c_ft-\rho e^{\vp t}-2L-3) }{\phi^\beta(x\cdot e_j-c_f t)}
\leq\frac{K_4^\beta}{K_3^\beta}\vp e^{c_f\beta(2L+3)}.
\end{align*}
For $x\in\bar {\mathcal H_j}\cap\bar\Omega$ such that $0\leq x\cdot e_j\leq L$, it follows from \eqref{estimate2} and $\phi'<0$ that
\begin{align*}
0\leq\lim_{t\to-\i}I_2(t,x)
&\leq\lim_{t\to-\i}\frac{\bar u_1(t,x)}{\phi^\beta(x\cdot e_j-c_f t)}\\
&\leq\lim_{t\to-\i}\frac{\vp \phi^\beta(-L-c_ft-\rho e^{\vp t}-2) }{\phi^\beta(L-c_f t)}
\leq\frac{K_4^\beta}{K_3^\beta}\vp e^{c_f\beta(2L+2)}.
\end{align*}
Together with the arbitrary of $\vp$, we arrive at \eqref{uuj}.

Since $\bar{\Omega\backslash\bigcup_{s=1}^m \mathcal H_s}\subset \bar {B(0,L)}$, then
some similar arguments  yield that
${u(t, x)}/{\phi^\beta(x\cdot e_j-c_f t)}\to0$ uniformly for $x\in\bar{\Omega\backslash\bigcup_{s=1}^m \mathcal H_s}$ for all $j\in J$ as $t\to-\i$.
As a result, we have proved that \eqref{te1.1} admits an entire solution $u$ satisfying \eqref{dm0}.

We now turn to prove that the solution $u$ is increasing in $t$. Since $\underline u_{t}\geq0$ for sufficiently negative $t$, which implies that $(u_n)_t(-n,x)\geq0$ for all $x\in\overline\Omega$ and large $n$. It then  follows from the comparison principle that
$\left(u_{n}\right)_{t}\geq0$ for $t>-n$ and $x \in \overline\Omega$. Passing   $n
\rightarrow+\infty$, one gets that $u_{t} \geq 0$ for $t \in \mathbb{R}$ and $x \in
\overline\Omega$. Again by the maximum principle and the Hopf boundary lemma, either $u_{t}>0$ or $u_{t} \equiv 0$ in $\mathbb R\times\overline\Omega$. Since
$ u(t, x) \rightarrow  \phi(-x\cdot e_i-c_f t)$ in $\overline{\mathcal{H}_{i}} \cap
\bar{\Omega}$ as $t \rightarrow-\infty$ for all $i\in I$ and $\phi'>0$, hence $u_{t} \equiv  0$
is impossible. Therefore, one has $u_{t}>0$ for $t \in \mathbb{R}$ and $x \in\overline \Omega$.

{\it Step 2: uniqueness of the entire solution.}
 Take any
 \begin{align}\label{ddddd}
 0<\delta \leq \min\left(\frac{\varepsilon_0}{2},\frac{\theta}{2},\frac{-f'(1)}{4},-\frac{\Lambda(-\beta c_f)}{4}\right).
 \end{align}
 By  \eqref{dm0}, there are $T_{\delta}<0$ small enough and $C_\delta>C$ large enough such that
\begin{align}\label{cd1}
&u(t,x)<\delta\text{ for }t\l T_\delta\text{ and }x\in\overline\Omega\cap\overline{\mathcal H_i}\text{ with }-x\cdot e_i-c_f t>C_\delta \text{ for any }i\in I\text{ or } x\in\overline{\Omega\backslash \bigcup_{i\in I}\mathcal H_i}
\end{align}
and
\begin{align}\label{cd2}
&u(t,x)>1-\delta\text{ for }t\l T_\delta\text{ and }x\in\overline\Omega\cap\overline{\mathcal H_i} \text{ with }
-x\cdot e_i-c_f t<-C_\delta\text{ for any }i\in I,
\end{align}
where the constant $C>0$ such that \eqref{betac2} and \eqref{defphi'} hold.  Even if it means decreasing $T_\delta$, one can assume that
\begin{align}\label{tdelta}
T_\delta\leq-\fr{2L+C}{c_f}.
\end{align}
Denote
\be\lb{ex1.2}
\Omega_\delta^i(t)=\left\{x \in \overline\Omega\cap\overline{\mathcal H_i}:\left|-x\cdot e_i-c_f  t\right| \leq C_\delta\right\}
\text{ and }\Omega_{\delta}(t)= \bigcup_{i \in I}\Omega_\delta^i(t).
\ee
We then claim that there is  $\hat\kappa>0$
such that
\be\lb{ex1.3}
u_t(t,x)\g \hat\kappa\ \ \text{for all $t\l T_\dl$ and $x\in\Omega_\dl(t)$.}
\ee
Assume to the contrary that there exists $i\in I$ such that
$$
\inf_{(t,x)\in(-\i, T_\delta]\times\Omega_\delta^i(t)}
u_t(t,x)=0.
$$
 Up to rotations, we can denote  $\mathcal H_i=\{(x_1,x'):x'\in \omega,\ x_1> 0\}$,
 where $\omega\subset\R^{N-1}$ is the open bounded connected section given by \eqref{defhi1}.
Then there exist two sequences $t_{n}\in\left(-\infty, T_{\dl}\right]$ and $x_{n}=\left(x_{n, 1}, x_{n}'\right) \in\Omega_{\dl}^i(t_n)$ such that
$$
u_{t} \left(t_{n}, x_{n}\right) \rightarrow   0  \text { as } n\rightarrow +\infty.
$$
Up to extraction of a subsequence,
we may assume that either $t_{n}$ converges to some $t^{*} \in$ $\left(-\infty,
T_{\dl}\right]$ or $t_{n} \rightarrow-\infty$ as $n \rightarrow \infty$. Since $x_n'\in\omega$ is bounded, one can assume that $x_n'\to x_*'$ as $n\to\i$. If
$t_n\to t_*$ as $n\to +\i$, it then follows from \eqref{ex1.2} that $x_{n, 1}$ is also bounded for large $n$.
Up to extraction of a subsequence, we can assume that $x_{n, 1} \rightarrow x_{1}^{*}$ as
$n \rightarrow +\infty$. Then $u_{t} \left(t^*,
x^*\right) = 0$ with $x^*=(x_1^*,x_*')$, this contradicts $u_t>0$ in $\R\times\bar{\Omega}$.

If $t_{n} \rightarrow-\infty$ as $n \rightarrow +\infty$,
it can be inferred from \eqref{ex1.2} that $x_{n,1}\to+\i$ as $n\to+\i$. Let $\mathcal H_i^n=\mathcal H_i-(x_{n,1},0)$,
then there is a cylinder $\mathcal H_i^\i$
such that $\mathcal H_i^n\to \mathcal H_i^\i$ as $n\to+\infty$.
Denote $u_n(t,x_1,x')= u(t+t_n,x_1+x_{n,1},x')$ for $t\in\R$ and $(x_1,x')\in\mathcal H_i^n$. Since $(0,x_n')\in
\mathcal H_i^n$, there holds $(u_n)_t(0,0,x_n')\to   0$ as $n\to+\i$.
By the parabolic estimates and up to extraction of a subsequence, there is a solution $u_\i(t,x)$ of \eqref{te1.1} such that
$u_n(t,x)\to u_\i(t,x)$ locally uniformly for $t\in\R$ and $x\in\mathcal H_i^\i$ as $n\to+\infty$. One immediately obtains that
$(u_\i)_t(0,0,x_*')= 0$  and $(u_\i)_t(t,x)\g   0$ for $t\in\R$ and $x\in\mathcal H_i^\i$ since $u_t>0$. By applying the maximum principle to $(u_\infty)_t$, one has
$(u_\i)_t(t,x)\equiv   0$ for $t\l 0$ and $x\in\mathcal H_i^\i$, this contradicts
the fact that $u$ satisfies \eqref{dm0}. The proof of the claim \eqref{ex1.3} is thereby complete.

Let us turn to prove the uniqueness of the entire solution.
Assume that there is another entire solution $v(t, x)$ satisfying \eqref{dm0}. Then, for any $0<\varepsilon<\delta/2$,
there is $t_{\varepsilon}<0$ such that
\begin{equation}\label{ex1.3+}
\begin{cases}
\e| \fr{u(t,x)-  v(t,x)}{\phi^\beta(-x\cdot e_i -c_f t)}\r| \leq \varepsilon
&\text { for all $t \leq t_{\varepsilon}$ and $x\in\bar{\mathcal H_i}\cap\bar\Omega$ with $i\in I$,}\\
\e| \fr{u(t,x)-  v(t,x)}{\phi^\beta(x\cdot e_j -c_f t)}\r| \leq \varepsilon
&\text { for all $t \leq t_{\varepsilon}$ and $x\in\bar{\mathcal H_j}\cap\bar\Omega$ with $j\in J$,}\\
\e| \fr{u(t,x)-  v(t,x)}{\phi^\beta(x\cdot e_j -c_f t)}\r| \leq \varepsilon
&\text { for all $t \leq t_{\varepsilon}$ and $x\in\overline{\Omega\backslash\bigcup_{s=1}^m{\mathcal H_s}}$ and for all $j\in J$.}\\
\end{cases}\end{equation}
 {Take any $t_{0} \leq \min \left(t_{\varepsilon}, T_{\delta}-\rho \varepsilon\right)$. For all $t\in\left[0, T_{\delta}-t_{0}-\rho\varepsilon\right]$, define the functions
\begin{align*}
u^\pm(t,x)=\bc
u(\alpha ^\pm(t),x)\pm\vp e^{-\dl t}\phi^\beta(\xi_i^\pm(t,x))\pm\vp e^{-\dl t}\phi^\beta(\chi_i^\pm(t,x))\\
~~~~~~~~~~~~~~~~~~~~~~~~~~~~~~~~~~~~~~~~~~~~~~~~~~~~~~~~~~
\text{ in }\bar{\{x\in{\mathcal H_i}:x\cdot e_i\geq L\}}\text{ with any }i\in I,\\
u(\alpha ^\pm(t),x)\pm2\vp e^{-\dl t}\phi^\beta(-L-c_f\alpha^\pm(t))
~~~~~~~~\text{ in }\bar{\Omega\backslash\bigcup_{i\in I}\{x\in\mathcal H_i:x\cdot e_i\geq L\}},
\ec
\end{align*}
where $\alpha^\pm(t)=t_{0}+t\pm\rho \varepsilon\e(1-e^{-\delta t}\r)$, $\xi_i^\pm(t,x)=-x\cdot e_i-c_f\alpha^\pm(t)$,
$\chi^\pm_i(t,x)=x\cdot e_i-c_f\alpha^\pm(t)-2L$,
  $\rho>0$ is a constant such that
\begin{align}\label{rrr}
\rho\hat\kappa \delta \geq 2 \dl+2\e\|\phi''/\phi\r\|_{L^\i(\mathbb R)}+ 2L_f
\end{align}
 with $\hat\kappa$  defined by \eqref{ex1.3} and $L_f$  defined by \eqref{deflf}.
We shall show that the function $u^{+}(t, x)$ is a supersolution of the problem satisfied by $  v\left(t_{0}+t, x\right)$
in $\left[0, T_{\delta}-t_{0}-\rho\varepsilon\right]\times{\Omega}$,
a subsolution $u^{-}(t, x)$   can be proved by some similar methods.
Fix any $i\in I$ and $j\in J$. Let us first verify the initial and boundary conditions.
By \eqref{ex1.3+}, one infers that
$$
v(t_0,x)\leq u(t_0,x)+\vp \phi^\beta(-x\cdot e_i-c_f t_0)\leq u^+(0,x)
$$
 for $x\in \bar{\mathcal H_i}$ such that $x\cdot e_i\geq L$. From \eqref{ex1.3+} and $\phi'<0$, it follows that
$$
v(t_0,x)\leq u(t_0,x)+\vp \phi^\beta(-x\cdot e_i-c_f t_0)\leq u(t_0,x)+\vp \phi^\beta(-L-c_f t_0)\leq u^+(0,x)
$$
 for $x\in \bar{\mathcal H_i}$ such that $x\cdot e_i\leq L$. Similarly, one can obtain that
$$
v(t_0,x)\leq u(t_0,x)+\vp \phi^\beta(x\cdot e_j-c_f t_0)\leq u(t_0,x)+\vp \phi^\beta(-L-c_f t_0)\leq u^+(0,x)
$$
for $x\in\bar{\mathcal H_j}$ such that $x\cdot e_j\geq L$ and $x\in\bar{\Omega\backslash\bigcup_{s=1}^m\{x\in\mathcal H_s:x\cdot e_s\geq L \}}$.
As a result, there holds
$v(t_0,x)\leq u^+(0,x)$ for $x\in\bar\Omega$.
In addition, since $\nu\cdot e_i=0$ for $x\in \p\mathcal H_i\cap\p\Omega$, then $u^{+}_\nu(t, x)= 0$ for $t\in\left[0, T_{\delta}-t_{0}-\rho \varepsilon\right]$ and $x\in\p\Omega$.}

 Since $u^+$ is $C^2$ for $x\in\mathcal H_i$ such that $x\cdot e_i\neq L$ and is $C^1$ for $x\in\bar\Omega$, in order to show that $u^+$ is a supersolution,
it suffices to check that $\mathscr L[u^+](t,x):=u^+_t(t,x)-\Delta u^+(t,x)-f(u^+(t,x))\g0$ for $t\in\left[0, T_{\delta}-t_{0}-\rho \varepsilon\right]$ and both for $x\in\bar{\mathcal H_i}\cap\overline{\Omega}$ with $x\cdot e_i\geq L$ and
$x\in\bar{\Omega}\backslash\bigcup_{i\in I}\bar{\{x\in\mathcal H_i:x\cdot e_i\geq L\}}$.
Note that $\alpha^+(t)\leq t_0+t+\rho\vp\leq T_\dl$.
We first consider the case $x\in\bar{\mathcal H_i}\cap\overline{\Omega}$ with $x\cdot e_i\geq L$.
After some computations, there holds
\begin{align*}
\mathscr L [u^+](t,x)
=&\rho\vp\dl e^{-\dl t}u_t(\alpha^+(t),x)-\vp \dl
 e^{-\dl t}  \phi^\beta (\xi^+_i(t,x))\\
&+\vp  e^{-\dl t}  \beta \phi^{\beta-1}(\xi^+_i(t,x))\phi'(\xi^+_i(t,x))
 (-c_f-c_f\rho\vp\dl e^{-\dl t})\\
&-\vp  e^{-\dl t}  \beta (\beta-1)\phi^{\beta-2}(\xi^+_i(t,x))(\phi'(\xi^+_i(t,x)))^2\\
&-\vp  e^{-\dl t}  \beta \phi^{\beta-1}(\xi^+_i(t,x))\phi''(\xi^+_i(t,x))
-\vp \dl
 e^{-\dl t}  \phi^\beta (\chi^+_i(t,x))\\
&+\vp  e^{-\dl t}  \beta \phi^{\beta-1}(\chi^+_i(t,x))\phi'(\chi^+_i(t,x))
 (-c_f-c_f\rho\vp\dl e^{-\dl t})\\
&-\vp  e^{-\dl t}  \beta (\beta-1)\phi^{\beta-2}(\chi^+_i(t,x))(\phi'(\chi^+_i(t,x)))^2\\
&-\vp  e^{-\dl t}  \beta \phi^{\beta-1}(\chi^+_i(t,x))\phi''(\chi^+_i(t,x))
+f (u(\alpha^+(t),x))-f (u^+(t,x)).
\end{align*}

 If $-C_\dl<\xi_i^+(t,x)<C_\delta$ (that is, $x\in \Omega^i_\dl(\alpha^+(t))$), then \eqref{ex1.3} implies that $u_t(\alpha^+(t),x)\geq\hat\kappa$. Since
$0\leq  u(\alpha^+(t),x)\leq u^+(t,x)\leq 1+2\vp\leq1+\vp_0$ by \eqref{ddddd} and $\phi<1$,
it then follows from   \eqref{deflf} and $\phi<1$  that
\begin{align*}
f (u(\alpha^+(t),x))-f (u^+(t,x))
\geq-L_f
\e(\vp e^{-\dl t}\phi^\beta(\xi_i^+(t,x))+\vp e^{-\dl t}\phi^\beta(\chi_i^+(t,x))\r)
\geq- 2L_f\vp e^{-\dl t}.
\end{align*}
 By $\phi'<0$, $0<\beta<1$, $0<\phi<1$ and \eqref{rrr}, one has
\begin{align*}
\mathscr L [u^+](t,x)\g &\rho\vp\dl e^{-\dl t} u_t(\alpha^+(t),x)-2\vp \dl
e^{-\dl t}
-2\vp  e^{-\dl t}\|\phi''/\phi\|_{L^\infty(\mathbb R)}
  -2 L_f\vp e^{-\dl t}\\
\g & \vp  e^{-\dl t}\e(\rho \hat \kappa\dl-2\dl-2\e\|\phi''/\phi\r\|_{L^\i(\R)}-2L_f
 \r)\\
\g& 0.
\end{align*}

 If $\xi_i^+(t,x)>C_\delta$, it follows from  \eqref{ddddd} and \eqref{cd1} that
 $0< u(\alpha^+(t),x)\leq u^{+}(t, x)\leq\dl+2\vp\leq2\dl\leq\theta$, hence  $f (u(\alpha^+(t),x))=0=f (u^+(t,x))$.
  Note that $\chi_i^+(t,x)\geq-c_f\alpha^+(t)-2L\geq-c_f T_\dl-2L\geq C$ by \eqref{tdelta}.
  Since  $u_t>0$, $\phi'<0$, one concludes from   \eqref{betac2} and
  \eqref{ddddd} that
that
\begin{align*}
\mathscr L [u^+](t,x)
\geq& \vp  e^{-\dl t}\phi^\beta (\xi^+_i(t,x))\e(\beta \e(\fr{\phi'(\xi^+_i(t,x))}{\phi(\xi^+_i(t,x))}\r)^2-\beta\fr{\phi''(\xi^+_i(t,x))}{\phi(\xi^+_i(t,x))}
 \r.\\
 &\e.-\beta^2  \e(\fr{\phi'(\xi^+_i(t,x))}{\phi(\xi^+_i(t,x))}\r)^2-c_f\beta \fr{\phi'(\xi^+_i(t,x))}{\phi(\xi^+_i(t,x))}-\dl \r)\\
 &+ \vp  e^{-\dl t}\phi^\beta (\chi^+_i(t,x))\e(\beta \e(\fr{\phi'(\chi^+_i(t,x))}{\phi(\chi^+_i(t,x))}\r)^2-\beta\fr{\phi''(\chi^+_i(t,x))}{\phi(\chi^+_i(t,x))}
 \r.\\
 &\e.-\beta^2  \e(\fr{\phi'(\chi^+_i(t,x))}{\phi(\chi^+_i(t,x))}\r)^2-c_f\beta \fr{\phi'(\chi^+_i(t,x))}{\phi(\chi^+_i(t,x))}-\dl \r)\\
    \g& \vp  e^{-\dl t}\phi^\beta (\xi^+_i(t,x))\e(-\fr{3\Lambda(-\beta c_f)}{4}+\fr{\Lambda(-\beta c_f)}{4}+\fr{\Lambda(-\beta c_f)}{4}\r)\\
&+ \vp  e^{-\dl t}\phi^\beta (\chi^+_i(t,x))\e(-\fr{3\Lambda(-\beta c_f)}{4}+\fr{\Lambda(-\beta c_f)}{4}+\fr{\Lambda(-\beta c_f)}{4}\r)\\
\g &0.
\end{align*}

If  $\xi^+_i(t,x)<-C_\delta$, then $1-\vp_0\leq1-\dl\l u(\alpha^+(t),x)\leq u^+(t,x)\l 1+2\vp\leq1+\vp_0$ by  \eqref{ddddd} and \eqref{cd2}.
From \eqref{vp0}, there holds
\begin{align*}
f (u(\alpha^+(t),x))-f (u^+(t,x))
\geq&-\fr34 f'(1)\left( \vp
 e^{-\dl t}\phi^\beta(\xi_i^+(t,x))+
 \vp
 e^{-\dl t}\phi^\beta(\chi_i^+(t,x))
 \right)\\
 \geq& -\fr34 f'(1) \vp
 e^{-\dl t}\phi^\beta(\xi_i^+(t,x)).
\end{align*}
Thus, one infers from $u_t>0$, $\beta<1$, $\phi'<0$, \eqref{betac2}, \eqref{defphi'} and
  \eqref{ddddd} that
\begin{align*}
\mathscr L [u^+](t,x)
\g&-\vp \dl e^{-\dl t}  \phi^\beta (\xi^+_i(t,x))
 -\fr34 f'(1)\vp e^{-\dl t}\phi^\beta(\xi^+_i(t,x))\\
 & -\vp  e^{-\dl t}  \beta \phi^{\beta-1}(\xi^+_i(t,x))\phi''(\xi^+_i(t,x))\\
&+ \vp  e^{-\dl t}\phi^\beta (\chi^+_i(t,x))\e(\beta \e(\fr{\phi'(\chi^+_i(t,x))}{\phi(\chi^+_i(t,x))}\r)^2-\beta\fr{\phi''(\chi^+_i(t,x))}{\phi(\chi^+_i(t,x))}
 \r.\\
 &\e.-\beta^2  \e(\fr{\phi'(\chi^+_i(t,x))}{\phi(\chi^+_i(t,x))}\r)^2-c_f\beta \fr{\phi'(\chi^+_i(t,x))}{\phi(\chi^+_i(t,x))}-\dl \r)\\
\g&\vp e^{-\dl t}\phi^\beta(\xi^+_i(t,x))\e(\fr{f'(1)}{4}+\fr{ f'(1)}{8}-\fr{3 f'(1)}{4}  \r)\\
&+ \vp  e^{-\dl t}\phi^\beta (\chi^+_i(t,x))\e(-\fr{3\Lambda(-\beta c_f)}{4}+\fr{\Lambda(-\beta c_f)}{4}+\fr{\Lambda(-\beta c_f)}{4}\r)\\
\g&0.
\end{align*}

Next, we consider the case $x\in\bar{\Omega}\backslash\bigcup_{i\in I}\bar{\{x\in\mathcal H_i:x\cdot e_i\geq L\}}$.
From \eqref{tdelta}, one obtains that $-L-c_f\alpha^+(t)\geq-L-c_fT_\delta\g C$.
  Since $-\vp_0\leq-2\vp\leq u(\alpha^+(t),x)\leq u^+(t,x)<\delta+2\vp\leq\theta$ by \eqref{ddddd} and \eqref{cd1},
  it then follows from \eqref{coma} and \eqref{ext} that $f(u(\alpha^+(t),x))=f(u^+(t,x))=0$.  Since $u_t>0$, $\beta<1$, $\phi'<0$, it follows from  \eqref{betac2} and \eqref{ddddd} that
\begin{align*}
\mathscr L [u^+](t,x)=
&f(u(\alpha^+(t),x))-f(u^+(t,x))+
\rho\vp\dl e^{-\dl t } u_t(\alpha^+(t),x)
-2\vp \dl e^{-\dl t}\phi^\beta(-L-c_f\alpha^+(t)) \\
 &+2\vp  e^{-\dl t}  \beta \phi^{\beta-1}(-L-c_f\alpha^+(t))\phi'(-L-c_f\alpha^+(t))
 (-c_f-c_f\rho\vp\dl e^{-\dl t})\\
 \geq&2\vp  e^{-\dl t}\phi^\beta(-L-c_f\alpha^+(t))\e(-\dl-c_f\beta\frac{\phi'(-L-c_f\alpha^+(t))}{\phi(-L-c_f\alpha^+(t))}\r)\\
 \geq&2\vp  e^{-\dl t}\phi^\beta(-L-c_f\alpha^+(t))\e(\frac{\Lambda(-\beta c_f)}{4}-\frac{3\Lambda(-\beta c_f)}{4}\r)\\
\geq&0.
\end{align*}

As a result, the function $u^{+}(t, x)$ is a supersolution of the problem satisfied by $  v\left(t_{0}+t, x\right)$ for
$t\in\left[0, T_{\delta}-t_{0}-\rho\varepsilon\right]$ and $x\in\Omega$.
By the comparison principle and $\phi<1$, one has that
$$
 u^{-}(t-t_0, x) \leq   v\left(t, x\right) \leq u^{+}(t-t_0, x)
$$
for $t \in\left[t_0, T_{\delta}-\rho \varepsilon\right]$ and
$x \in {\Omega}$.
Letting $t_{0} \to-\infty$ in above inequalities, there holds
\be\lb{ex1.4}
 u(t-\rho \varepsilon, x)\leq   v(t, x) \leq  u(t+\rho \varepsilon, x)\
\text{ for $t \in\left(-\infty, T_{\delta}-\rho \varepsilon\right]$ and $x \in \Omega$.}
\ee
Applying the comparison principle again, one gets that \eqref{ex1.4} holds for all $t \in \mathbb{R}$ and $x \in\Omega$. Since $\varepsilon$ is arbitrary given, one obtains that $  v(t, x) \equiv  u(t, x)$.
This completes the proof.
\end{pr}

\SE{Convergence of the entire solutions}\label{s3}
In this section, we consider the large time behavior of  the entire
solution  obtained in Theorem \ref{existence-}, that is, we
prove Theorem \ref{existence} in the sequel. To this end, we first give four auxiliary lemmas.
\begin{lem}\lb{lem2.1}
 {Let $u(t,x)$ be the entire solution obtained in Theorem \ref{existence-}.}
For any $\beta\in(0,1)$ and any small $\dl>0$, there exist $t_{1} \in \mathbb{R}$, $t_{2}
\in \mathbb{R}$, $\tau_{1} \in \mathbb{R}$, $\tau_2\in\R$ and  $\mu>0$  such that, for every $j \in J$,
\begin{align}\lb{eq2.1}
u(t, x)\leq&   \phi\left(x \cdot e_{j}-c_{f}t
+\tau_{1}\right)+\delta \phi^\beta\left(x \cdot e_{j}-c_{f}t+\tau_{1}\right)
  e^{-\delta\left(t
-t_{1}\right)}\nonumber \\
&  +\delta  \phi^\beta\left(x \cdot e_{j}-c_{f}t+\tau_{1}\right) e^{-\mu\left(x \cdot e_{j}-L\right)}
\end{align}
for $t \geq t_{1}$ and $x \in \overline{\mathcal{H}_{j}}$
such that $x \cdot e_{j} \geq L$, and
\begin{align}\lb{eq2.2}
u(t, x) \geq &   \phi\left(x \cdot e_{j}-c_{f}
t
-L\right)-\delta\phi^\beta\left(x \cdot e_{j}-c_{f}
t
+\tau_2\right) e^{-\delta\left(t-t_{2}\right)}\nonumber\\
&- \delta\phi^\beta\left(x \cdot e_{j}-c_{f}
t
+\tau_2\right)
e^{-\mu\left(x \cdot e_{j}-L\right)}
\end{align}
for $t \geq t_{2}$ and $x \in \overline{\mathcal{H}_{j}}$
such that $x \cdot e_{j} \geq L$.
\end{lem}
\begin{pr}
The proof is divided into three steps.

{\it Step 1: choice of some parameters.}
Fix any $\beta\in(0,1)$.
Choose $\delta>0$ such that
\be\lb{defdelta}
\delta<\min \left(c_f,\fr16,\fr{\varepsilon_0}{4},\fr\theta4,\frac{|\Lambda(-\beta c_f)|}{12},
\frac{\left|f'(1)\right|}{4},\frac{\left|f'(1)\right|}{4}c_f,\frac{|\Lambda(-\beta c_f)|}{4}c_f
\right).
\ee
Take $0<\mu<1$   such that
\be\lb{defmu}
\mu c_f=\dl.
\ee
Fix now an index $j \in J$. By \eqref{dm0},  there exists a real number $t_1<0$ such that
\be\lb{eq2.4}
u\left(t_{1}, x\right)\leq\delta\phi^\beta(x\cdot e_j-c_ft_1)   \text
{ for all } x \in \overline{\mathcal{H}_{j}} \text
{ such that } x \cdot e_{j} \geq L.
\ee
By the property of $\phi$, there exists $C>0$  such that
\begin{equation}\label{phic}
\begin{cases}
\phi\g   1- \fr\dl2, \ \
\e|\fr{\phi ' }{\phi}\r|\l \frac{\left|f'(1)\right|}{8}, \ \
\e|\fr{\phi '' }{\phi}\r|\l \frac{\left|f'(1)\right|}{4}\
&\text{in $(-\i,-C]$},\\
\phi^{1-\beta}\l \dl, \ \
-\fr32 c_f\leq\fr{\phi ' }{\phi }\leq-\frac{1}{2}c_f
&\text{in $[C,+\i)$},\\
\e|\e(\fr{\phi ' }{\phi }\r)^2-  \fr{\phi '' }{\phi }\r| \l -\frac1{4}{\Lambda(-\beta c_f)},\ \
c_f \beta\fr{\phi '}{\phi}+\beta^2\e(\fr{\phi '}{\phi}\r)^2\l \frac3{4}{\Lambda(-\beta c_f)}
&\text{in $[C,+\i)$}.
\end{cases}
\end{equation}
Define $\kappa>0$ such that
\be\lb{defkappa}
 \phi'
\l -\kappa\ \text{ in } [-C,C].
\ee
Let $\rho>0$  be
large enough such that
\begin{align}\label{defrho1}
\rho   \kappa   \g   \dl
+ \e\|{\phi '' }/{\phi}\r\|_{L^\i(\R)}+L_f + e^{-\dl t_1 }\e(2 \e\|{\phi ' }/{\phi}\r\|_{L^\i(\R)}+ \e\|{\phi '' }/{\phi }\r\|_{L^\i(\R)}+1+L_f\r).
\end{align}
Choose
\begin{align}\label{defr}
R \geq \rho+2 C+L-c_ft_1(>\rho+2 C).
\end{align}

{\it Step 2: proof of \eqref{eq2.1}.}
For $t \geq t_{1}$ and $x \in \overline{\mathcal{H}_{j}}$
with $x \cdot e_{j} \geq L$,  set
$$
\bar{  v}(t, x)=\min \left( \phi(\bar{\xi}(t, x))
+\delta e^{-\delta\left(t
-t_{1}\right)} \phi^\beta(\bar{\xi}(t, x))
+\delta   \phi^\beta(\bar{\xi}(t, x))
e^{-\mu\left(x \cdot e_{j}-L\right)},  1\right),
$$
where
$$
\bar{\xi}(t, x)=x \cdot e_{j}-c_f t
+\rho e^{-\delta\left(t-t_{1}\right)}-\rho-L-R+C.
$$
Let us verify that $\bar{  v}(t, x)$ is a supersolution
of the problem satisfied by $u(t, x)$ for $t \geq t_{1}$
and $x \in \overline{\mathcal{H}_{j}}$ with $x \cdot e_{j}
\geq L$.
We first verify the initial and boundary conditions.
By \eqref{defr}, one has
\[\bar{\xi}(t_1, x)= x \cdot e_{j}-c_ft_1-L-R+C\l x \cdot e_{j}-c_f  t_1.\]
It then follows from \eqref{eq2.4} and $\phi'<0$ that
\[\bar{  v}\left(t_{1}, x\right) \geq \delta \phi^\beta(\bar{\xi}(t_1, x)) \geq \dl \phi^\beta(x \cdot e_{j}-c_f  t_1) \g
 u\left(t_{1}, x\right)\]
for all $x \in \overline{
\mathcal{H}_{j}}$ with $x \cdot e_{j} \geq L$. On the
other hand, for $t \geq t_{1}$ and $x \in \overline{
\mathcal{H}_{j}}$ with $x \cdot e_{j}=L$, one has $
\bar{\xi}(t, x) \leq-c_ft_1-R +C \leq-C$ by \eqref{defr}.
It then follows from $\beta<1$, $\phi<1$ and \eqref{phic} that
$\phi^\beta\g \phi$ and $\phi \g 1-\dl/2$, which  {further} implies that
\[1\geq\bar{  v}(t, x) \geq
\min \left( 1
- \delta/2+ \delta\e(1
- \delta/2\r),   1\right)=  1.\]
Thus, $\bar{  v}(t, x)\geq u(t,x)$ for $t \geq t_{1}$ and $x \in \mathcal{H}_{j}$
with $x \cdot e_{j}=L$. Lastly, it is immediate to see that $\bar{  v}_{\nu}(t, x)
=  0$
for $t \geq t_{1}$ and $x \in \partial \mathcal{H}_{j}$
with $x \cdot e_{j}\geq L$ and $\bar{  v}(t, x)< 1$, since
$\nu\cdot e_{j}=0$ on $\partial\mathcal H_j\cap\p\Omega$.

Let us  check that
\[\mathscr{L}  [\bar{  v}](t, x):=
\bar{v }_{t}(t, x)-  \Delta \bar{v }(t, x)
-f (\bar{  v}(t, x))
\geq 0\]
for every $t \geq t_{1}$ and $x \in
\overline{\mathcal{H}_{j}}$ such that $x \cdot e_{j}
\geq L$ and $\bar{  v}(t, x)< 1$.
After a straightforward
computation, we get that
\begin{align*}
\mathscr{L}  [\bar{  v}](t, x)=& -\rho \dl \phi '(\bar \xi(t,x))
e^{-\dl (t-t_1)} -\dl^2 e^{-\dl(t-t_1)}\phi^{\beta}(\bar \xi(t,x))\\
& +\dl \beta \phi^{\beta-1}(\bar \xi(t,x))\phi '(\bar \xi(t,x))
\e(-c_f -\rho \dl e^{-\dl (t-t_1)}\r)e^{-\dl(t-t_1)}\\
& +\dl \beta \phi^{\beta-1}(\bar \xi(t,x))\phi '(\bar \xi(t,x))
\e(-c_f -\rho \dl e^{-\dl (t-t_1)}\r)e^{-\mu (x\cdot e_j-L)}\\
&-\dl \beta(\beta-1)\phi^{\beta-2}(\bar \xi(t,x))(\phi '(\bar \xi(t,x)))^2e^{-\dl(t-t_1)}\\
&- \dl\beta \phi^{\beta-1}(\bar \xi(t,x))\phi ''(\bar \xi(t,x)) e^{-\dl(t-t_1)}\\
&-\dl \beta(\beta-1)\phi^{\beta-2}(\bar \xi(t,x))(\phi '(\bar \xi(t,x)))^2e^{-\mu (x\cdot e_j-L)}\\
&+2\dl \mu \beta \phi^{\beta-1}(\bar \xi(t,x))\phi '(\bar \xi(t,x)) e^{-\mu (x\cdot e_j-L)}\\
&- \dl\beta \phi^{\beta-1}(\bar \xi(t,x))\phi ''(\bar \xi(t,x)) e^{-\mu (x\cdot e_j-L)}\\
&-\dl \mu^2 \phi^{\beta}(\bar \xi(t,x)) e^{-\mu (x\cdot e_j-L)}
+ f ( \phi(\bar{\xi}(t, x)))
-f (\bar{  v}(t, x)).
\end{align*}

If $\bar{\xi}(t, x)<-C$, one then gets from \eqref{phic} that $ 1>\bar{v}(t, x)\geq\phi(\bar{\xi}(t, x))
\geq   1- \delta/2\g1-\vp_0$. By \eqref{vp0}, there holds
\begin{align*}
f ( \phi(\bar{\xi}(t, x)))
-f (\bar{  v}(t, x))
\g& -\fr34 f'(1)\e(\delta\phi^\beta(\bar{\xi}(t, x))e^{-\delta\left(t-t_{1}\right)}
+\delta   \phi^\beta(\bar{\xi}(t, x))
e^{-\mu\left(x \cdot e_{j}-L\right)}\r).
\end{align*}
Since $\phi ^\prime<0$, $\beta<1$ and $\mu<1$,
it then can be inferred from \eqref{defdelta} and \eqref{phic} that
\begin{align*}
\mathscr{L}  [\bar{  v}](t, x)
&\g \delta  e^{-\delta\left(t-t_{1}\right)}\phi^\beta(\bar{\xi}(t, x))\e(-\fr34 f'(1)-\dl-\beta \fr{\phi ''(\bar \xi(t,x))}{\phi (\bar \xi(t,x))}\r)\\
&~~~+\delta   \phi^\beta(\bar{\xi}(t, x))
e^{-\mu\left(x \cdot e_{j}-L\right)}\e(-\fr34 f'(1)-\mu^2-\beta \fr{\phi ''(\bar \xi(t,x))}{\phi (\bar \xi(t,x))} +2\mu\beta \fr{\phi '(\bar \xi(t,x))}{\phi (\bar \xi(t,x))}\r)\\
&\g \delta  e^{-\delta\left(t-t_{1}\right)}\phi^\beta(\bar{\xi}(t, x))\e(-\fr34 f'(1)+\fr 14f'(1)+\fr 14f'(1)\r)\\
&~~~+\delta   \phi^\beta(\bar{\xi}(t, x))
e^{-\mu\left(x \cdot e_{j}-L\right)}\e(-\fr34 f'(1)+\fr 14f'(1)+\fr 14f'(1) +\fr 14f'(1)\r)\\
&\g0.
\end{align*}

If $\bar{\xi}(t, x)>C$, one then gets from \eqref{phic} that $\phi\leq\phi^{1-\beta}<\dl\l\theta$ and $\bar v(t,x)\l 3\dl\l \theta$, hence
$f ( \phi(\bar{\xi}(t, x)))
=0=f (\bar{  v}(t, x))$.
By \eqref{defdelta}, \eqref{defmu} and \eqref{phic}, there holds
\begin{align*}
\mathscr{L}  [\bar{  v}](t, x)
&\geq  \dl  e^{-\dl(t-t_1)}\phi^{\beta}(\bar \xi(t,x))\e(-\dl-\beta c_f \fr{\phi '(\bar \xi(t,x))}{\phi(\bar \xi(t,x))}
-\beta^2\e(\fr{\phi '(\bar \xi(t,x))}{\phi(\bar \xi(t,x))}\r)^2 \r.\\
&\e.+\beta \e(\fr{\phi '(\bar \xi(t,x))}{\phi(\bar \xi(t,x))}\r)^2-\beta \fr{\phi ''(\bar \xi(t,x))}{\phi(\bar \xi(t,x))}  \r)\\
& + \dl  e^{-\mu (x\cdot e_j-L)}\phi^{\beta}(\bar \xi(t,x))\e(-\mu^2-\beta c_f \fr{\phi '(\bar \xi(t,x))}{\phi(\bar \xi(t,x))}-\beta^2\e(\fr{\phi '(\bar \xi(t,x))}{\phi(\bar \xi(t,x))}\r)^2 \r.\\
&\e.+\beta \e(\fr{\phi '(\bar \xi(t,x))}{\phi(\bar \xi(t,x))}\r)^2-\beta \fr{\phi ''(\bar \xi(t,x))}{\phi(\bar \xi(t,x))} +2\mu\beta \fr{\phi '(\bar \xi(t,x))}{\phi(\bar \xi(t,x))} \r)\\
\g& \dl  e^{-\dl(t-t_1)}\phi^{\beta}(\bar \xi(t,x))\e(-\fr34 \Lambda(-\beta c_f)+\fr14 \Lambda(-\beta c_f) +\fr1{12} \Lambda(-\beta c_f)  \r)\\
&+\dl e^{-\mu (x\cdot e_j-L)}\phi^{\beta}(\bar \xi(t,x))\e(-\fr34 \Lambda(-\beta c_f)+\fr14 \Lambda(-\beta c_f) +\fr14 \Lambda(-\beta c_f) +\fr14 \Lambda(-\beta c_f) \r)\\
\g&0.
\end{align*}

If $-C\l \bar{\xi}(t, x)\l C$,
then \eqref{defr} implies that
$x \cdot e_{j}-L \geq c_f t-\rho e^{-\delta(t-t_0)}+\rho+R-2 C\geq c_ft$,
hence
\[e^{-\mu\left(x \cdot e_{j}-L\right)} \leq
e^{-\mu c_f  t}
= e^{-\dl\left(t-t_{1}\right)-\delta t_1 }.\]
Moreover, it follows from  \eqref{deflf} that
\begin{align*}
f ( \phi(\bar{\xi}(t, x)))
-f (\bar{  v}(t, x))
\geq
-L_f
\e(\delta e^{-\delta\left(t
-t_{1}\right)} \phi^\beta(\bar{\xi}(t, x))
+\delta   \phi^\beta(\bar{\xi}(t, x))
e^{-\mu\left(x \cdot e_{j}-L\right)}\r).
\end{align*}
Since $\phi'<0$, $\mu\in(0,1)$ and $\beta\in(0,1)$, one then obtains from \eqref{defkappa}  and \eqref{defrho1} that
\begin{align*}
\mathscr{L}  [\bar{  v}](t, x)
&\g -\rho \dl \phi '(\bar \xi(t,x))e^{-\dl (t-t_1)} -\dl^2 e^{-\dl(t-t_1)}\phi^{\beta}(\bar \xi(t,x))
- \dl\beta \phi^{\beta-1}(\bar \xi(t,x))\phi ''(\bar \xi(t,x)) e^{-\dl(t-t_1)}\\
&~~~+2\dl \mu \beta \phi^{\beta-1}(\bar \xi(t,x))\phi '(\bar \xi(t,x)) e^{-\mu (x\cdot e_j-L)}
 - \dl\beta \phi^{\beta-1}(\bar \xi(t,x))\phi ''(\bar \xi(t,x)) e^{-\mu (x\cdot e_j-L)}\\
&~~~-\dl \mu^2 \phi^{\beta}(\bar \xi(t,x)) e^{-\mu (x\cdot e_j-L)}
 + f ( \phi(\bar{\xi}(t, x)))-f (\bar{  v}(t, x))\\
&\g
\rho \dl \kappa e^{-\dl (t-t_1)} -\dl e^{-\dl (t-t_1)}\e(\dl
+ \e\|{\phi ''}/{\phi}\r\|_{L^\i(\R)}+L_f\r)\\
&~~~-\dl  e^{-\dl\left(t-t_{1}\right)-\delta t_1 }\e(2  \e\|{\phi '}/{\phi}\r\|_{L^\i(\R)}+ \e\|{\phi ''}/{\phi}\r\|_{L^\i(\R)}+1+L_f\r)
\\
&\g  0.
\end{align*}

According to the comparison principle, there holds
$$
\begin{aligned}
u(t, x) \leq \bar{  v}(t, x) \leq & \phi\left(x \cdot
e_{j}-c_ft +\rho e^{-\delta\left(t
-t_{1}\right)}-\rho-L-R+C\right)\\
&+\delta e^{-\delta\left(t
-t_{1}\right)} \phi^\beta\left(x \cdot
e_{j}-c_ft +\rho e^{-\delta\left(t
-t_{1}\right)}-\rho-L-R+C\right)\\
&+\delta   \phi^\beta \left(x \cdot
e_{j}-c_ft +\rho e^{-\delta\left(t
-t_{1}\right)}-\rho-L-R+C\right)
e^{-\mu\left(x \cdot e_{j}-L\right)}
\end{aligned}
$$
for all $t \geq t_{1}$ and $x \in \overline{\mathcal{H}_{j}}$
with $x \cdot e_{j} \geq L$.
Let
\[\tau_{1}=-\rho-L-R+C.\]
It then
can be inferred from $\phi '<0$ that \eqref{eq2.1} holds.

{\it Step 3: proof of \eqref{eq2.2}.} Since the propagation of $u$ is complete in the sense of \eqref{coms} and $\phi(-\i)=1$, one has $(1-u(t, \cdot))/\phi^\beta(x\cdot e_j-c_f t)
\rightarrow   0$ as $t\to +\i$ locally uniformly
in $x\in\overline{\mathcal{H}_{j}}\cap\bar\Omega$. Then there is $t_{2}>0$ such that, for every $j \in J$,
\begin{align}\label{ut2}
u(t_2, x) \geq   1-\delta \phi^\beta(x\cdot e_j-c_f t_2)\
\frac{}{}\text{ for  }t \geq t_{2}\text{ and }x \in \overline{\mathcal{H}_{j}}\text{ with }L \leq
x \cdot e_{j} \leq L+R+c_ft_2,
\end{align}
where the constant $R$ is defined as in \eqref{defr}. For
$t \geq t_{2}$ and $x \in \overline{\mathcal{H}_{j}}$
with $x \cdot e_{j} \geq L$,  define
$$
\underline{v}(t, x)=\max \left( \phi(\underline{\xi}
(t, x))-\delta   \phi^\beta (\underline{\xi}
(t, x)) e^{-\delta\left(t-t_{2}\right)}-
\delta    \phi^\beta (\underline{\xi}
(t, x)) e^{-\mu\left(x \cdot e_{j}-L\right)},  0\right),
$$
where
$$
\underline{\xi}(t, x)=x \cdot e_{j}-c_ft -\rho e^{-\delta\left(t-t_{2}\right)}
+\rho-L-R+C.
$$
Let us check that $\underline{  v}(t, x)$ is a subsolution
of the problem satisfied by $u(t, x)$ for $t
\geq t_{2}$ and $x \in \overline{\mathcal{H}_{j}}$
with $x \cdot e_{j} \geq L$.
By a similar argument to that of  step 2, one has
 \[\mathscr{L} [\underline{  v}](t,x)=
\underline{v}_{t}(t, x)-  \Delta
\underline{v }(t, x)-f (\underline{  v}(t, x)) \leq 0
\]
for all $t \geq t_{2}$ and $x \in \overline{\mathcal{H}_{j}}$
such that $x \cdot e_{j} \geq L$ and
$\underline{  v}(t, x)> 0$.
It suffices to verify the
initial and  boundary conditions. At time $t_{2}$,
it follows from $\phi'<0$, \eqref{defr} and \eqref{ut2} that
\[
\underline{  v}\left(t_{2},
x\right) \leq   1-\delta \phi^\beta(x\cdot e_j-c_f t_2)\leq u\left(t_{2}, x\right)
\]
for all $x \in \overline{\mathcal{H}_{j}}$ such that
$L\l x \cdot e_{j} \leq L+R+c_ft_2$. On the other hand,
for $x \in \overline{\mathcal{H}_{j}}$ such that
$x \cdot e_{j} \geq L+R+c_ft_2$, one has that
$\underline{\xi}\left(t_{2},
x\right) \geq L+R+c_ft_2-c_ft_2 -L-R+C=C$, it then follows from \eqref{phic} that $\phi^{1-\beta}(\bar \xi(t,x))\l \dl$, hence
\[
\underline{  v}
\left(t_{2}, x\right) \leq \max \left(\phi^\beta(\underline \xi(t,x)) (\phi^{1-\beta}(\underline \xi(t,x))-\dl)  -\delta \phi^\beta (\underline \xi(t,x))e^{-\mu\left(x
\cdot e_{j}-L\right)},   0\right)
=  0 \leq u\left(t_{2}, x\right).
\]
As a result, there holds
$\underline{  v}\left(t_{2}, x\right) \leq u\left(t_{2},
x\right)$ for all $x \in \overline{\mathcal{H}_{j}}$
with $x \cdot e_{j} \geq L$.
By $u_t>0$, $\phi'<0$, \eqref{defr} and \eqref{ut2}, one has
$\underline{v}(t, x)\l 1-\dl \phi^\beta(-R+\rho+C-c_ft)\l 1-\dl \phi^\beta(L-c_ft_2)\leq u(t_2,x)\leq u(t, x)$
for $t \geq t_{2}$ and $x \in \overline{\mathcal{H}_{j}}$ with
$x \cdot e_{j}=L$. Since $\nu\cdot e_j=0$ on $\partial\mathcal H_j\cap\p\Omega$,
then
$\underline{  v}_{\nu}(t, x)=  0$ for $t \geq t_{2}$ and
$x \in \partial \mathcal{H}_{j}$ with $x \cdot e_{j}\geq L$.

By \eqref{defr}, one has
$\underline{\xi}(t, x)\l x \cdot e_{j}-c_f
t +\rho-L-R+C\l x \cdot e_{j}-c_f
 t -L$
and
$
\underline{\xi}(t, x)\g x \cdot e_{j}-c_ft-L-R$.
Since $\phi'<0$, one gets that
\[
\phi(x \cdot e_{j}-c_f
t -L)\l
\phi (\underline \xi(t,x))\l \phi(x \cdot e_{j}-c_f
t -L-R),
\]
hence
\begin{align*}
\underline{  v}(t, x)\g &  \phi(x \cdot e_{j}-c_f
t -L)-\dl\phi^\beta(x \cdot e_{j}-c_f
 t-L-R) e^{-\dl(t-t_2)}\\
&-\dl  \phi^\beta(x \cdot e_{j}-c_f
t-L-R) e^{-\mu(x\cdot e_j-L)}
\end{align*}
 for $t\geq t_2$ and $x\in\bar{\mathcal H_j}$ such that $x\cdot e_j\geq L$. Let $\tau_2=-L-R(<0)$. According to the comparison principle, one gets that \eqref{eq2.2} holds. This completes the proof.
\end{pr}

\begin{lem}\lb{lem2.2}
For any $\beta\in(0,1)$ and $\varepsilon>0$, there exist
$t_{\varepsilon} \in \mathbb{R}$ and $\tau_\vp<0$ such that,
for every $j \in J$,
\begin{align*}
u(t, x) &\geq \phi\left(x \cdot e_{j}-c_f t
-L\right)-\varepsilon
\phi^\beta\left(x \cdot e_{j}-c_f t
+\tau_\vp\right) e^{-\delta\left(t-t_{\varepsilon}\right)}\\
&~~~-\varepsilon \phi^\beta\left(x \cdot e_{j}-c_f t
+\tau_\vp\right)
e^{-\mu\left(x \cdot e_{j}-L\right)}
\end{align*}
for all $t \geq t_{\varepsilon}$ and $x \in
\overline{\mathcal{H}_{j}}$ such that $x \cdot
e_{j} \geq L$, where the  constants $\delta>0$ and $\mu>0$ are defined by \eqref{defdelta}
and \eqref{defmu}, respectively.
\end{lem}
\begin{pr}
Fix any $\beta\in(0,1)$. Let $\dl>0$ be defined as in \eqref{defdelta}.
For any $\varepsilon>0$, define
\be\lb{defhatvp}
\hat{\varepsilon}={\varepsilon}/{\delta}.
\ee
It follows from  Lemma \ref{lem2.1} that the conclusion
of Lemma \ref{lem2.2} holds
for $\varepsilon \geq \delta$ with
$t_{\varepsilon}=t_{2}$ and $\tau_\vp=\tau_2$, where $t_2>0$ and $\tau_2<0$
are given in the proof of \eqref{eq2.2}. It remains
to consider the case $0<\varepsilon<\delta$. In this
case, there is $C_{\varepsilon} \geq C>0$ such that \eqref{phic} holds  with $C$ and $\delta$ replaced by $C_\vp$ and $\hat\vp\delta$, respectively. Let
 $R_{\varepsilon}$ be large enough so that
$R_{\varepsilon} \geq \hat{\varepsilon} \rho+
2 C_{\varepsilon}$ with $\rho$ defined in \eqref{defrho1}.
By some similar arguments as  the
proof of \eqref{eq2.2}, one can show that there exists $t_{\varepsilon}
\in \mathbb{R}$ such that, for every $j \in J$, the function
$$
\underline{ v}(t, x)=\max \left(
 \phi(\underline{\xi}(t, x))
-\hat{\varepsilon} \delta \phi^\beta (\underline{\xi}(t, x))
 e^{-\delta\left(t-
t_{\varepsilon}\right)}-\hat{\varepsilon} \delta
\phi^\beta(\underline{\xi}(t, x))
e^{-\mu\left(x \cdot e_{j}-L\right)},  0\right)
$$
is a subsolution of the problem satisfied by $u(t, x)$
for $t \geq t_{\varepsilon}$ and $x \in
\overline{\mathcal{H}_{j}}\cap\bar\Omega$ with $x \cdot e_{j} \geq L$,
where $\underline{\xi}(t, x)=x \cdot e_{j}-c_f  t -\hat{\varepsilon} \rho e^{-
\delta\left(t-t_{\varepsilon}\right)}+\hat{\varepsilon}
\rho-L-R_{\varepsilon}+C_{\varepsilon}$.
By $R_{\varepsilon} \geq \hat{\varepsilon} \rho+
2 C_{\varepsilon}$, one has
$x \cdot e_{j}-c_ft-L-R_\vp\leq\underline{\xi}(t, x)
 \l x \cdot e_{j}-c_f
 t -L$.
Together with $\phi'<0$, there holds
$
\phi(x \cdot e_{j}-c_f
t -L)\l
\phi (\underline \xi(t,x))\l \phi(x \cdot e_{j}-c_f
t -L-R_\vp)$. Let $\tau_\vp=-L-R_\vp(<0)$. By the comparison principle, one has
\begin{align*}
u(t,x)\geq\underline{  v}(t, x)\g &  \phi(x \cdot e_{j}-c_f
t -L)-\hat\vp\dl\phi^\beta(x \cdot e_{j}-c_f
 t+\tau_\vp) e^{-\dl(t-t_\vp)}\\
&-\hat\vp\dl  \phi^\beta(x \cdot e_{j}-c_f
t+\tau_\vp) e^{-\mu(x\cdot e_j-L)}
\end{align*}
 for $t\geq t_\vp$ and $x\in\bar{\mathcal H_j}$ such that $x\cdot e_j\geq L$.  The proof is complete.
\end{pr}

The next lemma is the stability of
the planar traveling front $\phi(x\cdot e_j-c_f t)$ in any branch $\mathcal{H}_{j}$.
\begin{lem}\lb{lem2.3}
 For any $\beta\in(0,1)$, $j \in J$ and $\varepsilon>0$,
if there are $t_{0} \in \mathbb{R}$ and
$\tau \in \mathbb{R}$ such that
\[\sup _{x \in \overline{\mathcal{H}_{j}}, x \cdot e_{j}
\geq L}\fr{\left|u\left(t_{0}, x\right)- \phi\left(x \cdot e_{j}
-c_f  t_{0}+\tau\right)\right|}{\phi^\beta(x \cdot e_{j}
-c_f  t_{0}+\tau)} \leq \varepsilon
\]
as well as  $(1-\phi\left(L-c_f  t_{0}+\tau\right))/\phi^\beta\left(L-c_f  t_{0}+\tau\right)
\l \varepsilon$ and ${1-u(t, x) }
\l  \phi^\beta(L-c_ft_0+\tau+\rho\hat\vp) $
for all $t \geq t_{0}$ and $x \in \overline{\mathcal{H}_{j}}$
with $x \cdot e_{j}=L$, then there is a positive constant
$\bar M$ such that
$$
\sup _{x \in \bar{\mathcal{H}_{j}}, x \cdot e_{j} \geq L}
\fr{\left|u(t, x)- \phi\left(x \cdot e_{j}-c_f  t+\tau\right)
\right|}{\phi^\beta\left(x \cdot e_{j}-c_f  t+\tau\right)} \leq \bar Me^{{c_f\beta\rho\vp/\delta}}\vp\ \
\text { for all } t \geq t_{0},
$$
where $\delta$, $\rho$ and
$\hat{\varepsilon}$ are
defined by  \eqref{defdelta},
\eqref{defrho1} and \eqref{defhatvp}, respectively.
\end{lem}
\begin{pr}
Let $\mu$ be
defined as in  \eqref{defmu}. For $t \geq
t_{0}$ and $x \in \overline{\mathcal{H}_{j}}\cap\bar\Omega$ with $x
\cdot e_{j} \geq L$, define
$$
\underline u(t,x)=\max \left( \phi\left(\underline\xi(t,x)\right)
-\hat{\varepsilon}
\delta \phi^\beta (\underline\xi(t,x))
e^{-\delta\left(t-t_{0}\right)}-\hat{\varepsilon}
\delta \phi^\beta(\underline\xi(t,x))
e^{-\mu\left(x \cdot e_{j}-L\right)}, 0\right)
$$
and
$$
\bar u(t,x)=\min \left(\phi\left( \bar\xi(t,x)\right)+\hat{\varepsilon}
\delta \phi^\beta(\bar\xi(t,x)) e^{-\delta\left(t-t_{0}\right)}
+\hat{\varepsilon}
\delta \phi^\beta(\bar\xi(t,x))e^{-\mu\left
(x \cdot e_{j}-L\right)}, 1\right),
$$
where
$\underline\xi(t,x)=x \cdot e_{j}-c_f  t-
\rho \hat{\varepsilon} e^{-\delta\left(t-t_{0}\right)}
+\rho \hat{\varepsilon}+\tau$
 and
$\bar\xi(t,x)=x \cdot e_{j}-c_f  t+\rho
\hat{\varepsilon} e^{-\delta\left(t-t_{0}\right)}-
\rho \hat{\varepsilon}+\tau$. Since $\sup _{x \in \overline{\mathcal{H}_{j}}, x \cdot e_{j}
\geq L}\left({\left|u\left(t_{0}, x\right)- \phi\left(x \cdot e_{j}
-c_f  t_{0}+\tau\right)\right|}/{\phi^\beta(x \cdot e_{j}
-c_f  t_{0}+\tau)}\right) \leq \varepsilon $, there holds
$$
\underline u(t_0,x)\leq u(t_0,x)\leq \bar u(t_0,x)\ \ \text{ for all  $x \in \overline{\mathcal{H}_{j}}\cap\bar\Omega$ with $x
\cdot e_{j} \geq L$}.
$$
 On the one hand, for $t \geq t_{0}$ and $x\in\bar{\mathcal H_j}$ with $x \cdot e_{j}=L$,
 it follows from  $\phi'<0$ and  $u(t, x)
\geq  1-\vp \phi^\beta(L-c_f t_0+\rho\hat\vp+\tau)$ that
\begin{align*}
\underline u(t,x)\leq\phi(\underline\xi(t,x))
-\hat{\varepsilon}
\delta \phi^\beta(\underline\xi(t,x))
e^{-\mu\left(x \cdot e_{j}-L\right)}\l   1-\vp \phi^\beta(L-c_f t_0+\rho\hat\vp+\tau)\l u(t,x).
\end{align*}
 On the other hand, it can be deduced from $\phi'<0$ and $\phi\left(L-c_f  t_{0}+\tau\right)\geq  1- \varepsilon \phi^\beta \left(L-c_f  t_{0}+\tau\right)$ that
\begin{align*}
\bar u(t,x)\geq \phi( \bar\xi(t,x))
+\hat{\varepsilon}
\delta \phi^\beta(\bar\xi(t,x))e^{-\mu
(x \cdot e_{j}-L)}
\geq\phi\left(L-c_f  t_{0}+\tau\right)+\vp \phi^\beta \left(L-c_f  t_{0}+\tau\right)
\geq   1\geq u(t,x)
\end{align*}
for $t \geq t_{0}$ and $x\in\bar{\mathcal H_i}$ with $x \cdot e_{j}=L$.
Together with some similar arguments as  the proof of Lemmas \ref{lem2.1}
and \ref{lem2.2}, one infers that  $\underline u(t,x)$ and $\bar u(t,x)$
are sub- and  supersolutions
of the problem satisfied by $u(t, x)$ for $t \geq
t_{0}$ and $x \in \overline{\mathcal{H}_{j}}\cap\bar\Omega$ with $x
\cdot e_{j} \geq L$, respectively.
By the comparison principle, one has
\begin{align}\label{buusu}
\underline u(t,x)
\leq u(t, x) \leq \bar u(t,x)\ \ \text{for $t \geq t_{0}$ and $x \in \overline{\mathcal{H}_{j}}\cap\bar\Omega$
such that $x \cdot e_{j} \geq L$}.
\end{align}

For $t \geq t_{0}$ and $x \in \overline{\mathcal{H}_{j}}\cap\bar\Omega$
such that $x \cdot e_{j} \geq L$,
it can be derived from  $\phi'<0$ and the mean value theorem that there is $\tilde \theta_1\in(0,1)$ such that
\begin{align}\label{uleq}
\bar u(t, x)&\leq \phi\left(x \cdot e_{j}-c_f  t-
\rho \hat{\varepsilon}+\tau\right)+2 \hat{\varepsilon}
\delta \phi^\beta\left(x \cdot e_{j}-c_f  t-
\rho \hat{\varepsilon}+\tau\right) \nonumber\\
&= \phi\left(x \cdot e_{j}
-c_f  t+\tau\right)-
\rho \hat{\varepsilon} \phi' (x \cdot e_{j}-c_f  t-\tilde\theta_1
\rho \hat{\varepsilon}+\tau) +2 \hat{\varepsilon} \delta \phi^\beta\left(x \cdot e_{j}-c_f  t-
\rho \hat{\varepsilon}+\tau\right).
\end{align}
For these $t$ and $x$, it is easy to check that there exists $\tilde C_0>1$ such that
\begin{align}\label{c00}
I(t,x):=\fr{ \phi^\beta\left(x \cdot e_{j}-c_f  t-
\rho \hat{\varepsilon}+\tau\right)}{\phi^\beta\left(x \cdot e_{j}-c_f  t +\tau\right)}\leq \tilde C_0e^{c_f\beta\rho\hat\vp}.
\end{align}
In fact, if $x\cdot e_j-c_ft-
\rho \hat{\varepsilon}+\tau>0$, then  $I(t,x)\leq K_4^\beta e^{c_f\beta\rho\hat\vp}/K_3^\beta$ by \eqref{estimate2}.
If $-
\rho \hat{\varepsilon}\leq x\cdot e_j-c_ft-
\rho \hat{\varepsilon}+\tau\leq0$, then \eqref{estimate2} and $\phi<1$ imply that $I(t,x)\leq e^{c_f\beta(x\cdot e_j-c_ft+\tau)}/K_3^\beta\leq e^{c_f\beta\rho\hat\vp}/K_3^\beta$. If $ x\cdot e_j-c_ft-\rho \hat{\varepsilon}+\tau\leq-\rho \hat{\varepsilon}$, it then follows from
$\phi<1$ and $\phi'<0$ that $I(t,x)\leq {1}/{\phi^\beta(0)}$.
Since $\tilde\theta_1\in(0,1)$, $\phi<\phi^\beta$ and $\phi'<0$, it then follows from \eqref{buusu}, \eqref{uleq} and \eqref{c00}  that
\begin{align*}
\fr{u(t, x)-\phi\left(x \cdot e_{j}-c_f  t+\tau\right)}
{\phi^\beta\left(x \cdot e_{j}-c_f  t+\tau\right)}
&\l -\rho \hat{\varepsilon}\fr{ \phi'\left(x \cdot e_{j}-c_f  t-\tilde\theta_1
\rho \hat{\varepsilon}+\tau\right)}{ \phi^\beta\left(x \cdot e_{j}-c_f  t+\tau\right)}+2 \hat{\varepsilon} \delta\fr{ \phi^\beta\left(x \cdot e_{j}-c_f  t-
\rho \hat{\varepsilon}+\tau\right)}{\phi^\beta\left(x \cdot e_{j}-c_f  t +\tau\right)}\\
&\l \rho \hat{\varepsilon}\left\|{\phi'}/{\phi}
\right\|_{L^{\infty}(\mathbb{R})}
\fr{ \phi\left(x \cdot e_{j}-c_f  t-
\tilde\theta_1\rho \hat{\varepsilon}+\tau\right)}{\phi^\beta\left(x \cdot e_{j}-c_f  t +\tau\right)}
+2 \hat{\varepsilon}\delta I(t,x)\\
&\l  \left(\rho \hat{\varepsilon}\left\|{\phi'}/{\phi}
\right\|_{L^{\infty}(\mathbb{R})}+2 \hat{\varepsilon}\delta\right) I(t,x)\\
&\l  \tilde C_0e^{{c_f\beta\rho\hat\vp}}\left(\rho \hat{\varepsilon}\left\|{\phi'}/{\phi}
\right\|_{L^{\infty}(\mathbb{R})}+2 \hat{\varepsilon}\delta\right)
\end{align*}
for $t \geq t_{0}$ and $x \in \overline{\mathcal{H}_{j}}\cap\bar\Omega$
such that $x \cdot e_{j} \geq L$.

For $t \geq t_{0}$ and $x \in \overline{\mathcal{H}_{j}}\cap\bar\Omega$
such that $x \cdot e_{j} \geq L$,
it can be derived from  $\phi'<0$ and the mean value theorem that there is $\tilde \theta_2\in(0,1)$ such that
\begin{align}\label{ugeq}
\underline u(t,x)& \geq \phi\left(x \cdot e_{j}-c_f  t+
\rho \hat{\varepsilon}+\tau\right)-2 \hat{\varepsilon}
\delta \phi^\beta\left(x \cdot e_{j}-c_f  t+\tau\right) \nonumber\\
&= \phi\left(x \cdot e_{j}
-c_f  t+\tau\right)+
\rho \hat{\varepsilon} \phi' \left(x \cdot e_{j}-c_f  t+\tilde\theta_2
\rho \hat{\varepsilon}+\tau\right) -2 \hat{\varepsilon} \delta \phi^\beta\left(x \cdot e_{j}-c_f  t
+\tau\right).
\end{align}
For these $t$ and $x$, since $\tilde\theta_2\in(0,1)$, $\phi<\phi^\beta$ and $\phi'<0$, one can deduced from \eqref{buusu} and \eqref{ugeq} that
\begin{align*}
\fr{u(t, x)-\phi\left(x \cdot e_{j}-c_f  t+\tau\right)}
{\phi^\beta\left(x \cdot e_{j}-c_f  t+\tau\right)}
&\geq \rho \hat{\varepsilon}\fr{ \phi'\left(x \cdot e_{j}-c_f  t+\tilde\theta_2
\rho \hat{\varepsilon}+\tau\right)}{ \phi^\beta\left(x \cdot e_{j}-c_f  t+\tau\right)}-2 \hat{\varepsilon} \delta\\
&\geq -\rho \hat{\varepsilon}\left\|{\phi'}/{\phi}
\right\|_{L^{\infty}(\mathbb{R})}{ \phi^{1-\beta}\left(x \cdot e_{j}-c_f  t+\tau\right)}
-2 \hat{\varepsilon}\delta\\
&\geq -\rho \hat{\varepsilon}\left\|{\phi'}/{\phi}
\right\|_{L^{\infty}(\mathbb{R})}
-2 \hat{\varepsilon}\delta
.
\end{align*}
Note that $\hat\vp=\vp/\dl$. As a consequence, there holds
\begin{align*}
\sup _{x \in
\overline{\mathcal{H}_{j}}, x \cdot e_{j} \geq L}
 \fr{\left|u(t, x)- \phi\left(x \cdot e_{j}-c_f t+\tau\right)
\right|}{\phi^\beta\left(x \cdot e_{j}
-c_f  t+\tau\right)}
\l  \tilde C_0e^{{c_f\beta\rho\hat\vp}}\left(\rho \hat{\varepsilon}\left\|{\phi'}/{\phi}
\right\|_{L^{\infty}(\mathbb{R})}+2 \hat{\varepsilon}\delta\right)
= \bar M e^{{c_f\beta\rho\vp/\delta}}\vp,
\end{align*}
where $
\bar M=\tilde C_0\e(\fr{\rho}{\delta}\left\|{ \phi'}/{\phi}
\right\|_{L^{\infty}(\mathbb{R})}+2\r)$ does not depend
on $j$, $\varepsilon$, $t_{0}$ and $\tau$. This completes the proof.
\end{pr}

\vspace{0.3cm}
\begin{pr}[Proof of Theorem \ref{existence}] We divide the proof into three steps.

{\it Step 1: choice of some parameters and notations.} Let $u$ be the entire solution obtained in Theorem \ref{existence-}. Take any $\beta\in(0,1)$. By Lemma \ref{lem2.1}, For any small $\dl>0$, there exist $t_{1} \in \mathbb{R}$, $t_{2} \in \mathbb{R}$,
$\tau_{1} \in \mathbb{R}$, $\tau_2\in\R$ and $\mu>0$
such that
\begin{align}\lb{eqth1-1}
& \phi\left(x \cdot e_{j}-c_f  t
-L\right)-\delta \phi^\beta\left(x \cdot e_{j}-c_f  t
+\tau_2\right) e^{-\delta\left(t-t_{2}\right)}\nonumber\\
&- \delta\phi^\beta\left(x \cdot e_{j}-c_f  t
+\tau_2\right)  e^{-\mu\left(x \cdot e_{j}-L\right)}\nonumber\\
\leq &  u(t, x) \leq  \phi\left(x \cdot
e_{j}-c_f  t +\tau_{1}\right)
+\delta \phi^\beta\left(x \cdot e_{j}-c_f  t +\tau_{1}\right)  e^{-\delta\left(t-t_{1}\right)}\nonumber\\
&
+\delta \phi^\beta\left(x \cdot e_{j}-c_f  t +\tau_{1}\right)
e^{-\mu\left(x \cdot e_{j}-L\right)}
\end{align}
for every $j \in J$, $t \geq\max \left(t_{1}, t_{2}\right)$ and $x \in
\overline{\mathcal{H}_{j}}$ with $x \cdot e_{j}\geq L$.
Consider now any sequence $\left(t_{k}\right)_
{k \in \mathbb{N}}$ such that $t_{k} \rightarrow
+\infty$ as $k \rightarrow+\infty$, and consider
any $j \in J$. Remember that $\mathcal{H}_{j}=
\mathcal{H}_{e_{j}, \omega, x_{j}}$ is a
straight half-cylinder as in \eqref{defhi1}. For every
$k \in \mathbb{N}$, let $\mathcal{H}_{j}^{k}=
\mathcal{H}_{j}-c_ft_{k} e_{j}$ be the shifted
half-cylinder in the direction $-e_{j}$. The
half-cylinders $\mathcal{H}_{j}^{k}$ converge to
a straight open cylinder $\mathcal{H}_{j}^{\infty}$
parallel to $e_{j}$ as $k \rightarrow+\infty$. Define
$$
 u_{k}(t, y)= u\left(t+t_{k}, y+c_ft_{k} e_{j}\right)\ \ \text{for $(t,y)\in\mathbb{R} \times \overline{\mathcal H_j^k}$}.
$$
From the standard parabolic estimates, up to extraction
of a subsequence, the functions
$u_{k}(t, y)$
 converge locally uniformly for $(t, y) \in \mathbb{R} \times \overline{
\mathcal{H}_{j}^{\infty}}$ as $k\to+\i$
to a solution
$ u_{\infty}(t, y)$ of
$$
\begin{cases}\left( u_{\infty}\right)_{t}-
  \Delta  u_{\infty}=f\left( u_{\infty}\right),
& t \in \mathbb{R}, \ y \in \mathcal{H}_{j}^{\infty},
\\ \left( u_{\infty}\right)_{\nu}=  0, & t \in \mathbb{R},\
y \in \partial \mathcal{H}_{j}^{\infty}.\end{cases}
$$
By \eqref{eqth1-1}, there holds
$$
 \phi\left(y \cdot e_{j}-c_f t
-L\right) \leq  u_{\infty}(t, y) \leq  \phi\left(y
\cdot e_{j}-c_f t +\tau_{1}\right)\ \
\text{ for all $(t, y) \in \mathbb{R} \times
\overline{\mathcal{H}_{j}^{\infty}}$.}
$$
Thus, $ u_{\infty}$ is a transition
front connecting $ 0$ and $ 1$ for \eqref{te1.1} in the
straight cylinder $\mathcal{H}_{j}^{\infty}$ with the sets
$$
\Gamma_{\i,t}=\e\{y\in\mathcal H_j^\i:y\cdot e_j=c_f t\r\}\ \ \text{ and }\ \
\Omega^\pm_{\i,t}=\e\{y\in\mathcal H_j^\i:\pm(y\cdot e_j-c_f t)<0\r\}.
$$
By some similar arguments as the proof of \cite[Lemma 3.5]{ghs}, one can show that $u_\i$ is a planar traveling front, i.e.,
there exists $\tau_j^* \in \mathbb{R}$
such that $ u_{\infty}(t, y)= \phi(y \cdot e_{j}
-c_ft+\tau_j^*)$ for all $(t, y) \in \mathbb{R}
\times \overline{\mathcal{H}_{j}^{\infty}}$. Therefore,
\be\lb{eqth1-2}
 u_{k}(t, y) \rightarrow  \phi\left(y \cdot e_{j}-c_f
t+\tau_j^*\right)\ \  \text{locally uniformly in $\mathbb{R}
\times \overline{\mathcal{H}_{j}^{\infty}}$ as
$k \rightarrow+\infty$.}
\ee

By a similar argument as the proof of \eqref{c00}, one gets that there is a constant $\tilde C_1\geq1$ such that
\begin{align}\label{ppp}
\frac{\phi^\beta(y\cdot e_j+\tau_1)}{\phi^\beta(y\cdot e_j+\tau_j^*)}\leq \tilde C_1\ \ \text{for all }y\in\bar{\mathcal H_j^k}.
\end{align}
Pick now any $0<\varepsilon<1$, and let $K_\vp>0$ such
that
 \begin{align}\label{Kvp}
 \phi \geq  1-\fr{\varepsilon}3\geq\left(1-\frac{\vp}{2}\right) \ \
\text{in $(-\infty, -K_\vp]$}\ \ \ \ \text{and} \ \ \ \
 \phi^{1-\beta} \leq \fr{\varepsilon}{3\tilde C_1}\e(<\frac{\vp}{2}\r) \ \
\text{in $[K_\vp,+\infty)$.}
\end{align}
Define
\begin{align*}
K_{1,\vp}=\max \left(K_\vp-\tau_{1}, K_\vp
-\tau_{j}^*\right)\ \ \text{and} \ \ K_{2,\vp}=\min \left(-K_\vp+L,-K_\vp-\tau_{j}^*\right)<K_{1,\vp}.
\end{align*}

{\it Step 2: stability of $\phi(x\cdot e_{j}-c_ft+\tau_{j}^*)$ for all $t\geq t_k$ and $x\in\bar{\mathcal H_j}$ such that $x\cdot e_j\geq L$.}
We first verify that for $k$ large enough,
\begin{align}\label{u-p}
\sup_{x \in
\overline{\mathcal{H}_{j}},
x \cdot e_{j} \geq L}\fr{| u(t_{k}, x)
- \phi(x \cdot e_{j}
-c_ft_{k}+\tau_{j}^*)|}{\phi^\beta (x \cdot e_{j}-c_ft_{k}
+\tau_j^*)} \leq   \varepsilon.
\end{align}
We claim that
\begin{equation}\label{clm}
\begin{cases}  0\l  \fr{u_{k}(0, y)}{\phi^\beta(y\cdot e_j+\tau_j^*)} \leq \frac{\varepsilon}{2} &
\text { for all } y \in \overline{\mathcal{H}_{j}^{k}}
\text { such that } y \cdot e_{j} \geq K_{1,\vp}, \\  1
-\varepsilon \leq  \fr{u_{k}(0, y)}{\phi^\beta(y\cdot e_j +\tau_j^*)}\l 1+\frac{\vp}{2}
& \text { for all }
y \in \overline{\mathcal{H}_{j}^{k}} \text { such
that } K_{2,\vp}-\frac{c_f }{2} t_{k} \leq y \cdot
e_{j} \leq K_{2,\vp}\end{cases}
\end{equation}
 for $k$ large enough. In fact, since $K_{1,\vp}+\tau_1\geq K_\vp$, it follows from \eqref{eqth1-1} and \eqref{Kvp} that
\begin{align*}
\fr{u_{k}(0, y)}{\phi^\beta(y\cdot e_j+\tau_1)}\l\phi^{1-\beta}(y\cdot e_j+\tau_1)+\delta e^{-\dl (t_k-t_1)}+\dl e^{-\mu(y\cdot e_j+c_f t_k-L)}\l \fr{\vp}{2\tilde C_1}
\end{align*}
for   $y \in \overline{\mathcal{H}_{j}^{k}}$
such that $y \cdot e_{j} \geq K_{1,\vp}$ and for $k$ large enough. Together with \eqref{ppp},  one has
\begin{align*}
0\leq\fr{u_{k}(0, y)}{\phi^\beta(y\cdot e_j+\tau_j^*)}=\fr{u_{k}(0, y)}{\phi^\beta(y\cdot e_j+\tau_1)}
\fr{\phi^\beta(y\cdot e_j+\tau_1)}{\phi^\beta(y\cdot e_j+\tau_j^*)}\l \frac{\vp}{2}
\end{align*}
for $y \in \overline{\mathcal{H}_{j}^{k}}$
such that $y \cdot e_{j} \geq K_{1,\vp}$ and for $k$ large enough. Hence, the first inequality of \eqref{clm} is true. We now turn to prove the second inequality of \eqref{clm}. On the one hand, since $\phi<1$  and $K_{2,\vp}-L\leq-K_\vp$, one obtains from $\phi^\beta<1$, \eqref{eqth1-1} and \eqref{Kvp} that
\[
\fr{u_{k}(0, y)}{\phi^\beta(y\cdot e_j +\tau_j^*)}
\geq \phi(y\cdot e_j-L)-\dl\phi^\beta(y\cdot e_j+\tau_2)e^{-\dl(t_k-t_2)}-\dl\phi^\beta(y\cdot e_j+\tau_2)e^{-\mu(y\cdot e_j+c_ft_k-L)}
\geq 1
-\varepsilon
\]
for all
$y \in \overline{\mathcal{H}_{j}^{k}}$ such
that $K_{2,\vp}-\frac{c_f }{2} t_{k} \leq y \cdot
e_{j} \leq K_{2,\vp}$ and for $k$ large enough. On the other hand, since $K_{2,\vp}+\tau_j^*\leq-K_\vp$, $\vp\in(0,1)$ and $\phi\leq\phi^\beta$, one then infers from \eqref{Kvp} that
\[
\fr{u_{k}(0, y)}{\phi^\beta(y\cdot e_j +\tau_j^*)}
\l\frac{1}{\phi(y\cdot e_j +\tau_j^*)}\l \fr{1}{1-\vp/3}\l 1+\frac{\vp}{2}
\]
for all
$y \in \overline{\mathcal{H}_{j}^{k}}$ such
that $K_{2,\vp}-{c_f t_{k}}/{2} \leq y \cdot
e_{j} \leq K_{2,\vp}$. The proof of the  claim \eqref{clm} is  complete.

Since $K_{1,\vp} +\tau_{j}^*\geq K_\vp$ and $K_{2,\vp}+\tau_{j}^* \leq-K_\vp$, one gets from \eqref{Kvp} that
\begin{align*}
 0\l   \phi^{1-\beta}(y \cdot e_{j}+
\tau_{j}^*) \leq  \varepsilon/2  \ \ \text{
for all $y \in \overline{\mathcal{H}_{j}^{k}}$
with $y \cdot e_{j} \geq K_{1,\vp}$}
\end{align*}
and
\begin{align*}
 1-\varepsilon/2  \leq
\phi(y \cdot e_{j}
+\tau_{j}^*) \leq
\phi^{1-\beta}(y \cdot e_{j}
+\tau_{j}^*)\l 1 \ \ \text{ for all $y \in \overline{
\mathcal{H}_{j}^{k}}$ with $y \cdot e_{j} \leq K_{2,\vp}$.}
\end{align*}
Combining with \eqref{clm}, one gets that, for $k$ large enough,
\begin{align*}
\fr{| u_{k}(0, y)- \phi(y \cdot e_{j}+\tau_{j}^*
)|}{\phi^\beta(y \cdot e_{j}
+\tau_{j}^*)}  \leq \varepsilon
\ \text { for all } y \in \overline{\mathcal{H}_{j}^{k}}
\text { with } y \cdot e_{j} \geq K_{1,\vp}  \text { or }
K_{2,\vp}-\frac{c_f }{2} t_{k} \leq y \cdot e_{j} \leq K_{2,\vp}.
\end{align*}
In addition, it follows from \eqref{eqth1-2} that
\begin{align*}
\sup _{ y \in \overline{\mathcal{H}_{j}^{k}},K_{2,\vp} \leq y \cdot e_{j} \leq K_{1,\vp}}\fr{
| u_{k}(0, y)- \phi(y \cdot e_{j}+\tau_{j}^*)
| }{\phi^\beta(y \cdot e_{j}+\tau_{j}^*)}\leq \varepsilon \ \  \text { for } k \text
{ large enough. }
\end{align*}
Therefore, one can deduce from the definitions of
$ u_{k}(t, y)$ and $\mathcal{H}_{j}^{k}$ that, for $k$ large enough,
\begin{align*}
&\fr{| u(t_{k}, x)- \phi(x \cdot e_{j}
-c_ft_{k}+\tau_{j}^*)|}{\phi^\beta(x \cdot e_{j}
-c_ft_{k}+\tau_{j}^*)} \leq \varepsilon \ \
\text { for all } x \in \overline{\mathcal{H}_{j}}
\text { such that } x \cdot e_{j} \geq K_{2,\vp}+
\frac{c_f }{2} t_{k}.
\end{align*}

Now, let us  show that, for $k$ large enough,
\begin{align}\lb{eqth1-4}
\fr{| u(t_{k}, x)
- \phi(x \cdot e_{j}
-c_ft_{k}+\tau_{j}^*)|}{\phi^\beta (x \cdot e_{j}-c_ft_{k}
+\tau_j^*)} \leq  \varepsilon
\ \ \text{
for all $x \in \overline{\mathcal{H}_{j}}$ such that $L \leq
x \cdot e_{j} \leq K_{2,\vp}+\frac{c_f }{2} t_{k}$}.
\end{align}
By Lemma \ref{lem2.2}, for $x \in
\overline{\mathcal{H}_{j}}$ with $L \leq x \cdot e_{j}
\leq K_{2,\vp}+c_ft_{k} / 2$ and for large $k$, there exist $t_\vp\in\R$ and $\tau_\vp<0$ such that
\begin{align*} u\left(t_{k}, x\right) \geq&  \phi
\left(x \cdot e_{j}-c_ft_{k}
-L\right)-\fr{\varepsilon}{4}  \phi^\beta
\left(x \cdot e_{j}-c_ft_{k}
+\tau_\vp\right) e^{-\delta\left(t_{k}-t_{\varepsilon}
\right)}\\& -\fr{\varepsilon}{4}  \phi^\beta
\left(x \cdot e_{j}-c_ft_{k}
+\tau_\vp\right) e^{
-\mu\left(x \cdot e_{j}-L\right)}.
\end{align*}
Since $t_k\to+\i$ as $k\to+\i$, $\phi(-\i)=1$ and $\phi'(-\i)=0$, one obtains from  the mean value theorem that there exists $\theta\in(0,1)$ such that
\begin{align*}
&~~\fr{u(t_{k}, x)
- \phi (x \cdot e_{j}
-c_ft_{k}+\tau_{j}^* )}{\phi^\beta  (x \cdot e_{j}-c_ft_{k}
+\tau_j^* )}\\
&\g \fr{\phi
 (x \cdot e_{j}-c_ft_{k}
-L )
- \phi (x \cdot e_{j}
-c_ft_{k}+\tau_{j}^* )}{\phi^\beta  (x \cdot e_{j}-c_ft_{k}
+\tau_j^* )}
-\fr{\varepsilon}{4} \fr{\phi^\beta
 (x \cdot e_{j}-c_ft_{k}+\tau_\vp)}{\phi^\beta  (x \cdot e_{j}-c_ft_{k}
+\tau_j^* )}e^{-\delta (t_{k}-t_{\varepsilon})}\\
&~~~~-\fr{\varepsilon}{4} \fr{\phi^\beta
 (x \cdot e_{j}-c_ft_{k}+\tau_\vp)}{\phi^\beta  (x \cdot e_{j}-c_ft_{k}
+\tau_j^* )}e^{
-\mu (x \cdot e_{j}-L )}\\
&\g -\fr{\left|(\tau_{j}^*
+L)\phi'
 (x \cdot e_{j}-c_ft_{k}+\tau_{j}^*-\theta(\tau_{j}^*
+L) )\right|}{\phi^\beta  (x \cdot e_{j}-c_ft_{k}
+\tau_j^* )}-\frac{\vp}{4}-\frac{\vp}{2}\\
&\g  -\vp
\end{align*}
for all $x \in
\overline{\mathcal{H}_{j}}$ with $L \leq x \cdot e_{j}
\leq K_{2,\vp}+c_ft_{k} / 2$ and for $k$ large enough.
Moreover, one has
\begin{align*}
 \fr{u(t_{k}, x)
- \phi(x \cdot e_{j}
-c_ft_{k}+\tau_{j}^*)}{\phi^\beta (x \cdot e_{j}-c_ft_{k}
+\tau_j^*)} \l  \fr{1
- \phi(x \cdot e_{j}
-c_ft_{k}+\tau_{j}^*)}{\phi^\beta (x \cdot e_{j}-c_ft_{k}
+\tau_j^*)} \l \vp
\end{align*}
for all $x \in
\overline{\mathcal{H}_{j}}$ with $L \leq x \cdot e_{j}
\leq K_{2,\vp}+c_ft_{k} / 2$ and for $k$ large enough since $\phi(-\i)=1$. Therefore, we have proved \eqref{eqth1-4}, hence \eqref{u-p} is true.

Since $ \phi(-\infty)= 1$ and
$ u(t, \cdot) \rightarrow  1$
as $t \rightarrow+\infty$ locally uniformly
in $\bar{\Omega}$, for $k$ large
enough, there holds  $({1-\phi(L-c_ft_{k}+\tau_{j}^*))}/{\phi^\beta(L-c_ft_{k}+\tau_{j}^*)} \l \varepsilon $ and ${(1-u(t, x))/}{\phi^\beta(L-c_f t_{k}+\rho\hat\vp+\tau_{j}^*)} \l
\varepsilon$
for all $t \geq t_{k}$ and $x \in \overline{\mathcal{H}_{j}}$
with $x \cdot e_{j}=L$. Together with \eqref{u-p}, one obtains from
Lemma \ref{lem2.3} that, for $k$ large enough,
\begin{align}\label{+il}
\fr{| u(t, x)- \phi(x
\cdot e_{j}-c_ft+\tau_{j}^*)|}{\phi^\beta(x
\cdot e_{j}-c_ft+\tau_{j}^*)} \leq \bar Me^{{c_f\beta\rho\vp/\delta}}\vp \ \ \text{for all $t \geq t_{k}$ and $x \in
\overline{\mathcal{H}_{j}}$ with $x \cdot e_{j} \geq L$,}
\end{align}
where $\rho>0$ is defined by \eqref{defrho1} (note that $\rho$ is independent of $\vp$).

{\it Step 3: proof of \eqref{dm3}.}
Since $\vp$ is arbitrary, it can be inferred from \eqref{+il} that
\begin{align*}
 \fr{  u(t, x)- \phi(x
\cdot e_{j}-c_ft+\tau_{j}^*) }{\phi^\beta(x
\cdot e_{j}-c_ft+\tau_{j}^*)}
{\rightarrow}   0\ \ \text{as $t \rightarrow+\infty$ uniformly for $x \in \overline{\mathcal{H}_{j}}$ such
that $x \cdot e_{j} \geq L$.}
\end{align*}
In addition, it also can be derived from $ \phi(-\infty)= 1$ and the complete propagation of  $u$ that
\begin{align*}
 \fr{  u(t, x)- \phi(x
\cdot e_{j}-c_ft+\tau_{j}^*) }{\phi^\beta(x
\cdot e_{j}-c_ft+\tau_{j}^*)}
{\rightarrow}   0\ \ \text{as $t \rightarrow+\infty$ uniformly for $x \in \overline{\mathcal{H}_{j}}$ such
that $0\leq x \cdot e_{j} \leq L$.}
\end{align*}
Since $j \in J$ was arbitrary, then
$(u(t, x)- \phi(x \cdot e_{j}-c_ft+\tau_{j}^*))/\phi^\beta(x \cdot e_{j}-c_f  t+\tau_{j}^*)
\rightarrow  0$ uniformly in $\overline{
\mathcal{H}_{j}} \cap \bar{\Omega}$ as
$t \rightarrow+\infty$ for any $j\in J$.
By $ \phi(-\infty)= 1$ and
$ u(t, \cdot) \rightarrow  1$
as $t \rightarrow+\infty$ locally uniformly
in $\bar{\Omega}$, it can be inferred from $x\cdot e_i>0$ for all $x\in\mathcal{H}_{i}$ that $ (1-u(t, x))/\phi^\beta\left(-x \cdot e_{i}-c_f  t\right) \rightarrow  0$ as
$t \rightarrow+\infty$ uniformly in  $\overline{\Omega
\backslash \bigcup_{s=1}^m \mathcal{H}_{s}}$ for all $i\in I$.
Fix any $i\in I$. If there exists $(t_0,x_0)\in\R\times (\bar{\mathcal H_i}\cap\bar\Omega)$ such that $u(t_0,x_0)=1$, then the maximum principle and the uniqueness of the solution of initial-boundary problem imply that $u\equiv1$ in $\mathbb R\times\bar\Omega $. Thus, it follows from $\phi(-\i)=1$  that  $ (1-u(t, \cdot))/\phi^\beta\left(-x \cdot e_{i}-c_f  t\right)=  0$ as
$t \rightarrow+\infty$ uniformly in $\bar{\mathcal H_i}\cap\overline{\Omega}$. Now we consider the case $u(t,x)<1$ for all $t\in\R$ and $x\in\bar{\mathcal H_i}\cap\bar\Omega$. Since $ u$ is
increasing in time $t$ (proved in Theorem \ref{existence-})
and satisfies \eqref{dm0}, then for any small $\vp>0$, there exists $t^-<0<t^+$ such that $u(t^+,x)\geq u(t^-,x)+\vp$ and $u(t^-,x)\geq\phi(-x\cdot e_i-c_ft^-)-\vp$ for all $x\in\bar{\mathcal H_i}\cap\bar\Omega$. Hence, there holds $u(t^+,x)\geq\phi(-x\cdot e_i-c_ft^-)$ for all $x\in\bar{\mathcal H_i}\cap\bar\Omega$. Since $\phi(-\i)=1$, it then follows from the comparison principle that $u(t,x)\to1$ uniformly for $x\in\bar{\mathcal H_i}\cap\bar\Omega$ as $t\to+\i$. Hence, $ (1-u(t, x))/\phi^\beta\left(-x \cdot e_{i}-c_f  t\right) \rightarrow  0$ as
$t \rightarrow+\infty$ uniformly in $\bar{\mathcal H_i}\cap\overline{\Omega}$. As a consequence, one gets that \eqref{dm3} holds. Finally, it is easy to check that $ u$ is a transition
front connecting $0$ and $1$ in the sense of Definition \ref{td1.1} with sets
$\Gamma_{t}$ and $\Omega_{t}^{\pm}$ defined by \eqref{dm1}
and \eqref{dm2}. The proof of Theorem \ref{existence} is thereby complete.
\end{pr}
\begin{rem}\lb{lem2.5}
By the same arguments as the proof of \cite[Lemma 3.6]{ghs}, we obtain that for every $0<a \leq  b<1$, there holds
\begin{align*}
\inf_{(t, x) \in \mathbb{R} \times \bar{\Omega}, \  a
\leq  u(t, x) \leq   b }  u_{t}(t, x)>  0,
\end{align*}
where $I$, $J$ and $u$ are defined as in Theorem \ref{existence}.
\end{rem}

\SE{Uniqueness of the global mean speed}\label{s4}
In this section, we prove the existence
and uniqueness of  the global mean speed of
transition front connecting $0$ and $1$ in cylindrical domains with multiple branches.
Recall that $L>0$ is given in \eqref{defmulti}. We take $\delta>0$  defined by \eqref{defdelta} in the sequel.

\begin{lem}\lb{lem3.1} For any $\varepsilon \in\left(0,
c_f \right)$, there exist some real numbers $\delta_\vp>0$,
$L_{\varepsilon}>L$ and $R_{\varepsilon}>0$,
such that for any $i \in\{1, \cdots, m\}$ and
any $l \geq R_{\varepsilon}+L_{\varepsilon}$, the solution $v_{i, l, R_{\varepsilon}}(t, x)$ of  Cauchy problem
\be\lb{eq2-1}
\begin{cases}\left(  v_{i, l, R_\vp}\right)_{t}(t,x)
-  \Delta   v_{i, l, R_\vp}(t,x)=f\left(  v_{i, l, R_\vp}(t,x)\right), &
t>0,\  x \in \bar{\Omega}, \\ \left(
v_{i, l, R_\vp}\right)_{\nu}(t,x)=  0, &
t>0,\  x \in \partial \Omega\end{cases}
\ee
with initial condition
\be\lb{eq2-2}
  v_{i, l, R_\vp}(0, x)= \begin{cases}
  1-\delta_\vp   & \text { for }
x \in \overline{\mathcal{H}_{i}} \text { with } l-R<x
\cdot e_{i}<l+R, \\   0 &
\text { elsewhere in } \bar{\Omega}
\end{cases}
\ee  satisfies
\be\lb{eq2-3}
\begin{array}{ll}
  v_{i, l, R_{\varepsilon}}(t, x) \geq  1-2
\delta_\vp  \text { for all } 0 \leq t \leq
T_{\varepsilon}=\frac{l-R_{\varepsilon}-L_{\varepsilon}}
{c_f -\varepsilon}
\text { and } x \in \overline{\mathcal{H}_{i}}
\text { with }\left|x \cdot e_{i}-l\right| \leq\left(c_f
-\varepsilon\right) t
\end{array}
\ee
$($notice that these points $x$ belong to
$\overline{\mathcal{H}_{i}} \cap \bar{\Omega}$
since $x \cdot e_{i} \geq -l+\left(c_f -\varepsilon\right)
T_{\varepsilon}=R_{\varepsilon}+L_{\varepsilon}>L$$)$ and
\be\lb{eq2-4}
\begin{array}{ll}
  v_{i, l, R_{\varepsilon}}(t, x) \geq
1-3 \delta_\vp
\text { for all } t \geq T_{\varepsilon}
\text { and } x \in \overline{\mathcal{H}_{i}}
\text { with } R_{\varepsilon}+L_{\varepsilon} \leq
x \cdot e_{i} \leq l+\left(c_f -\varepsilon\right) t.
\end{array}
\ee
\end{lem}
\begin{pr}
The proof is divided into three steps.

\textit{Step 1. choice of some parameters}. Take any $\vp\in(0,c_f)$. By the property of the function $\phi$, there exists a constant $C>0$ such that
\begin{equation}\label{cc}
\begin{cases}
\phi\geq1-\delta \text{ and }\phi''<0\text{ in }(-\infty,-C],\ \phi\leq\delta\text{ and }\phi''>0\text{ in } [C,+\infty),\\
c_f \beta\fr{\phi '}{\phi}+\beta^2\e(\fr{\phi '}{\phi}\r)^2\l \frac3{4}{\Lambda(-\beta c_f)},
\end{cases}
\end{equation}
where $\beta\in(0,1)$ is any constant.
Since $\phi'<0$ in $\R$, there are two real numbers $\kappa>0$ and $\delta_\varepsilon>0$ satisfying
\begin{align}\label{d-v}
\phi'\leq -\kappa\text{ in }[-C,C]
\ \text{ and }\
\delta_\varepsilon=\min\left(\vp,\frac{\varepsilon \kappa}{6L_f},\frac{\varepsilon \kappa}
{3(L_f + \|\phi''/\phi\|_{L^\i(\R)})},\frac{\delta}{2}\right),
\end{align}
where $L_f$ is defined by \eqref{deflf}.
 It is easy to see that there  exists  $C_{\varepsilon}>C$ such that
\begin{equation}\label{c-v}
\phi\geq1-\delta_\varepsilon\text{ in }(-\infty,-C_\vp]\
\text{ and }\
\phi\leq\delta_\varepsilon\text{ in }[C_{\varepsilon},+\infty).
\end{equation}
Since $\phi'<0$ in $\mathbb R$, there exists a constant $\kappa_\vp>0$ such that
\begin{align}\label{k-v}
\phi'\leq-\kappa_\vp\ \text{ in }\ [C,C_\vp].
\end{align}
Take
\begin{align}\label{al-v}
\alpha_\vp=\min\left(\fr{1}{2},\fr{2 \kappa_\vp}{3|f'(1)|},\fr{\varepsilon\kappa}{3|f'(1)|}\right).
\end{align}
Then there also exists a constant $\tilde C_\vp>C_\vp$ such that
\begin{align}\label{t-c-v}
\phi\geq1-\alpha_\vp\delta_\vp\ \text{ in }\ (-\infty,-\tilde C_\vp].
\end{align}
Let $h_\varepsilon:[0,\infty)\to(0,\infty)$ be a $C^2$ function satisfying
\begin{align}\label{h-v}
\begin{cases}
0\leq h_\varepsilon'(r)\leq1\text{ and } h_\varepsilon''(r)\leq\frac{\varepsilon}{2}&\text{in }[0,\infty),\\
h_\varepsilon'=0&\text{in a neighborhood of }0,\\
h_\varepsilon(r)=r&\text{in }[H_\varepsilon,\infty)\text { for some }H_\varepsilon>0.
\end{cases}
\end{align}
Notice in particular that
\begin{align}\label{h-v1}
 r\leq h_\varepsilon(r)\leq r+h_\varepsilon(0)\ \text{ for all }\ r\geq0.
\end{align}
Furthermore, we define
\begin{align}\label{r-l}
R_\vp=H_\vp+h_\vp(0)+C+C_\vp+\tilde C_\vp\ \text{ and }\
L_\vp= L+C_\vp-C(\geq L).
\end{align}
Finally, for any $i\in\{1,\cdots,m\}$ and any $l\geq R_\vp+L_\vp$, define $T_\vp=(l-R_{\varepsilon}-L_{\varepsilon})/
(c_f -\varepsilon)$.

\textit{Step 2: proof of \eqref{eq2-3}}. Notice that \eqref{eq2-3} is immediately if $T_\vp=0$. Let now
\begin{equation*}
\underline v(t,x)=
\begin{cases}
\max(\phi(\underline\zeta(t,x))-\delta_\vp,0),  &t\in[0,T_\vp]\text{ and }x\in \bar{\mathcal H_i}\cap\bar{\Omega},\\
0,&t\in[0,T_\vp]\text{ and }x\in \bar{\Omega\backslash\mathcal H_i},
\end{cases}
\end{equation*}
where
\begin{align*}
\underline\zeta(t,x)=h_\vp(|x\cdot e_i-l|)-(c_f-\vp)t-R_\vp+C.
\end{align*}
We show that the
function $\underline v$ is a subsolution of the problem \eqref{eq2-1} and \eqref{eq2-2} for $0\leq t\leq T_\vp$ and $x\in\Omega$.
 Since $l-L\geq l-L_\vp\geq R_\vp\geq H_\vp$, it then follows from \eqref{h-v1} and \eqref{r-l} that
 for $x\in \bar{\mathcal H_i}\cap\bar{\Omega}$ with $x\cdot e_i\leq L$ and $t\in[0,T_\vp]$,
 there holds $\underline\zeta(t,x)\geq l-L-(c_f-\vp)T_\vp-R_\vp+C=-L+L_\vp+C\geq C_\vp$.
 For these $t$ and $x$, one  infers from \eqref{c-v} that $\phi(\underline\zeta(t,x))\leq\delta_\vp$,
 which implies that $\underline v(t,x)=0$. Together with the definition of $\underline v$, one gets that $\underline v$ is continuous.
 Observe also that $\underline v$ is of class $C^{1,2}$ in the set where it is positive.

We now check the initial and boundary conditions. For $x\in\bar{\mathcal H_i}\cap\overline\Omega$ such that $l-R\leq x\cdot e_i\leq l+R$,
 there holds $\underline v(0,x)\leq \max(1-\delta_\varepsilon,0)=1-\delta_\vp=v_{i,l,R_\vp}(0,x)$.
  For  $x$ belongs to elsewhere in $\bar\Omega$, one has $\underline v(0,x)=0=v_{i,l,R_\vp}(0,x)$.
  On the other hand, there holds $\nu\cdot e_i=0$ on $\partial\mathcal H_i\cap\p\Omega$,
  then $\underline v_\nu(t,x)=0$ for all $t\in [0,T_\vp]$ and $x\in\partial\Omega$.

It suffices to verify that
$
\mathscr L[\underline v](t,x):=\underline v_t(t,x)-\Delta  \underline v(t,x)-f(\underline v(t,x))\leq0
$
for all $t\in[0,T_\vp]$ and $x\in \overline\Omega$ such that $\underline v(t,x)>0$. By direct calculations, there holds
\begin{align*}
\mathscr L[\underline v](t,x)
&=\fr\vp2\phi'(\underline\zeta(t,x))
+\left(1-(h_\vp'(|x\cdot e_i-l|))^2\right)\phi''(\underline\zeta(t,x))\\
&~~~+\left(\fr\vp2-h_\vp''(|x\cdot e_i-l|)\right)\phi'(\underline\zeta(t,x))
 +f(\phi(\underline\zeta(t,x)))-f(\underline v(t,x)).
\end{align*}
If $\underline\zeta(t,x)\geq C$, then  $0\leq\underline v(t,x)\leq\phi(\underline\zeta(t,x))\leq\delta\leq\theta$ by \eqref{cc},
hence $f(\underline v(t,x))=f(\phi(\underline\zeta(t,x)))=0$. It follows from \eqref{h-v1} and \eqref{r-l} that $|x\cdot e_i-l|\geq-h_\vp(0)+R_\vp+(c_f-\vp)t\geq H_\vp$,
then
$(h_\vp'(|x\cdot e_i-l|))^2=1$ by \eqref{h-v}. According to $\phi'<0$ and \eqref{h-v}, one has
\begin{align*}
\mathscr L[\underline v](t,x)
&=\fr\vp2\phi'(\underline\zeta(t,x))
+\left(\fr\vp2-h_\vp''(|x\cdot e_i-l|)\right)\phi'(\underline\zeta(t,x))\leq0.
\end{align*}
If $\underline\zeta(t,x)\leq -C$, then  \eqref{cc} yields that $1-2\vp_0\leq1-\delta-\delta_\vp\leq\underline v(t,x)\leq\phi(\underline\zeta(t,x))\leq1$. By \eqref{vp0} and the mean value theorem, there exists $\tau\in(0,1)$ such that
\begin{align*}
f(\phi(\underline\zeta(t,x)))-f(\underline v(t,x))=f'(\underline v(t,x)+\tau(\phi(\underline\zeta(t,x))-\underline v(t,x))\delta_\vp\leq\fr{3 f'(1)}{4}\delta_\vp.
\end{align*}
Since $\phi'<0$, it can be inferred from \eqref{cc} and \eqref{h-v} that
\begin{align*}
\mathscr L[\underline v](t,x)
&\leq\fr\vp2\phi'(\underline\zeta(t,x))
+\left(\fr\vp2-h_\vp''(|x\cdot e_i-l|)\right)\phi'(\underline\zeta(t,x))+\frac{3 f'(1)}{4}\delta_\vp\leq0.
\end{align*}
If $-C\leq\underline\zeta(t,x)\leq C$, then \eqref{h-v1} and  \eqref{r-l} imply that $|x\cdot e_i-l|\geq-h_\vp(0)+R_\vp+(c_f-\vp)t-2C\geq H_\vp$, hence  $(h_\vp'(|x\cdot e_i-l|))^2=1$ by \eqref{h-v}. Since $0\leq\underline v(t,x)\leq\phi(\underline\zeta(t,x))\leq1$, it can be inferred from \eqref{deflf} and the mean value theorem that there exists $\tau\in(0,1)$ such that
\begin{align*}
f(\phi(\underline\zeta(t,x)))-f(\underline v(t,x))=f'(\underline v(t,x)+\tau(\underline v(t,x)-\phi(\underline\zeta(t,x)))\delta_\vp\leq L_f\delta_\vp.
\end{align*}
By \eqref{d-v} and \eqref{h-v}, there holds
\begin{align*}
\mathscr L[\underline v](t,x)
&\leq-\fr\vp 2\kappa
+L_f\delta_\vp\leq0.
\end{align*}

For all $t\in[0,T_\vp]$ and $x\in\mathcal H_i$ with $|x\cdot e_i-l|\leq(c_f-\vp)t$, it follows from \eqref{h-v1} and \eqref{r-l} that
$$
\underline\zeta(t,x)
\leq h_\vp(0)+|x\cdot e_i-l|-(c_f-\vp)t-R_\vp+C
\leq h_\vp(0)-R_\vp+C\leq-C_\vp.
$$
By the comparison principle and \eqref{c-v} that
$$
v_{i,l,R_\vp}(t,x)\geq \underline v(t,x)\geq 1-2\delta_\vp\text{ for all $t\in[0,T_\vp]$ and $x\in\mathcal H_i$ with $|x\cdot e_i-l|\leq(c_f-\vp)t$}.
$$

{\it Step 3: proof of \eqref{eq2-4}}. For $t\geq T_\vp$, we set
\begin{equation*}
\underline w(t,x)=
\begin{cases}
\max\left(\phi(\xi_1(x))+\phi(\xi_2(t,x))-1-2\delta_\vp,0\right)
&\text{for  }x\in\overline{\mathcal H_i}\cap\bar{\Omega},\\
0 &\text{for }x\in\bar{\Omega}\backslash\bar{\mathcal H_i},
\end{cases}
\end{equation*}
where
$$
\xi_1(x)=-x\cdot e_i+L_\vp+h_\vp(0)+C\ \text{ and }\
\xi_2(t,x)=x\cdot e_i-(c_f-\vp)t-l+h_\vp(0)+C-R_\vp.
$$
For all $t\geq T_\vp$ and $x\in\mathcal H_i$ with $x\cdot e_i\leq L$, it follows from \eqref{r-l} that $\xi_1(x)\geq-L+L_\vp+h_\vp(0)+C\geq C_\vp$,
then \eqref{c-v} yields that
$\phi(\xi_1(x))+\phi_2(\xi_2(t,x))-1-2\delta_\vp
\leq\delta_\vp+1-1-2\delta_\vp\leq0$.
Therefore, $\underline w$ is continuous in $[T_\vp,\infty)\times\bar\Omega$ and is of class $C^{1,2}$ in the set where it is positive.

Let us prove that $\underline w(t,x)$ is a subsolution of the problem \eqref{eq2-1} and \eqref{eq2-2} in $[T_\vp,+\i)\times\Omega$. We now check the initial (at time $t=T_\vp$) and  boundary conditions. Let $\underline v(t,x)$ and $\underline\zeta(t,x)$ be defined as in step 2.
If $x\in\mathcal H_i$ with $x\cdot e_i<l$, then
$\underline\zeta(T_\vp,x)\leq l-x\cdot e_i+h_\vp(0)-l+R_\vp+L_\vp-R_\vp+C=\xi_1(x)$ by \eqref{h-v1}.
Since $\phi'<0$, then $\phi(\xi_1(x))\leq\phi(\underline\zeta(T_\vp,x))$, which yields that
$\underline w(T_\vp,x)\leq\underline v(T_\vp,x)\leq v_{i,l,R_\vp}(T_\vp,x)$
 for all $x\in\mathcal H_i$ with $x\cdot e_i< l$.
If $x\in\mathcal H_i$ with $x\cdot e_i\geq l$, then
$\underline\zeta(T_\vp,x)\leq x\cdot e_i-l+h_\vp(0)-(c_f-\vp)T_\vp-R_\vp+C=\xi_2(T_\vp,x)$ by \eqref{h-v1}.
Similarly, one has
$\underline w(T_\vp,x)\leq\underline v(T_\vp,x)\leq v_{i,l,R_\vp}(T_\vp,x)$
 for all $x\in\mathcal H_i$ such that $x\cdot e_i\geq l$.
Therefore, one infers that
$\underline w(T_\vp,x)\leq v_{i,l,R_\vp}(T_\vp,x)$
for all $x\in\bar\Omega$.
Moreover, since $\mathcal H_i$ is parallel to $e_i$, we have $\underline w_\nu(t,x)=0$ for all $(t,x)\in[T_\vp,+\i)\times\partial\Omega$.

We now turn to verify
$\mathscr L[\underline w](t,x):=\underline w_t(t,x)-\Delta \underline w(t,x)-f(\underline w(t,x))\leq0$
for all $t\geq T_\vp$ and $x\in\bar{\mathcal H_i}\cap\bar\Omega$ such that $\underline w(t,x)>0$. By some  direct calculations, there holds
\begin{align*}
\mathscr L[\underline w](t,x)
&=-(c_f-\vp)\phi'(\xi_2(t,x))-\phi''(\xi_1(x))-\phi''(\xi_2(t,x))-f(\underline w(t,x))\\
&=f(\phi(\xi_1(x)))+f(\phi(\xi_2(t,x)))-f(\underline w(t,x))+c_f\phi'(\xi_1(x))+\vp\phi'(\xi_2(t,x)).
\end{align*}
We first handle the case $x\in\mathcal H_i$ such that $x\cdot e_i<l$. By \eqref{r-l}, one has $\xi_2(t,x)\leq h_\vp(0)+C-R_\vp\leq-\tilde C_\vp$,  then \eqref{al-v} and \eqref{t-c-v} yield that $\phi(\xi_2(t,x))\geq1-\alpha_\vp\dl_\vp\geq1-2\vp_0$. By the mean value theorem and \eqref{vp0}, there exists $\tau\in(0,1)$ such that
$$
f(\phi(\xi_2(t,x)))=
-f'(1+\tau(\phi(\xi_2(t,x))-1))(1-\phi(\xi_2(t,x)))
\leq\frac{3|f'(1)|}{2}\alpha_\vp\delta_\vp.
$$
If $\xi_1(x)\geq C_\vp$,  then  \eqref{c-v} yields that  $\phi(\xi_1(x))\leq\delta_\vp$, hence $\underline w(t,x)=0$. If $C\leq\xi_1(t,x)\leq C_\vp$, one infers from \eqref{cc} that $\underline w(t,x)\leq\phi(\xi_1(t,x))\leq\delta\leq\theta$, hence $f(\underline w(t,x))=f(\phi(\xi_1(t,x)))=0$. Since $\phi'<0$ and $\delta<c_f$, it follows from \eqref{k-v} and \eqref{al-v} that
$$
\mathscr L[\underline w](t,x)
\leq\frac{3|f'(1)|}{2}\alpha_\vp\delta_\vp-c_f\kappa_\vp\leq0.
$$
If $-C\leq\xi_1(x)\leq C$, then $0\leq\underline w(t,x)\leq \phi(\xi_1(x))\leq1$.  Since $\phi(\xi_2(t,x))\geq1-\alpha_\vp\dl_\vp\geq1-\dl_\vp$, it can be derived from \eqref{deflf} and the mean value theorem that there exists $\tau\in(0,1)$ such that
\begin{align*}
f(\phi(\xi_1(x)))-f(\underline w(t,x))
&=f'\left(\underline w(t,x)+\tau\left(\phi(\xi_1(x))-\underline w(t,x)\right)\right)
\left(1+2\delta_\vp-\phi(\xi_2(t,x))\right)\leq3L_f\delta_\vp.
\end{align*}
Since $\phi'<0$, $\vp<c_f$ and $\delta_\vp<1$, it follows from \eqref{d-v} and \eqref{al-v} that
\begin{align*}
\mathscr L[\underline w](t,x)
\leq3 L_f\delta_\vp-c_f\kappa+\frac{3|f'(1)|}{2}\alpha_\vp\delta_\vp
\leq\frac{\vp\kappa}{2}-c_f\kappa+\frac{\vp\kappa}{2}
\leq0.
\end{align*}
If $\xi_1(x)\leq -C$, it then can be deduced from \eqref{defdelta} and \eqref{cc} that $1-2\vp_0\leq1-\delta\leq\phi(\xi_1(x))\leq1$ and $1-2\vp_0\leq1-4\dl\leq1-\delta+1-\alpha_\vp\dl_\vp-1-2\dl_\vp\leq\underline w(t,x)\leq1$. Since $\phi(\xi_2(t,x))\leq1$, one gets from \eqref{vp0} and the mean value theorem that there exists $\tau\in(0,1)$ such that
\begin{align*}
f(\phi(\xi_1(x)))-f(\underline w(t,x))
&=f'\left(\underline w(t,x)+\tau\left(\phi(\xi_1(x))-\underline w(t,x)\right)\right)
\left(1+2\delta_\vp-\phi(\xi_2(t,x))\right)\leq\frac{3}{2}f'(1)\delta_\vp.
\end{align*}
Thus, it follows from $\phi'<0$ and \eqref{al-v} that
\begin{align*}
\mathscr L[\underline w](t,x)
&\leq\frac{3}{2}f'(1)\delta_\vp+\frac{3|f'(1)|}{2}\alpha_\vp\delta_\vp
\leq0.
\end{align*}

We then consider the case $x\in\mathcal H_i$ such that $x\cdot e_i\geq l$. By \eqref{r-l} and $l\geq R_\vp+L_\vp$, there holds $\xi_1(x)\leq-R_\vp-L_\vp+L_\vp+h_\vp(0)+C\leq-\tilde C_\vp$, then \eqref{t-c-v} yields that $\phi(\xi_1(x))\geq1-\alpha_\vp\dl_\vp\geq1-2\vp_0$. By the mean value theorem and \eqref{vp0}, there exists $\tau\in(0,1)$ such that
$$
f(\phi(\xi_1(x)))=
-f'(1+\tau(\phi(\xi_1(x))-1))(1-\phi(\xi_1(x)))
\leq\frac{3|f'(1)|}{2}\alpha_\vp\delta_\vp.
$$
If $\xi_2(t,x)\geq C_\vp$, one gets from \eqref{c-v} that  $\phi(\xi_2(t,x))\leq\delta_\vp$ and then $\underline w(t,x)=0$. If $C\leq\xi_2(t,x)\leq C_\vp$, one infers from \eqref{cc} that $\underline w(t,x)\leq\phi(\xi_2(t,x))\leq\delta\leq\theta$, hence $f(\underline w(t,x))=f(\phi(\xi_2(t,x)))=0$. Since $\phi'<0$ and $\delta_\vp<\vp$, one infers from  \eqref{k-v} and \eqref{al-v} that
$$
\mathscr L[\underline w](t,x)
\leq\frac{3|f'(1)|}{2}\alpha_\vp\delta_\vp-\vp\kappa_\vp\leq0.
$$
If $-C\leq\xi_2(t,x)\leq C$, there holds $0\leq\underline w(t,x)\leq \phi(\xi_2(t,x))\leq1$. Since $\phi(\xi_1(x))\geq1-\alpha_\vp\dl_\vp\geq1-\dl_\vp$, one then derives from  \eqref{deflf} and the mean value theorem, there exists $\tau\in(0,1)$ such that
\begin{align*}
f(\phi(\xi_2(t,x)))-f(\underline w(t,x))
&=f'\left(\underline w(t,x)+\tau\left(\phi(\xi_1(x))-\underline w(t,x)\right)\right)
\left(1+2\delta_\vp-\phi(\xi_1(x))\right)\leq3L_f\delta_\vp.
\end{align*}
Since $\phi'<0$ and $\delta_\vp<\vp$, it follows from  \eqref{d-v} and \eqref{al-v} that
\begin{align*}
\mathscr L[\underline w](t,x)
\leq3 L_f\delta_\vp-\vp\kappa+\frac{3|f'(1)|}{2}\alpha_\vp\delta_\vp
\leq\frac{\vp\kappa}{2}-\vp\kappa+\frac{\vp\kappa}{2}
\leq0.
\end{align*}
If $\xi_2(t,x)\leq -C$, it then can be deduced from \eqref{cc} that $1-2\vp_0\leq1-\delta\leq\phi(\xi_2(t,x))\leq1$ and $1-2\vp_0\leq1-4\dl\leq1-\delta+1-\alpha_\vp\dl_\vp-1-2\dl_\vp\leq\underline w(t,x)\leq1$. Since $\phi(\xi_1(x))\leq1$, one gets from \eqref{vp0} and the mean value theorem that there exists $\tau\in(0,1)$ such that
\begin{align*}
f(\phi(\xi_2(t,x)))-f(\underline w(t,x))
&=f'\left(\underline w(t,x)+\tau\left(\phi(\xi_2(t,x))-\underline w(t,x)\right)\right)
\left(1+2\delta_\vp-\phi(\xi_1(x))\right)\leq\frac{3}{2}f'(1)\delta_\vp.
\end{align*}
Thus, it follows from $\phi'<0$ and \eqref{al-v} that
\begin{align*}
\mathscr L[\underline w](t,x)
&\leq\frac{3}{2}f'(1)\delta_\vp+\frac{3|f'(1)|}{2}\alpha_\vp\delta_\vp
\leq0.
\end{align*}

 For all $t\geq T_\vp$ and $x\in\bar{\mathcal H_i}\cap\bar{\Omega}$ with $R_\vp+L_\vp\leq x\cdot e_i\leq l+(c_f-\vp)t$, one gets from \eqref{r-l} that $$\xi_1(x)\leq-R_\vp-L_\vp+L_\vp+h_\vp(0)+C\leq-\tilde C_\vp$$
 and
$$\xi_2(t,x)\leq l+(c_f-\vp)t-(c_f-\vp)t-l+h_\vp(0)+C-R_\vp\leq-\tilde C_\vp.$$
 It follows from \eqref{t-c-v} that $\phi(\xi_1(x))\geq1-\alpha_\vp\dl_\vp$ and $\phi(\xi_2(t,x))\geq1-\alpha_\vp\dl_\vp$. Together with the comparison principle and \eqref{al-v}, there holds
$$
v_{i,l,R_\vp}(t,x)\geq\underline w(t,x)\geq1-\alpha_\vp\dl_\vp+1-\alpha_\vp\dl_\vp-1-2\delta_\vp\geq1-3\delta_\vp
$$
 {for all $t\geq T_\vp$ and $x\in\bar{\mathcal H_i}\cap\bar{\Omega}$ with $R_\vp+L_\vp\leq x\cdot e_i\leq l+(c_f-\vp)t$.}
The proof is complete.
\end{pr}

\begin{lem}\lb{lem3.2}
For any $\varepsilon \in\left(0, c_f \right)$, there exist some real numbers $\delta_\vp>0$, $L_{\varepsilon}>L$ and $R_{\varepsilon}>0$ such that, for any $i \in\{1, \cdots, m\}$, any $R>R_{\varepsilon}$ and any $l \geq R+L_{\varepsilon}$, the solution $  w_{i, l, R}(t, x)$ of Cauchy problem
\be\lb{defwi}
\left\{\begin{array}{ll}
\left(  w_{i, l, R}\right)_{t}(t,x)-  \Delta   w_{i, l, R}(t,x)= f\left(  w_{i, l, R}(t,x)\right), & t>0, x \in \bar{\Omega}, \\
\left(  w_{i, l, R}\right)_{\nu}(t,x)= 0, & t>0, x \in \partial \Omega
\end{array}\right.
\ee
with initial condition
\be\lb{defiwi}
  w_{i, l, R}(0, x)=
\bc
\delta_\vp &\text{for $x \in \overline{\mathcal{H}_{i}}$ such that $l-R<x \cdot e_{i}<l+R$},\\
  1 &\text{elsewhere in $\bar{\Omega}$}
\ec
\ee
satisfies
$$
\begin{array}{l}
  w_{i, l, R}(t, x) \leq 2 \delta_\vp \  \text { for all } 0 \leq t \leq T_{\varepsilon}=\frac{R-R_{\varepsilon}}{c_f +\varepsilon}
\text { and } x \in \overline{\mathcal{H}_{i}} \text { with }\left|x \cdot e_{i}-l\right| \leq R-R_{\varepsilon}-\left(c_f +\varepsilon\right) t.
\end{array}
$$
\end{lem}
\begin{pr}
Notice that these points $x$ belong to
$\overline{\mathcal{H}_{i}} \cap \bar{\Omega}$
since $x \cdot e_{i} \geq l-R+R_{\varepsilon}+\left(c_f +\varepsilon\right)
T_{\varepsilon}=l>L$.
Take any $\varepsilon \in\left(0, c _f\right)$. Let $C$, $\delta_\vp$, $C_\vp$, $\kappa$, $h_\vp$, $R_\vp$ and $L_\vp$ all be defined as in Lemma \ref{lem3.1}.
Consider now any $i \in\{1, \cdots, m\}$, any $R>R_{\varepsilon}$ and any $l \geq R+L_{\varepsilon} .$ Set $T_{\varepsilon}=\left(R-R_{\varepsilon}\right) /\left(c_f +\varepsilon\right)$. Define the function
$$
\bar{ w}(t, x)=\left\{\begin{array}{ll}
\min \left(\phi(\bar{\zeta}(t, x))+\delta_\vp,  1\right) & \text { for } t
\geq 0 \text { and } x \in \overline{\mathcal{H}_{i}} \cap \bar{\Omega}, \\
1  & \text { for } t \geq 0 \text { and } x \in \bar{\Omega} \backslash \overline{\mathcal{H}_{i}},
\end{array}\right.
$$
where
$$
\bar{\zeta}(t, x)=-h_{\varepsilon}\left(\left|x \cdot e_{i}-l\right|\right)-\left(c_f +\varepsilon\right) t+R-C_\vp.
$$
If $x
\cdot e_{i} \leq L$, then  \eqref{r-l} implies that  $\left|x \cdot e_{i}-l\right|\geq l-L \geq l-L_{\varepsilon} \geq R>$ $R_{\varepsilon} \geq H_{\varepsilon}$, hence $h_{\varepsilon}\left(\left|x \cdot e_{i}-l\right|\right)=l-x \cdot e_{i}$ by \eqref{h-v}. For all $0 \leq t \leq T_{\varepsilon}$ and $x \in
\overline{\mathcal{H}_{i}} \cap \bar{\Omega}$ such that $x\cdot e_{i} \leq L$, there holds
$\bar{\zeta}(t, x) \leq-l+x \cdot e_{i}+R-C_\vp \leq-L_{\varepsilon}+L-C_\vp\leq-C_{\varepsilon}$,
it then follows from \eqref{c-v} that $\phi(\bar{\zeta}(t, x)) \geq  1-\delta_{\varepsilon}$, which implies that $\bar{ w}(t, x)= 1$.
As a result, the function $\bar{ w}$ is  continuous in $[0,T_\vp] \times \bar{\Omega}$. It is easy to see that $\bar{ w}$ is of class $C^{1,2}$ in the set where
is  positive.

We now turn to show that $\bar{ w}$
is a supersolution of the problem \eqref{defwi} and \eqref{defiwi} in $[0,T_{\varepsilon}]\times{\Omega}$.
Let us first check the initial and boundary conditions.
If $x \in \overline{\mathcal{H}_{i}}$
such that  $l-R<x \cdot e_{i}<l+R$,  there clearly holds
$\bar{ w}(0, x)\g \dl_\vp
= w_{i, l, R}(0, x)$. If $x \in \overline{\mathcal{H}_{i}}$ such that $|x \cdot e_{i}-l|
\geq R(\geq R_\vp)$, then \eqref{h-v1} and \eqref{r-l} imply that   $h_{\varepsilon}\left(\left|x \cdot
e_{i}-l\right|\right)=\left|x \cdot e_{i}-l\right|$
and $\bar{\zeta}(t, x) \leq-C_\vp$,
hence $ \phi(\bar{\zeta}(t, x)) \geq   1-\delta_\vp$ by \eqref{c-v}. Thus,  $\bar{ w}(0,
x)=1\g  w_{i, l, R}(0, x)$ for $x \in \overline{\mathcal{H}_{i}}$ such that $|x \cdot e_{i}-l|
\geq R$. As a result, one has that $
\bar{ w}(0,x)\g  w_{i, l, R}(0, x)$
for all $x \in \bar{\Omega}$.
In addition, it can be inferred that $\bar{ w}_\nu(t,x)= 0$ for $t\in[0,T_\vp]$ and $x\in\p\Omega$ due to
$\nu\cdot e_i=0$ on $\partial\mathcal H_i\cap\p\Omega$.

Let us now prove that
$\mathscr{L}[\bar{ w}](t, x):=\bar{w}_t(t,x)-\Delta
\bar{w}(t,x) -f(\bar{ w}(t, x))\g 0$ for $0 \leq t \leq
T_{\varepsilon}$ and $x \in
\overline{\mathcal{H}_{i}} \cap \bar{\Omega}$ such that
$\bar{ w}(t, x)< 1$.
After some computations, one has
\begin{align*}
\mathscr{L}[\bar{ w}](t, x)
&=-\fr\vp2\phi'(\bar \zeta(t,x))+f(\phi(\bar\zeta(t,x)))-f(\bar {w}(t,x))\\
&~~~+\e(1-\e(h_\vp'(|x\cdot e_i-l|)\r)^2\r)\phi''(\bar \zeta(t,x))-\e(\fr\vp2-h_\vp''(|x\cdot e_i-l|\r)\phi'(\bar \zeta(t,x)).
\end{align*}
If $\bar \zeta(t,x)<-C$, then \eqref{h-v1} and \eqref{r-l} yield that $|x\cdot e_i-l|>-h_\vp(0)-(c_f+\vp)T_\vp+R-C_\vp+C\geq H_\vp$, hence
$h_\vp'(|x\cdot e_i-l|)=1$ by \eqref{h-v}. From \eqref{defdelta} and \eqref{cc}, one has
$1-2\vp_0\leq 1-\dl \l\phi(\bar{\zeta}(t, x))\l \bar{ w}(t, x) \l 1$. It then follows from \eqref{vp0} and the mean value theorem
that there is $\tau\in(0,1)$ such that
\begin{align*}
f (\phi(\bar\zeta(t,x)))-f (\bar{ w}(t, x))
=-f'(\bar{ w}(t, x)
+\tau(\phi(\bar\zeta(t,x))-\bar{ w}(t, x)))\delta_\vp
\g  -\fr34f'(1)\delta_\vp.
\end{align*}
By $\phi'<0$ and \eqref{h-v}, one has
\begin{align*}
\mathscr{L}[\bar{ w}](t, x)\geq f(\phi(\bar\zeta(t,x)))-f(\bar {w}(t,x))
-\e(\fr\vp2-h_\vp''(|x\cdot e_i-l|\r)\phi'(\bar \zeta(t,x))\geq-\fr34f'(1)\delta_\vp\geq0.
\end{align*}
If $\bar\zeta(t,x)>C$, then $0<\phi(\bar\zeta(t,x))<\bar w(t,x)<\dl+\delta_\vp<\theta$, hence $f(\phi(\bar\zeta(t,x)))=f(\bar w(t,x))=0$.
Since $\phi'<0$, one infers from \eqref{cc} and \eqref{h-v} that
\begin{align*}
\mathscr{L}[\bar{ w}](t, x)
&\geq\e(1-\e(h_\vp'(\bar \zeta(t,x))\r)^2\r)\phi''(\bar \zeta(t,x))-\e(\fr\vp2-h_\vp''(|x\cdot e_i-l|\r)\phi'(\bar \zeta(t,x))\geq0.
\end{align*}
If $-C\l \bar\zeta(t,x)\l C$, then
 $|x\cdot e_i-l|>-h_\vp(0)-(c_f+\vp)T_\vp+R_\vp-C_\vp-C
\g H_\vp$ owing to \eqref{h-v1} and \eqref{r-l}, hence $h_\vp'(|x\cdot e_i-l|)=1$.
By \eqref{deflf} and the mean value theorem, there is $\tau\in(0,1)$ such that
\begin{align*}
f (\phi(\bar\zeta(t,x)))-f (\bar{ w}(t, x))=-f'(\bar{ w}(t, x)
+\tau(\phi(\bar\zeta(t,x))-\bar{ w}(t, x)))\delta_\vp\geq-L_f\delta_\vp.
\end{align*}
Since $\phi'<0$, one deduces from \eqref{d-v} and \eqref{h-v} that
\begin{align*}
\mathscr{L}[\bar{ w}](t, x)
&\geq\fr{\vp\kappa}{ 2}-L_f\delta_\vp-\e(\fr\vp2-h_\vp''(|x\cdot e_i-l|\r)\phi'(\bar \zeta(t,x))\geq0.
\end{align*}

According to the comparison principle, there holds
$ w_{i, l, R_{\varepsilon}}(t, x) \leq \bar{ w}(t, x)
$ for all $0 \leq t \leq T_{\varepsilon}$ and $x \in
\bar{\Omega}$. For any $0 \leq t \leq
T_{\varepsilon}$ and $x \in \overline{\mathcal{H}_{i}}$
such that $\left|x \cdot e_{i}-l\right| \leq
R-R_{\varepsilon}-\left(c_f +\varepsilon\right) t$,
it follows from \eqref{h-v1} and  \eqref{r-l} that $
\bar{\zeta}(t, x) \geq-
\left|x \cdot e_{i}-l\right|-h_{\varepsilon}(0)
-\left(c_f +\varepsilon\right) t+R-C_\vp \geq
R_{\varepsilon}-h_{\varepsilon}(0)-C_\vp \geq
C_{\varepsilon}$.
By \eqref{d-v}, then $\phi(\bar{\zeta}(t, x)) \leq
\delta_{\varepsilon}$, which implies
$$
 w_{i, l, R_{\varepsilon}}(t, x) \leq \bar{ w}(t, x)
\leq \delta_{\varepsilon}  +\delta_{\varepsilon}   \leq 2 \delta_\vp.
$$
This completes the proof.
\end{pr}

\begin{lem}\lb{lem3.3}
For any $\varepsilon \in\left(0, c_f \right)$, there exist $\delta_\vp>0$ and $R_{\varepsilon}>0$ such that, for any $R \geq R_{\varepsilon}+L$, the solution
$\widetilde{ w}_{R}(t, x)$  of  Cauchy problem
\be\lb{wr}
\left\{\begin{array}{ll}
\left(\widetilde{ w}_{R}\right)_{t}(t,x)-  \Delta \widetilde{ w}_{R}(t,x)= f\left(\widetilde{ w}_{R}(t,x)\right), & t>0, x \in \bar{\Omega}, \\
\left(\widetilde{ w}_{R}\right)_{\nu}(t,x)= 0, & t>0, x \in \partial \Omega
\end{array}\right.
\ee
with initial condition
\be\lb{wr0}
\widetilde{ w}_{R}(0, x)=\left\{\begin{array}{ll}
\delta_\vp & \text { for } x \in \bar{\Omega} \cap\left(\overline { B ( 0 , L ) } \cup \bigcup\limits _ { i = 1 } ^ { m } \left\{x \in \overline{\mathcal{H}_{i}}: x\cdot e_i<R\right\}\right), \\
 1 & \text { elsewhere in } \bar{\Omega}
\end{array}\right.
\ee
satisfies
\begin{align*}
\widetilde{ w}_{R}(t, x) \leq  3 \dl_\vp &\ \ \text { for all } 0 \leq t \leq T_{\varepsilon}=\frac{R-R_{\varepsilon}-L}{c_f +\varepsilon} \\
&\ \ \text { and } x \in \bar{\Omega} \cap\left(\overline{B(0, L)} \cup \bigcup_{i=1}^{m}\left\{x \in \overline{\mathcal{H}_{i}}: x \cdot e_{i} \leq R-R_{\varepsilon}-\left(c_f +\varepsilon\right) t\right\}\right) .
\end{align*}
\end{lem}
\begin{pr} {Take any $\beta\in(0,1)$ and  $\varepsilon \in\left(0, c _f\right)$. Let $C$, $\dl_\vp$, $C_\vp$ and $\kappa$  be some positive constants defined as in Lemma \ref{lem3.1}.
 Consider a $C^{2}$ nondecreasing function
$\hat{h}_{\varepsilon}: \mathbb{R} \rightarrow[0,1]$ such
that
\begin{equation}\lb{defhathvp}
\hat{h}_{\varepsilon}=1\text { in }[L+1, +\i)\ \ \text{ and }\ \
\ \hat{h}_{\varepsilon}=0 \text { in }(-\i,L].
\end{equation}
Since $\phi(-\i)=0$, by increasing $C>0$ defined as in \eqref{cc} if necessary, one can assume that
\begin{align}\label{hhhh}
\phi(\xi)\left(\|\hat{h}_{\varepsilon}''\|_{L^\i(\R)}
+2\|\hat{h}_{\varepsilon}''\|_{L^\i(\R)}
\|\phi'/\phi\|_{L^\i(\R)}\right)
\l -\fr14\dl_\vp \Lambda(-\beta c_f)\ \ \text{ for }\xi\geq C.
\end{align}
 Set
\be\lb{defrvp1}
R_{\varepsilon}=2 C+2C_\vp+1.
\ee
Take any $R \geq R_{\varepsilon}+L(\geq L+1)$. Define $T_{\varepsilon}=\left(R-R_{\varepsilon}-L\right) /\left(c_f +\varepsilon\right)$.}

 {For all $t \in\left[0, T_{\varepsilon}\right]$ and $x \in \bar{\Omega}$, let us set
$$
\bar{ w}(t, x)=
\begin{cases}
\min \left(\hat{h}_{\varepsilon}\left(x\cdot e_i\right) \phi\left(\xi_{i}(t, x)\right)+\delta_{\varepsilon} \phi^\beta(\xi(t))
+\dl_\vp,  1\right)\\
~~~~~~~~~~~~~~~~~~~~~~~~~~~~~~~~~~~~~~\text{for  } x \in \overline{\mathcal{H}_{i}}  \text { such that } x \cdot e_{i}>L, \ i=1, \cdots, m, \\
\delta_{\varepsilon} \phi^\beta(\xi(t))
+\dl_\vp
\ \ \ \ ~~~~~~~~~~~~~ \ \text{ for  } x \in \bar\Omega\cap\overline{B(0,L)},
\end{cases}
$$
where
$$
\xi_{i}(t, x)=-x \cdot e_{i}-\left(c_f +\varepsilon\right) t+R-C-C_\vp
\ \text{ and }\
\xi(t)=-L-1-\left(c_f +\varepsilon\right) t+R-C-C_\vp.
$$
We shall show that $\bar{ w}$ is a supersolution of the
problem \eqref{wr} and \eqref{wr0} in $\left[0,
T_{\varepsilon}\right] \times \bar{\Omega}$.
Since $2\dl_\vp\leq1$ and $\hat h_\vp(x\cdot e_i)=0$  for $x \in \overline{\mathcal{H}_{i}}
\cap \bar{\Omega}$ such that $x \cdot e_{i} \leq L$  with any $i \in\{1, \cdots, m\}$  by \eqref{defhathvp},
hence the function $\bar{w}$ is well defined
in $\left[0, T_{\varepsilon}\right] \times \bar{\Omega}$.}

 Let us now verify the initial and boundary conditions. If $x \in
\overline{\mathcal{H}_{i}} \cap \overline{\Omega}$ such that
$x \cdot e_{i}<R$ for some  $i \in\{1, \cdots, m\}$, or if  $x \in \bar{\Omega} \cap \overline{B(0, L)}$, we get from $\hat h_\vp\geq0$ and $\phi>0$ that
$\bar{ w}(0, x)\geq\delta_\vp \geq \widetilde{ w}_{R}(0, x)$.
If $x \in \overline{\mathcal{H}_{i}}$
and $x \cdot e_{i} \geq R(\geq L+1)$ for some $i \in\{1,
\cdots, m\}$, then  $\xi_{i}(0, x)=-x\cdot e_i+R-C-C_\vp
\leq-C_\vp$, hence one obtains from $\phi>0$, \eqref{c-v} and  \eqref{defhathvp}  that
 \[\bar{ w}(0, x) \geq \min
\left( 1-\delta_\vp +\dl_\vp,  1\right)= 1
\geq \widetilde{ w}_{R}(0, x).
\]
As a result, one obtains that $\bar{ w}(0,
\cdot) \geq \widetilde{ w}_{R}(0, \cdot)$ in $\bar{\Omega}$.
In addition, since each $\mathcal{H}_{i}$ is parallel to $e_{i}$, there holds $\bar{
w}_{\nu}(t, x)= 0$ for all $(t, x) \in\left[0,
T_{\varepsilon}\right] \times \partial \Omega$.

We then check that
$$
\mathscr{L} [ \bar{ w}](t, x):=\bar{w} _{t}(t, x)- \Delta \bar{w} (t, x)-f (\bar{ w}(t, x)) \geq 0
$$
for  $0 \leq t \leq T_{\varepsilon}$ and $x \in
\bar{\Omega}$ such that $\bar{ w}(t, x)<1 $.
For $x \in \overline{B(0, L)}$, it follows from $\phi<1$, \eqref{defdelta} and \eqref{d-v} that
$0\leq\bar{ w}(t,x)\le2\delta_\vp\leq\theta$, which implies that $f(\overline w(t,x))=0$. After some direct computations, one gets from $\phi'<0$ that
\begin{align*}
\mathscr{L}[ \bar{ w}](t, x)
&=
-\e(c_f+\vp\r) \dl_\vp\beta\phi^{\beta-1}(\xi(t))\phi'(\xi(t))-f(\overline w(t,x))
\g 0.
\end{align*}

 For $x \in \overline{\mathcal{H}_{i}}$ such that $x \cdot e_{i}\geq L+1$  with some $i \in\{1, \cdots, m\}$,
one infers from \eqref{defhathvp} that $\hat h_\vp(x\cdot e_i)=1$,
 hence $ \bar w(t,x)=\phi\left(\xi_{i}(t, x)\right)+\delta_{\varepsilon} \phi^\beta(\xi(t))
+\dl_\vp$.
After some calculations, one has
\begin{align*}
\mathscr{L}[ \bar{ w}](t, x)
&=f\left(\phi\left(\xi_{i}(t, x)\right)\right)-f(\bar{w}(t, x))-\vp\phi'\left(\xi_{i}(t, x)\right)
   -\delta_\vp\beta\left(c_f+\vp\right)\phi^{\beta-1}\left(\xi(t)\right)\phi'\left(\xi(t)\right).
\end{align*}
If $\xi_{i}(t, x)<-C$, then \eqref{defdelta} and \eqref{cc} lead to
$1>\bar{ w}(t, x)\geq
\phi\left(\xi_{i}(t, x)\right) \geq  1-\delta
\g  1-\vp_0$.
Together with \eqref{vp0} and $\phi'<0$, one has
\begin{align*}
 \mathscr{L}[ \bar{ w}](t, x)\geq
 f \left(\phi\left(\xi_{i}(t,
x)\right)\right)-f (\bar{ w}(t, x))
\geq-\fr34f'(1) \e(\delta_{\varepsilon}\phi^\beta\left(\xi(t)\right)
 +\dl_\vp \r)
\g 0.
\end{align*}
If $-C \leq \xi_{i}(t, x)<C$,  then
$
0\leq\bar{ w}(t, x)\leq1$ by $ \phi>0$.
From \eqref{deflf} and $\phi<1$, one obtains that
\begin{align*}
 f\left(\phi\left(\xi_{i}(t,
x)\right)\right)-f(\bar{w}(t, x))
 \geq -L_f\e(\delta_{\varepsilon} \phi^\beta(\xi(t))+\dl_\vp\r)
 \geq-2L_f\delta_\vp.
\end{align*}
By a direct calculation, it can be inferred from $\phi'<0$, $\phi<1$, $\beta<1$ and \eqref{d-v} that
\begin{align*}
\mathscr{L} [ \bar{ w}](t, x)
\geq  -2L_f\delta_\vp+\vp \kappa
\geq  0.
\end{align*}
If $\xi_{i}(t, x) \geq C$, one then gets $f\left(\phi\left(\xi_{i}(t, x)\right)\right)
= f (\bar{ w}(t, x))=0$,
since  \eqref{defdelta}, \eqref{d-v} and \eqref{c-v} imply that
$\phi\left(\xi_{i}(t, x)\right) \leq \bar{ w}(t, x)\leq \delta+3\delta_\vp  \l \theta$.
Combining with $\phi'<0$, one has $\mathscr{L} [ \bar{ w}](t, x)\geq0$.

 For
$x \in \overline{\mathcal{H}_{i}}$ such that $L\leq x \cdot e_{i}\leq L+1$  with some $i \in\{1, \cdots, m\}$,
it is easy to see that $ \bar w(t,x)=\hat h_\vp(x\cdot e_i)\phi\left(\xi_{i}(t, x)\right)+\delta_{\varepsilon} \phi^\beta(\xi(t))
+\dl_\vp$. By \eqref{defrvp1} and the definition of $T_\vp$,
there holds $\xi_{i}(t, x)\geq\xi(t)\geq-L-1-\left(c_f +\varepsilon\right) T_\vp+R-C-C_\vp\geq C$.
From \eqref{defdelta}, \eqref{d-v} and \eqref{c-v}, one gets that
$0\leq\phi\left(\xi_{i}(t, x)\right)\leq\dl\leq\theta$ and $0 \leq\bar{ w}(t, x)\leq \delta+2\delta_\vp  \l \theta$,
hence  $f\left(\phi\left(\xi_{i}(t, x)\right)\right)
= f (\bar{ w}(t, x))=0$.
In addition, it follows from $\phi'<0$ that $\phi(\xi_{i}(t, x))\leq\phi(\xi(t))$.
By $\phi'<0$,
\eqref{defdelta}, \eqref{cc} and \eqref{hhhh}, one obtains that
\begin{align*}
\mathscr{L} [ \bar{ w}](t, x)
=&-(c_f+\vp)\dl_\vp\beta \phi^{\beta-1}(\xi(t))\phi'(\xi(t))
-\vp\hat h_\vp(x\cdot e_i)\phi'(\xi_{i}(t, x))
-\hat{h}_{\varepsilon}'' (x\cdot e_i)\phi(\xi_i(t,x))\\
&+2\hat{h}_{\varepsilon}'(x\cdot e_i)\phi'(\xi_{i}(t, x))
+\hat h_\vp(\xi_{i}(t, x))f(\phi(\xi_{i}(t, x)))-f(\bar w(t,x))\\
\g&\phi^{\beta}(\xi(t))\left(-\dl_\vp c_f\beta\frac{ \phi'(\xi(t))}{\phi(\xi(t))}
-\phi^{1-\beta}(\xi_i(t,x))\left(\|\hat{h}_{\varepsilon}''\|_{L^\i(\R)}
+2\|\hat{h}_{\varepsilon}'\|_{L^\i(\R)}\|\phi'/\phi\|_{L^\i(\R)}\right)\right)\\
\g& \phi^\beta(\xi(t))\e( -\fr34\dl_\vp\Lambda(-\beta c_f)+\fr14\dl_\vp\Lambda(-\beta c_f)\r)
\\
\geq&0.
\end{align*}

According to the comparison principle, one has $\widetilde{ w}_{R} (t,x)\leq \bar{ w}(t,x)$ for  $(t,x)\in\left[0, T_{\varepsilon}\right] \times \bar{\Omega}$.
If $x \in \bar\Omega\cap\overline{B(0, L)}$, it follows from $\delta_\vp<1$ and $\phi<1$   that
\[
\widetilde{
w}_{R}(t, x) \leq \bar{ w}(t, x)\leq\dl_\vp\phi^\beta(\xi(t))+\dl_\vp
 \leq \delta_\vp
+\delta_{\varepsilon}
  \leq 2\delta_\vp.
\]
If $x \in
\overline{\mathcal{H}_{i}}$ with $L<x \cdot e_{i} \leq
R-R_{\varepsilon}-\left(c_f +\varepsilon\right) t$ for some $i\in\{1,\cdots,m\}$, one infers from \eqref{defrvp1} that $\xi_{i}(t, x)
\geq-R+R_{\varepsilon}-(c_f+\vp)T_\vp+R-C-C_\vp\geq C_\vp$. Combining with  \eqref{defhathvp}, $\hat h_\vp\leq1$ and $\phi<1$, there holds
\[\widetilde{ w}_{R}(t, x) \leq \bar{ w}(t, x)
\leq \phi(\xi_i(t,x))+\dl_\vp\phi^\beta(\xi(t))+\dl_\vp\leq
\dl_\vp+\delta_\vp +\delta_{\varepsilon} \leq
3 \delta_\vp\]
for  $x \in
\overline{\mathcal{H}_{i}}$ such that $L<x \cdot e_{i} \leq
R-R_{\varepsilon}-\left(c_f +\varepsilon\right) t$ with $i\in\{1,\cdots,m\}$.
The proof is  complete.
\end{pr}

\begin{rem}\lb{rem1}{\rm It follows from the proofs of Lemmas \ref{lem3.1}-\ref{lem3.3} that one can choose the same $L_{\varepsilon}>L$ and $R_{\varepsilon}>0$ such that all three conclusions hold.}
\end{rem}

\begin{lem}\lb{lem3.4} Let $u$ be any  transition front of \eqref{te1.1} connecting $0$ and $1$ in the sense of Definition \ref{td1.1}. Under the same assumptions of Theorem \ref{gms},  the propagation of $u$ is complete in the sense of \eqref{coms} and, for every $\varrho \geq 0$, there exist some real numbers $T_{1}<T_{2}$ such that
\be\lb{defvarrho}
\bc
\Omega \cap B(0, L+\varrho) \subset \Omega_{t}^{-}\ \ \text{ for all $t \leq T_{1}$}, \\
\Omega \cap B(0, L+\varrho) \subset \Omega_{t}^{+} \ \ \ \text{for all $t \geq T_{2} $.}
\ec
\ee
\end{lem}
\begin{pr} We divide the proof into three steps.

{\it Step 1: complete propagation of $u $.}  Take any $0<\vp<\min(\vp_0,1/2)$, where $\vp_0\in(0,\theta/2)$ is defined as in \eqref{vp0}. Since the time increasing solution $u^i$ of \eqref{defhi} propagates completely, it then follows from Theorem \ref{existence} that $u^i$ satisfies \eqref{dm3}  with $I=\{i\}$ and $j\in\{1,\cdots,m\}\backslash\{ i\}$. By Remark \ref{lem2.5}, one has
\begin{align}\label{kkkv}
\kappa_\vp=\inf_{(t,x)\in\R\times\bar\Omega,\vp\leq u^i(t,x)\leq 1-\vp }u^i_t(t,x)>0.
\end{align}
 Choose any
\begin{align}\label{ogg}
0<\hat \vp<\min\e(\frac{\vp\kappa_\vp}{L_f},\vp\r),
\end{align}
 where $L_f>0$ is defined by \eqref{deflf}.  By \eqref{eq1.4},
every $\Omega_{t}^{+}$ must  contain a half-infinite
branch, that is, for every $t \in \mathbb{R}$, there exist
$\tilde R_{t}>L$ and $i_{t} \in\{1, \cdots, m\}$ such that
\begin{align*}
\left\{x \in \mathcal{H}_{i_{t}}: x \cdot e_{i_{t}}
\geq \tilde R_{t}\right\} \subset \Omega_{t}^{+}.
\end{align*}
 In particular, by denoting $i=i_{0}$,
Definition \ref{td1.1} then implies that
$u (0, x) \geq  1-\hat\vp $
for all $x \in
\overline{\mathcal{H}_{i}}$ with $x\cdot e_i\geq\tilde R_0+M_{\hat\vp}$. Since $u^i(t,x)$ satisfies \eqref{defhi} and $\phi(+\i)=0$, there exists $T_\vp<0$ small enough such that $u^i(T_\vp,x)<\hat\vp$ for all
$x \in
\overline{\mathcal{H}_{i}}$ with $x\cdot e_i<\tilde R_0+M_{\hat\vp}$ and for all $x\in\bar{\Omega\backslash\mathcal H_i}$.

 For all $t\geq0$ and $x\in\bar\Omega$, we define
\begin{align*}
\underline{u }(t, x) = \max \left(u ^{i}\left((1-\vp) t +T_\vp, x\right)  -\hat\vp,  0\right).
\end{align*}
We shall show that the function $\underline{u }(t, x)$ is  a subsolution of the problem satisfied by the transition front $u$ in $[0,+\i] \times
\bar{\Omega}$.  It is easy to see that $u (0, \cdot)
\geq \underline{u }(0, \cdot)$ in $\bar{\Omega}$ and
$\underline{u }_{\nu}(t, x)=0$ for  $t\g0$ and $x
\in \partial \Omega$. Let us then prove that
$$\mathscr{L} [ \underline{u }](t, x):= \underline{u} _{t}(t,
x)- \Delta \underline{u} (t, x)-f (\underline{u }(t, x)) \l 0
$$
 for $t\geq0$ and $x \in
 \bar{\Omega}$ such that
$\underline{u } (t, x)>0$.
By a direct calculation, one has
\begin{align*}
\mathscr{L} [ \underline{u }](t, x)
&=-\vp u^{i}_{t}\left((1-\vp)t+T_\vp, x\right)
   +f \left(u ^{i}\left((1-\vp)t +T_\vp,
x\right)\right)-f (\underline{u }(t, x)).
\end{align*}
If $u ^{i}\left((1-\vp) t +T_\vp,
x\right)>1-\vp$, then
$1\geq u ^{i}\left((1-\vp)t +T_\vp,
x\right)>\underline{u}(t, x)\geq1-2\vp\g1-2\vp_0$. By \eqref{vp0} and the mean value theorem,
  there is $\tau \in(0,1)$ such that
\begin{align*}
f \left(u ^{i}\left((1-\vp) t +T_\vp,x\right)\right)-f (\underline{u}(t, x))
=\hat\vp f'\left(\underline{u }(t, x)+\tau\hat\vp\right)\l \fr34f'(1)  \hat\vp.
\end{align*}
Since $u ^i_t>0$,
then
\begin{align*}
\mathscr{L} [ \underline{u }](t, x)\l
-\vp u^{i}_{t}\left((1-\vp) t +T_\vp, x\right)
+\fr34f'(1)  \hat\vp
\l 0.
\end{align*}
If $u ^{i}\left((1-\vp) t +T_\vp,
x\right)<\vp$, then $\vp<\vp_0$ yields that
$\underline{u }(t,
x)\l u ^{i}\left((1-\vp)t+T_\vp,
x\right) \leq \vp_0\leq\theta$, it implies that
 $f \left(u ^{i}\left((1-\vp)t +T_\vp,
x\right)\right)=f (\underline{u }(t, x))=0$.
By $u ^i_t>0$, there holds
\begin{align*}
\mathscr{L} [ \underline{u }](t, x)
=-\vp u^{i}_{t}\left((1-\vp) t +T_\vp, x\right)
\l 0.
\end{align*}
If $\vp\leq u ^{i}\left((1-\vp) t +T_\vp,
x\right)\leq1-\vp$, it then follows from
 \eqref{deflf} and the mean value theorem that there exists $\tau\in(0,1)$ such that
\begin{align*}
f \left(u ^{i}\left((1-\vp) t +T_\vp,x\right)\right)-f (\underline{u}(t, x))
=\hat\vp f'\left(\underline{u }(t, x)+\tau\hat\vp\right)\l \hat\vp L_f .
\end{align*}
By \eqref{ogg}, one has
\begin{align*}
\mathscr{L} [ \underline{u }](t, x)\l
-\vp \kappa_\vp
+ \hat\vp L_f
\l 0.
\end{align*}

Thanks to the
comparison principle, one has
\begin{align*}
u (t, x) \geq \underline{u }(t, x)\geq u ^{i}\left( (1-\vp)t +T_\vp, x\right)-\hat\vp
\text{ for all $t\g0$ and
$x \in \bar{\Omega}$.}
\end{align*}
Since $u^i$ propagates completely, it then follows from \eqref{ogg} that
\begin{align*}
\lim\limits_{t\to+\infty}u (t, x) \geq 1-\hat\vp\geq1-\vp\ \  \text{ locally uniformly for }x\in\bar\Omega.
\end{align*}
By the arbitrariness of $\vp$,   one gets that $u (t, x) \rightarrow  1$ as
$t \rightarrow+\infty$ locally uniformly in $x \in
\bar{\Omega}$, that is, $u $ propagates completely in the sense of \eqref{coms}.

{\it Step 2: proof of the first formula of \eqref{defvarrho}.}
It suffices to show the assertion for any $\varrho
\geq 0$ large enough, then it will hold automatically
for any $\varrho \geq 0$. Consider any $0<\varepsilon
<\min(c_f,\vp_0/2)$, where $\vp_0$ is defined as in \eqref{vp0}. Let $\dl$ be given in \eqref{defdelta} and let $\delta_\vp$, $L_{\varepsilon}$
and $R_{\varepsilon}$ be defined as in Lemma \ref{lem3.1}.  Even if it means decreasing $\delta_\vp$, we can assume that $\delta_\vp<{\kappa_\vp^2}/{(9L_f^2)}$, where $L_f$ and $\kappa_\vp$ are defined as in \eqref{deflf} and \eqref{kkkv}, respectively. Let
$\varrho$ be any large positive number such that
\be\lb{eq1.2+1}
\bigcup_{i=1}^{m}\left\{y \in \overline{\mathcal{H}_{i}}:
L_{\varepsilon}+R_{\varepsilon} \leq y \cdot e_{i} \leq
L_{\varepsilon}+3 R_{\varepsilon}+M_{\delta_\vp}\right\} \subset
\bar{\Omega} \cap B(0, L+\varrho).
\ee
Consider such $\varrho$ and assume  that the first
statement of  \eqref{defvarrho}
is false for that $\varrho$. It then follows
from \eqref{eq1.3} and up to extraction of a subsequence that two cases may occur:
either there is a sequence
$\left(t_{n}\right)_{n \in \mathbb{N}} $ with $t_n\rightarrow-\infty$ as $n\to+\i$
such that $\Omega \cap B(0, L+\varrho) \cap \Gamma_{t_{n}}
\neq \emptyset$ for each $n \in \mathbb{N}$,
or there is a
sequence $\left(t_{n}\right)_{n \in \mathbb{N}}$ with $t_n\rightarrow-\infty$ as $n\to+\i$ such that $\Omega \cap B(0, L+\varrho)
\subset \Omega_{t_{n}}^{+}$ for each $n \in \mathbb{N}$.

{\it Case 1: $\Omega \cap B(0, L+\varrho) \cap \Gamma_{t_{n}}
\neq \emptyset$ for each $n \in \mathbb{N}$.} For each $n \in
\mathbb{N}$, pick a point $x_{n} \in \Omega \cap B(0,
L+\varrho) \cap \Gamma_{t_{n}}$. It follows from
\eqref{eq1.5} that, for any
$R>0$, there is $r_{R+M_{\dl_\vp}}>0$ such that, for each $n \in \mathbb{N}$,
there exists $y_{n} \in\R^N$ such that
\[y_{n} \in \Omega_{t_{n}}^{+} ,\ \
 d_{\Omega}\left(x_{n}, y_{n}\right) \leq r_{M_{\dl_\vp}+R} \ \ \text{and}\ \
 d_{\Omega}\left(y_{n}, x_{n}\right) \geq
d_{\Omega}\left(y_{n}, \Gamma_{t_{n}}\right) \geq
M_{\delta_\vp}+R .
\] Thus, up to extraction of a subsequence
and by taking $R \geq R_{\varepsilon}$ large enough
independently of $n$, there is $i \in\{1, \cdots, m\}$
such that, for each $n \in \mathbb{N}$,
\[y_{n} \in
\overline{\mathcal{H}_{i}},\ y_{n} \cdot e_{i} \geq
R_{\varepsilon}+L_{\varepsilon}, \ E_{n}:=\left\{y \in
\mathcal{H}_{i}:\left|y \cdot e_{i}-y_{n} \cdot e_{i}
\right| \leq R_{\varepsilon}\right\} \subset
\Omega_{t_{n}}^{+} \ \text{and}\  d_{\Omega}\left(E_{n},
\Gamma_{t_{n}}\right) \geq M_{\delta_\vp}.
\]
Since the number $m$ of branches is
finite, then
$i$ can be chosen independently of $n$ up to extraction
of a subsequence. This implies that
\[u \left(t_{n}, y\right)
\geq  1-\delta_\vp  \ \ \text{ for all
$y \in \overline{\mathcal{H}_{i}}$
such that $\left|y \cdot e_{i}-y_{n} \cdot e_{i}\right|
\leq R_{\varepsilon}$.}
\]
By the comparison
principle and Lemma \ref{lem3.1},
for every $n \in \mathbb{N}$
and $t \geq t_{n}+\frac{y_{n} \cdot
e_{i}-R_{\varepsilon}-L_{\varepsilon}}
{c_f -\varepsilon}$, there holds
$$
u (t, y) \geq  v_{i, y_{n} \cdot e_{i},
R_{\varepsilon}}\left(t-t_{n}, y\right) \geq  1-3
\delta_\vp
$$
for all $y \in \overline{\mathcal{H}_{i}}$ such that
$L_{\varepsilon}+R_{\varepsilon} \leq y \cdot e_{i}
\leq y_{n} \cdot e_{i}+\left(c_f -\varepsilon\right)
\left(t-t_{n}\right)$.
Since
$\left(x_{n}\right)_{n \in \mathbb{N}}$
and $\left(d_\Omega(x_n,y_n)\right)_{n\in\mathbb N}$ are bounded, then $\left(y_{n}\right)_{n \in \mathbb{N}}$ is also bounded. By
letting $n \rightarrow+\infty$ in the
above inequality, one gets
\[u (t, y) \geq  1-3 \delta_\vp \
 \ \text{
for all $t \in \mathbb{R}$ and $y
\in \overline{\mathcal{H}_{i}}$ such that $y \cdot e_{i} \geq
L_{\varepsilon}+R_{\varepsilon}$.}
\]
In particular, $u \left(0,
y\right)\geq  1-3 \delta_\vp $ for all $y \in
\overline{\mathcal{H}_{i}}$ such that $y \cdot e_{i} \geq
L_{\varepsilon}+R_{\varepsilon}$.
 Since $u^i$ satisfies \eqref{defhi} and $\phi(+\i)=0$, there exists $\bar T_\vp <0$ small enough such that $u^i({{\bar T_\vp}},y)<3\delta_\vp$ for all
$y\in
\overline{\mathcal{H}_{i}}$ with $y\cdot e_i<L_{\varepsilon}+R_{\varepsilon}$ and for all $y\in\bar{\Omega\backslash\mathcal H_i}$.
For all $t\geq0$ and $y\in\bar{\Omega}$, define the function
\begin{align*}
\underline{u }(t, y)&=\max \left(u ^{i}((1-\sqrt{\delta_\vp})t+{\bar T_\vp},y)-3\delta_\vp,0\right).
\end{align*}
It is easy to see that $u (0, \cdot)
\geq \underline{u }(0, \cdot)$ in $\bar{\Omega}$ and
$\underline{u }_{\nu}=0$ in  $[0,+\i)\times \partial \Omega$. By a similar argument  in step 1, one can infer that the function $\underline{u }$
is a subsolution of the problem satisfied by $u$ in
$[0,+\i) \times \bar{\Omega}$.
By the comparison
principle, there holds
\begin{align*}
u(t,y)\geq\underline{u}(t, y)\geq u ^{i}((1-\sqrt{\delta_\vp})t+\bar T_\vp ,y)-3\delta_\vp
 \end{align*}
for all $t \geq0$
and $y\in \bar{\Omega}$.
Since $u ^{i}$ propagates completely, one gets that $u (t, y) \geq  1-4\dl_\vp$ for all $t\geq0$ large enough and locally uniformly for
$ y \in \bar{\Omega}$, which is a
contradiction with Definition \ref{td1.1}.

{\it Case 2: $\Omega \cap B(0, L+\varrho) \subset
\Omega_{t_{n}}^{+}$ for each $n \in \mathbb{N}$.}
Owing to the property \eqref{eq1.2+1} satisfied by $\varrho$,
it follows that, for each $n \in \mathbb{N}$ and $i \in\{1,
\cdots, m\}$, there holds $u \left(t_{n}, y\right) \geq
 1-\delta_\vp$ for all
$y \in \overline{\mathcal{H}_{i}}$ such
that $L_{\varepsilon}+R_{\varepsilon} \leq$ $y \cdot e_{i}
\leq L_{\varepsilon}+3 R_{\varepsilon}$. By a similar argument
 as in case 1, one can reach a contradiction. As a result, the proof of the
first part of \eqref{defvarrho}
 has been finished.

{\it Step 3: proof of the second formula of \eqref{defvarrho}.}
Let us prove that, for any $\varrho \geq 0$, there is
$T_{2} \in \mathbb{R}$ such that
\be\lb{varrho+}
\Omega \cap B(0, L+\varrho) \subset \Omega_{t}^{+}\ \
\text {for all } t \geq T_{2}.
\ee
According to the Harnack inequality,  there is a
constant $C>0$ such that $u (t, y) \geq C u (t-1, y)$ for
all $(t, y) \in \mathbb{R} \times \bar{\Omega}$.
Fix now any $\varrho \geq 0$. Assume to the contrary that
\eqref{varrho+} does not hold for any $T_{2} \in \mathbb{R}$.
Then, there exists a sequence  of real numbers  $\left(t_{n}\right)_{n \in
\mathbb{N}}$ such that $t_{n}
\rightarrow+\infty$ as $n \rightarrow+\infty$, and
 a sequence
$\left(x_{n}\right)_{n \in \mathbb{N}}$ in
$\Omega \cap B(0, L+\varrho)$ such that $x_{n}
\notin \Omega_{t_{n}}^{+}$ for every $n \in \mathbb{N}$. By
\eqref{eq1.5} and Definition \ref{td1.1}, there is a sequence $\left(y_{n}\right)_{n
\in \mathbb{N}}$ such that $y_{n} \in \Omega_{t_{n}}^{-}$, $\sup _{n \in \mathbb{N}}
d_{\Omega}\left(y_{n}, x_{n}\right)<+\infty$ together with
 $u \left(t_{n}, y_{n}\right)
\leq C/2   $ for every $n \in
\mathbb{N}$, then $u \left(t_{n}-1, y_{n}\right) \leq
1/2 $ for all $n \in \mathbb{N}$.
On the other hand, since the sequences
$\left(x_{n}\right)_{n \in \mathbb{N}}$ and
$\left(d_{\Omega}\left(y_{n}, x_{n}\right)\right)_{n \in
\mathbb{N}}$ are bounded, it then follows from
the complete propagation of $u$
that $\liminf _{n
\rightarrow+\infty} u \left(t_{n}-1, y_{n}\right) \geq
3/4$, this is a contradiction. The proof is complete.
\end{pr}

 After proving Lemmas \ref{lem3.1}-\ref{lem3.4},  we can prove Theorem \ref{gms} by some similar arguments as in the proof of \cite[Theorem 1.8]{ghs}.
The details of the proof are omitted here.

\SE{Complete   propagation}\label{s6}
In this section, we show some geometrical conditions on $\Omega$ to ensure that the complete propagation assumptions in Theorems \ref{existence} and \ref{gms} are reasonable.
\begin{lem}\lb{lemu1}Assume that,
for every $i \in$ $\{1,\cdots,m\}$, the time increasing
solution $u^{i}$ of \eqref{defhi} propagates completely
in the sense of \eqref{coms}.
Then the entire solution $u(t, x)$ obtained in Theorem \ref{existence-} propagates completely.
\end{lem}
\begin{pr}
Take any $0<\vp<\min(\vp_0,1)$, where  $\vp_0$ is defined as in \eqref{vp0}. Fix any $i \in I$, where the set $I$ is given by Theorem \ref{existence-}.  Since $u^i$ propagates completely, it can be inferred from Theorem \ref{existence} and Remark \ref{lem2.5} that
\begin{align*}
\kappa_\vp=\inf\limits_{(t,x)\in\R\times\bar\Omega,\ \vp\l u^i(t,x)\l1-\vp}u^i_t(t,x)>0.
\end{align*}
Pick any
\begin{align*}
0<\hat\vp<\min\left(\frac{\vp\kappa_\vp}{L_f},\vp\right),
\end{align*}
where $L_f$ is defined as in \eqref{deflf}. By \eqref{defhi} and $\phi^\beta\in(0,1)$, there is $T_\vp<0$ small enough such that
\be\lb{ueq2}
\bc
\left|{u^i(t, x)-\phi\left(-x\cdot e_i-c_ft\right)}\right|
 \leq \fr {\hat\vp}{2} \ \ & \text { for } t \leq T_\vp \text { and } x \in\overline{\mathcal{H}_{i}}\cap\bar\Omega, \\
{u^i(t, x)} \leq \fr{\hat\vp}{2} & \text { for } t \leq T_\vp \text { and } x\in \overline{\Omega\backslash\mathcal{H}_{i}}.
 \ec
\ee
Even if it means decreasing $T_1<0$, one gets from \eqref{dm0} that
\be\lb{ueq3}
\left|u\left(T_\vp, x\right)-\phi\left(-x\cdot e_i-c_f T_\vp \right)\right| \leq\frac{\hat\vp}{2} \  \  \text { for } x \in \bar\Omega\cap\overline{\mathcal{H}_{i}}.
\ee
For all $t\geq0$ and $x\in\bar\Omega$, define
\begin{equation*}
\underline u(t,x)= \max \left(u^i((1-\vp)t+T_\vp, x)-\hat\vp, 0\right).
\end{equation*}
We shall show that $\underline{u}(t, x)$ is a subsolution of the problem satisfied by $u(t, x)$ for $t\geq0$ and $x\in\bar\Omega$.
Let us first verify the initial and boundary conditions. By  \eqref{ueq2} and \eqref{ueq3}, there holds
\begin{align*}
\underline{u}(0, x)
= \max \left(u^i\left(T_\vp, x\right)-\hat\vp , 0\right)
\leq\max\left( \phi\left(-x\cdot e_i-c_fT_\vp\right)-{\hat\vp}/{2},0\right)
\leq u\left(T_\vp, x\right)
\end{align*}
 for $x \in\bar\Omega\cap\overline{\mathcal{H}_{i}}$, and
 \begin{align*}
\underline{u}(0, x)
= \max \left(u^i\left(T_\vp, x\right)-\hat\vp , 0\right)\leq \max \left(\hat\vp/2-\hat\vp , 0\right)=0
 \leq  u\left(T_\vp, x\right)
\end{align*}
 for $x \in \overline{\Omega \backslash \mathcal{H}_i}$. Moreover, it is obvious that $\nu \cdot  \nabla \underline{u}=0$ on $x \in \partial \Omega$.
By a similar argument in  step 1 of the proof of Lemma \ref{lem3.4}, one obtains that
$
\mathscr L[\underline{u}](t, x):=\underline{u}_t(t,x)-\Delta \underline{u}(t,x)-f( \underline{u}(t,x)) \leq 0
$
 for  $t \geq 0$ and $x \in \bar{\Omega}$ such that $\underline{u}(t, x)>0$.

Applying the comparison principle, one concludes that
$$u(t, x) \geq \underline{u}\left(t-T_\vp, x\right) \geq u_i\left((1-\vp)t+\vp T_\vp, x\right)-\hat\vp
$$
for $t \geq T_\vp$ and $x \in \bar{\Omega}$.
Since $u^i(t, x)$ propagates completely, one has that $1-\hat\vp\leq\lim_{t\to+\i}u(t, x) \leq1$ locally uniformly for $x\in\bar{\Omega}$.
By the arbitrariness of $\hat\vp$,
one obtains that $u(t, x)\to1$ as $t\to+\i$ locally uniformly for $x\in\bar\Omega$.
 The proof is complete.
\end{pr}
\vspace{0.2cm}

\begin{pr}[Proof of Theorems \ref{suffc1} and \ref{suffc2}]
We  {need only to} prove that under the condition of Theorems \ref{suffc1} and \ref{suffc2}, for every $i\in\{1,\cdots,m\}$, the propagation of $u^i$ is  complete. Let $g$ be any  $C^{1}([0,1], \mathbb{R})$ function  satisfying the following properties:
\begin{equation*}
\bc
g(0)=g(\vartheta)=g(1)=0, g'(0)<0, g'(\vartheta)>0, g'(1)<0, \\
g<0 \text { in }(0, \vartheta), \ 0<g \leq f \text { in }(\vartheta, 1),\
\int_{0}^{1} g(s) d s>0,
\ec
\end{equation*}
where $\vartheta>\theta$ is a constant with $\theta \in(0,1)$ being the ignition temperature given in \eqref{coma}.
Note that $g\l f$ in $[0,1]$, and that $g$ is a bistable function. For every $i\in \{1,\cdots,m\}$, let $v^i:\R\times\bar\Omega\to(0,1)$ be the time-increasing front-like solution of
\be\label{tg1.1}\bc
v^i_t=\Delta v^i+g(v^i), &t\in\R, x\in\Omega,\\
v^i_\nu=0,  &t\in\R, x\in\p\Omega
\ec
\ee
with the condition
\begin{align}\label{tg1.2}
\bc
v^i(t,x)-\phi(-x\cdot e_i-c_ft)\to0&\text{ uniformly for }x\in\bar{\mathcal H_i}\cap\bar\Omega,\\
v^i(t,x)\to0&\text{ uniformly for }x\in\bar{\Omega\backslash\mathcal H_i}
\ec
\ \ \text{ as }t\to-\i.
\end{align}
Since $g\l f$ and $u^i$ satisfies \eqref{defhi}, then $u^i$ is a subsolution of the problem \eqref{tg1.1} and \eqref{tg1.2}. By providing the conditions of Theorems \ref{suffc1} and \ref{suffc2}, it then follows from the comparison principle that $1\geq u^i(t,x)\geq v^i(t,x)$ for $x\in\bar\Omega$ and $t\in\R$. According to \cite[Corollaries 1.11 and 1.12]{ghs}, one has $v^i(t,x)\to1$ locally uniformly for $x\in\bar\Omega$ as $t\to+\i$ under the conditions of Theorems \ref{suffc1} and \ref{suffc2}, then the propagation of $u^i$ is also  complete. Together with Lemma \ref{lemu1}, one gets that the entire solution obtained in Theorem \ref{existence-} propagates completely.  The proof is thereby complete.
\end{pr}

\section*{Acknowledgments}
The first author's work was partially supported by NSF
of China (11971128) and
by the Heilongjiang Provincial Natural Science
Foundation of China (LH2020A003).
The second author's work was partially
supported by NSF of China (12171120).
The last author's work was partially
supported by NSF of China (12071193, 11731005).

\section*{Date availability statements}
We do not analyze or generate any datasets, because our work proceeds within a theoretical and 
mathematical approach.

\section*{Conflict of interest}
There is no conflict of interest to declare.

\section*{Declaration of generative AI use}
The generative AI was not used in this paper.

\end{document}